\documentclass[12pt, a4paper]{amsart}

\usepackage{amsmath, amsthm, amssymb}
\usepackage{graphicx}
\usepackage{enumerate}
\usepackage[top=2.5cm, bottom=2.5cm, left=2.5cm, right=2.5cm]{geometry}

\newcommand{\x}{\mathbf x}

\newtheorem{proposition}{\bf Proposition}

\renewenvironment{proof}{\noindent {\bf Proof: }}{\rm\\}
\theoremstyle{definition}
\newtheorem{remark}{Remark}{\rm}
{\rm}

\usepackage{color}
\usepackage{stmaryrd}

\usepackage[T1]{fontenc}
\usepackage[export]{adjustbox}

\usepackage{bbm}

\usepackage{mathrsfs}
\usepackage{mathtools}

\newcommand{\vertiii}[1]{{\left\vert\kern-0.25ex\left\vert\kern-0.25ex\left\vert #1 
    \right\vert\kern-0.25ex\right\vert\kern-0.25ex\right\vert}}

    \usepackage{etoolbox}
    
    \usepackage{arydshln}

\DeclarePairedDelimiterX{\normi}[1]
  {|\!|\!|}
  {|\!|\!|}
  {\ifblank{#1}{\:\cdot\:}{#1}}

\usepackage{algorithm,algpseudocode}
\usepackage{caption}
\DeclareCaptionLabelFormat{algnonumber}{Algorithm}
\begin{document}

\title[$H_\infty$ control]{Robust control of  systems with hyperbolic partial differential equations}
\author{Pierre APKARIAN$^1$}
\author{Dominikus NOLL$^{2}$}
\thanks{$^1$ONERA, Department of System Dynamics, Toulouse, France}
\thanks{$^2$Institut de Math\'ematiques, Universit\'e de Toulouse, France}

\begin{abstract}
We discuss strategies to bring $H_\infty$-control techniques into play when the system dynamics are modeled by hyperbolic partial differential equations, or more generally, by systems with non-sectorial pole pattern.
\\

\\

\end{abstract}

\maketitle

\section{Introduction}
Highly oscillatory systems arise commonly in aerospace control and other high technology fields, and their modeling potentially calls for
hyperbolic partial differential equations. For practitioners this rises the question whether such infinite-dimensional oscillatory models are convenient for control, or 
whether they are better advised to continue to rely on more traditional
finite-dimensional approximations based e.g. on identification routines, or 
reduced-order models. 

The answer may depend on whether
$H_\infty$-control strategies, fairly well-established for finite-dimensional systems, may still be brought to 
work for such  infinite-dimensional systems.  
Here we understand $H_\infty$-control as of embedding the system $G$ in a {generalized} plant $P$ and solving a multi-objective 
optimization problem of
the form
\begin{eqnarray}
\label{problem}
\begin{array}{ll}
\mbox{minimize} & \displaystyle\max_{i\in S} \| T_{w_iz_i}(P,K) \|_\infty \\
\mbox{subject to} & \displaystyle\max_{i\in H} \| T_{w_iz_i}(P,K)\|_\infty \leq 1,\\
& K \in \mathscr K \mbox{ stabilizes $G$}
\end{array}
\end{eqnarray}
with $w_i \to z_i$ designer-chosen robustness or performance channels built into $P$, divided into soft $i\in S$ and hard $i\in H$ constraints, and where $\mathscr K$
is a designer-specified class of structured, typically conveniently implemented finite-dimensional controllers; cf. \cite{AN15}. 
The purpose of this note is to demonstrate that while the infinite-dimensional line encounters no principled difficulties for parabolic equations,
hyperbolic systems may still be handled reliably 
under suitable precautions.

The main obstruction to frequency analysis and controller synthesis of hyperbolic systems is that they exhibit an infinity of poles arranged on vertical strips in the complex plane.
When unstable, this prevents use of standard tools like the Nyquist test or the use of the system spectral abscissa \cite{BO:1994}, often used to find 
stabilizing controllers. But even when this string of poles is on the stable
half plane, proximity to the imaginary axis still leads to strange behavior causing difficulties in synthesis. Finite-dimensional approximations of such systems $G$ inevitably miss
high frequency resonant poles with non-negligible magnitude, which is why practitioners may have a point
in considering these $G$ unrealistic. 
If the use of hyperbolic equations for control is to pass muster, it has to demonstrate its ability to provide practical controllers which are robust with regard 
to such highly  oscillatory modes. 
 
When systems of PDEs are considered, the distinction between hyperbolic and parabolic equations is no longer helpful,
and a better way to describe the situation is to distinguish between sectorial and non-sectorial operators, or semi-groups. When the operator is not sectorial, 
then the mentioned difficulties caused by poles with arbitrary high frequencies close to the imaginary axis are felt, while in a sectorial system exceedingly high frequency dynamics 
die out quickly. Here we present a general algorithmic approach capable to deal with such non-sectorial systems, and
then demonstrate its ability by controlling a Timoshenko and an Euler-Bernoulli beam.

The structure of the paper is as follows. We present our algorithm in Section 
\ref{sect_algo} and comment on the individual steps in the subsections, 
highlighting potential difficulties caused by non-sectorial pole pattern.
In the sequel, we use two studies, boundary control of a cantilever Timoshenko beam, and piezo-electric control
of an Euler-Bernoulli beam, to demonstrate the mentioned difficulties 
with these hyperbolic systems.
Section \ref{sect-Timo} briefly presents the cantilever Timoshenko beam model, where we 
prepare three settings, undamped, viscous damping, and Kelvin-Voigt damping. This leads to 
$2 \times 2$ MIMO-control problems. 
Section \ref{sect_Euler} recalls the Euler-Bernoulli beam model, controlled 
by a collocated piezoelectric sensor-actuator pair, again with the options undamped, viscous, and
Kelvin-Voigt damping. 
$H_\infty$-synthesis  for the Timoshenko beam is discussed in Section \ref{sect_timo_control},
and for the Euler-Bernoulli beam in Section \ref{sect_euler_control}, showing how the mentioned difficulties can be overcome. Conclusions
are drawn in Section \ref{sect_conclusion}.

\section{Algorithmic scheme \label{sect-algo}}
\label{sect_algo}
In \cite{AN:18,apkarian:19,an1,peak,ANR:17}, we have developed a general frequency-based approach to $H_\infty$- or $H_2$-control based on (\ref{problem}), 
which can be presented as follows:
 \begin{algorithm}[!ht]
 \begin{algorithmic}[1]
\caption{\!\!{\bf :} $H_\infty$-control of infinite-dimensional systems $G$}
\State{\bf Steady-state.} 
Compute steady state of non-linear system $G_{\rm nl}$, shift it to origin, and obtain linearization $G$.
\State {\bf Transfer function.} Compute transfer function $G(s)$.  
\State{\bf Controller structure.} 
Choose practical controller structure $\mathscr K$ and find initial controller $K_0\in \mathscr K$ stabilizing the loop $(G,K_0)$.
\State{\bf Plant.}
Embed $G$ into plant $P$ with closed-loop performance and robustness specifications. Possibly give
special attention to strong non-linearity in
$G_{\rm nl}$.
\State{\bf Optimize.}
Use non-smooth optimization to solve the multidisk $H_\infty$-optimization program (\ref{problem}), maintaining stability of
the loop $(G,K)$ at iterates $K \in  \mathscr K$. Obtain optimal structured $H_\infty$-controller $K^*\in \mathscr K$.
\State{\bf Simulation.}
Simulate linear closed loop $(G,K^*)$ to verify whether robustness and performances are satisfactory. If not,
modify plant $P$ and specifications, and go back to step 5. 
\State {\bf Non-linear simulation.} Simulate non-linear closed loop $(G_{\rm nl},K^*)$ to verify whether design is satisfactory.
\end{algorithmic}
\end{algorithm}

We would naturally hope that this scheme remains to a large degree general, with only minimal amendments in a given particular case.
This is indeed the case for sectorial systems as for instance seen with parabolic PDEs.
The scheme become more case-dependent when the pole pattern is not sectorial, and in particular,
when
hyperbolic partial differential equations contribute to the dynamics.

In the remainder of this section, we go through the steps of the algorithm,  discuss their implementation, and
comment on the challenges caused by non-sectorial dynamics.

\subsection{Comments on Step 2}
We start by 
noticing that step 2 is not always available analytically. Rather shall we have to compute $G(s)$ using the numerical solution of an 
elliptic boundary value problem for a sufficiently dense set  $s_\nu = j\omega_\nu$ of frequencies.  More formally, if a boundary control system is given under the form
\begin{align}
\label{direct}
G : \left\{
\begin{array}{rl}
\dot{x} \!\!\!&= Ax \\
Px\! \!\!&= u\\
y \!\!\!&= Cx
\end{array}
\right.
\end{align}
with suitable operators $A,P,C$, see \cite[Sect. 3.3]{CurtainZwart1995}, then computation of a single $G(s)$ is obtained by applying the Laplace transform to (\ref{direct})  
and solving the elliptic boundary value problem
\begin{align}
\label{transf}
G(s) : \left\{
\begin{array}{rl}
sx(s)\!\!\! &= Ax(s) \\
Px(s) \!\!\!&= u(s) \\
y(s) \!\!\! &= Cx(s)
\end{array}
\right.
\end{align}
for a large set of $s_\nu=j\omega_\nu$.
This may be time consuming, but can be performed off-line in a pre-computation phase. 
This step will be explained via examples in our experimental section.

The numerical solution of the complex boundary value problem (\ref{transf}) may 
eventually encounter
numerical difficulties for frequencies $\omega$ beyond a certain 
high frequency limit $\overline{\omega}$. When the open-loop system has sufficient roll-off,  
values $G(j\omega)$ beyond $\overline{\omega}$ are usually irrelevant, but for non-sectorial
systems $\bar{\omega}$ may be very large. In that case it is crucial that roll-off
in $L=GK$ 
be generated by the controller, so that Nyquist test and $H_\infty$-norm estimates 
remain reliable. Some knowledge of $\overline{\omega}$ is needed for the theoretical estimates
in \cite{AN:18}. 
Note, however, 
that even in the difficult neutral case the limit $\overline{\omega}$ gives
a much better resolution  than the one we may hope to reach by 
approximations $(E,A,B,C,D)$  based on finite elements or finite differences. 
We may also understand $\overline{\omega}$ as an indicator of up to what 
resolution the infinite-dimensional system may be reliably simulated.
In addition to having to accept such a cut-off frequency $\overline{\omega}$,
there is also the challenge to not miss resonant frequencies 
in the range $[0,\overline{\omega}]$ in a highly oscillatory system.

\subsection{Comments on Step 3}
A difficulty may arise in step 3 of the algorithm. Since the chosen controller structure $\mathscr K$ is motivated by practical considerations like implementability,
simplicity, 
experience with distributed control architectures, it may be hard to
obtain a certified initial stabilizing controller $K\in \mathscr K$ for $G$, 
as this has to be proved for $G$  in infinite dimensions. 

If $G$ has only a finite number $n_p$ of unstable poles,
a finite-dimensional reduced-order approximation for stability $G_r$ of $G$ accurate in the unstable part is typically  available. This is for instance the case
for systems satisfying the spectrum decomposition condition \cite{CurtainZwart1995,engel_nagel}.
Then we may use the following heuristic: 
Compute a stabilizing controller
$K_0 \in \mathscr K$ for $G_r$, and use the Nyquist test together with knowledge of $n_p$ to check whether $K_0$  also stabilizes  $G$. This has good
chances of success in practice, \cite{AN:18,ANR:17,an1,apkarian:19,peak}. Note that the requirement that $G_r$ be accurate in the unstable part may render it 
useful for stabilization, but
better approximations $G_{\rm perf}$  are typically required when performances and robustness have to be addressed.

The true difficulty in step 3 occurs if the open loop
system is not sectorial and has  
infinitely many  unstable poles, typically located in a vertical strip of the right half plane. 
Then no finite-dimensional approximation
$G_r$ of $G$ for control is available, because  any such approximation
has to be exact in the unstable part.  Here traditional
PDE control techniques may offer
ways to find theoretically stabilizing control laws, which may then be used as starting points in optimization. In the case of our studies stability results are for instance
\cite{kim_renardy,liu1998spectrum,tian2017stability,mercier2015non,guiver2010non,mehrvarz2018vibration,morgul1992dynamic,grobbelaar2010uniform}  for the Timoshenko beam, and
\cite{ammari2000stabilization,awada2022optimized,pang1992active,liu2019novel} for the Euler-Bernoulli beam.

Assume for instance that a simple
stabilizing controller $K_0$ for  $G$ or $G_{\rm nl}$ based on a Lyapunov argument is known. Then we proceed as follows: We consider the stable closed loop
$G_0=(G,K_0)$
along with candidate controllers $K$ satisfying $K_0+K\in \mathscr K$. Since we have stability loop equivalence
$(G,K_0+K)\simeq ((G,K_0),K)=(G_0,K)$, we can hope to approximate $G_0=(G,K_0)$ by a stable finite-dimensional system $G_{0,\rm stab}$, and apply the above method
to the structure $\mathscr K'=\mathscr K-K_0$. If $K'\in \mathscr K'$ stabilizing $G_{0,\rm stab}$ is found and certified to stabilize $G_0$
via the Nyquist test, then $K=K'+K_0\in \mathscr K$ is the solution of step 3, as it stabilizes $G$.
For hyperbolic systems the role of the pre-stabilizer $K_0$ is to shift the string of vertical poles to the left half plane, from where on the frequency methods may bear.

\subsection{Comments on Step 5}
In order to perform
optimization in step 5, structured controllers $K\in\mathscr K$ are represented
by a finite set of tunable parameters, which we express by the notation  $K(\x)$ for $\x\in \mathbb R^{k}$.  Due to non-smoothness of the $H_\infty$-norm,
program (\ref{problem}) is then addressed
by a non-smooth optimization technique, supplemented by a method to maintain closed-loop stability
of the iterates $(G,K(\x))$. In finite dimensions the latter may be arranged by
the spectral abscissa $\alpha(A_{\rm cl})$ of the closed-loop system matrix. This takes the form of  a mathematical programming constraint,
$\alpha\left( A_{\rm cl}(\x)\right) \leq -\epsilon$ for some small $\epsilon > 0$, where $A_{\rm cl}(\x)$ is the system matrix
of the loop $(G,K(\x))$.  See \cite{BLO:gradient,BO:1994}.
However,  this 
function is usually not available for infinite-dimensional plants.

In cases where a finite-dimensional approximation $G_{\rm perf}$ for performance can be used, we may follow standard lines, as now implemented 
in the {\tt hinfstruct} \cite{CT2021b,GA2011a} and {\tt systune} \cite{RCT2021b,AGB:2014} MATLAB functions. The optimal $H_\infty$-controller $K^* \in \mathscr K$ will in the end undergo verification
with the true infinite-dimensional model as in \cite{AN:18}.

When no finite-dimensional approximation $G_{\rm perf}$ of $G$  for control of robustness and performances is available, or 
when the available ones hit computationally intractable dimensions, 
the truly infinite-dimensional frequency domain optimization procedure
of \cite{AN:18} is required. Here closed-loop stability of  iterates 
is verified using the Nyquist test as explained in \cite{AN:18}. In addition, a repelling technique 
is required to prevent iterates from repeatedly trying to go outside the hidden domain of stabilizing $K(\x)$.  
The fact that the Nyquist test may be brought to work here is rendered possible by the
preliminary stabilization step, which gives $n_p=0$. When initially $n_p=\infty$, we cannot use the Nyquist test.
Step 5 is critical if a very large number of frequencies $s_\nu = j\omega_\nu$ is required to represent the system.

\subsection{Comment on Lyapunov stability}
When $n_p=\infty$ 
preliminary stabilization requires a Lyapunov function, which in PDE boundary or distributed control often derives from physical knowledge
under the form of an energy functional 
$E(X,K)$, depending on $K\in \mathscr K$ and a string of parameters $X$. 
When the Lyapunov property of  $E$ can  be expressed as a mathematical programming constraint
$g(X,K) \leq 0$, then it can in principle be included in (\ref{problem}). This may give rise to a more stringent stability constraint, e.g. 
when $g(X,K)\leq 0$ implies exponential stability of the loop $(G,K)$. However, this method is limited by the fact that in realistic situations the number of additional variables
gathered in $X$ may grossly exceed the number of truly relevant unknowns $\x$ in $K(\x)$, rendering the numerical method difficult. 
This occurs already for finite-dimensional LTI systems, when
typically $X$ determines a quadratic Lyapunov function $x^\top X x$ and stability is turned into a Bilinear Matrix Inequality (BMI) $g(X,K) \preceq 0$ \cite{GSP:1995}. In that case, 
the number of unknown parameters in $X$ is of the order $(n+n_K)^2$ for $n$ the order of the system and $n_K$ the order of $K$. In practice, this may lead to 
BMIs with several thousands of entries $X$, even when $k$ the number of decision variables $\x$ in $K$ stays way below $100$. Choosing sparse $X$ is not a valid option in practice, as
demonstrated in  \cite{thevenet1,thevenet2}, as this leads to severe conservatism in stability,  performance and robustness specifications.
Furthermore, numerical testing indicates that decision variables $\x$ and Lyapunov variables $X$ 
often differ by several orders of magnitudes, which makes BMI optimization
for control highly ill-conditioned.

\subsection{Comments on Steps 6 and 7}
Simulation of the linearized system is standard when  a finite-dimensional approximation is available. When $G(s)$ is computed formally or
by a succession of boundary value problems, then simulation may be based on the inverse Laplace transform \cite{trefethen2007,cohen2007numerical}. 

For non-linear simulations a state-space approach is inevitable, and here a problem occurs if no finite-dimensional state-space model $G_{\rm perf}$ 
capturing the unstable part of $G$ is available, or if the discrepancy between linear simulations based on $G_{\rm perf}$ and those based on
$G$ is significant. Then it will be hard to decide whether a detected instability or loss of performance  is
caused by non-linearity, or has numerical reasons caused by instability of the open-loop
state-space model. 

\subsection{Comment on stability}
For linear systems 
stability in the $H_\infty$-sense is absence of unstable poles in tandem with boundedness of the closed-loop system on $j\mathbb R$. 
There is much to suggest that this is the most natural form of stability due to its physical relevance. Namely, it means that if we open the loop at two
arbitrary break points $w$ and $z$ and add a $L_2$ source signal at $w$, then we will still receive a $L_2$-signal at $z$. For finite-dimensional systems this 
notion of {\it internal stability} implies
exponential stability, but this need no longer be the case in infinite dimensions.  
However, we have the following quite satisfactory substitute:

\begin{proposition}
\label{morris}
Suppose the closed-loop $(G,K)$ is $H_\infty$-stable, $K$ is finite-dimensional, and $G$ is
exponentially stabilizable and detectable. Then the closed loop is exponentially stable.
\end{proposition}

\begin{proof}
This follows  because $K$ is also exponentially stabilizable and detectable, hence 
by Staffans \cite[Lemma 8.2.7]{Staffans_book} 
so is the loop $(G,K)$, and then by a result of Morris \cite[Theorem 5.2]{Morris:1999}, $(G,K)$ is even exponentially stable. 
\hfill $\square$
\end{proof}

Note that exponential stabilizability and detectability may
be certified with infinite-dimensional controllers which need not be practical. 
The analogue statement concerning strong stability 
is 
\begin{proposition}
\label{staffans}
Suppose the closed loop $(G,K)$
is $H_\infty$-stable, $K$ is finite-dimensional, and $G$ is strongly stabilizable and detectable. Then the closed loop is strongly stable.
\end{proposition}

\begin{proof}
The proof again
uses
\cite[Lemma 8.2.7]{Staffans_book}, which guarantees
that the loop $(G,K)$ is strongly stabilizable and detectable when the components are.
When the loop is $H_\infty$-stable, we infer now using
\cite[Theorem 8.2.11 (ii)]{Staffans_book} that the closed loop is then also strongly stable.
\hfill $\square$
\end{proof}

\section{First study: Cantilever Timoshenko beam \label{sect-Timo}}
Our first study uses a Timoshenko beam model \cite{kim_renardy}, \cite[§ 5.2]{meirovitch} of the form
\begin{align}
\label{timoshenko}
\begin{split}
    \rho \frac{\partial^2w}{\partial t^2} - K \frac{\partial w^2}{\partial x^2} + K \frac{\partial \phi}{\partial x}
    &= f\\
    I_\rho \frac{\partial^2\phi}{\partial t^2} -EI \frac{\partial^2\phi}{\partial x^2} + K \left(\phi-\frac{\partial w}{\partial x}  \right) &= g
\end{split}
\end{align}
where $w(x,t)$ is the total deflection, $\phi(x,t)$ the angle of rotation, $\rho$ the mass per unit length, $I_\rho$ the mass moment of inertia of the cross-section, $E$ Young's modulus, $I$ the moment of inertia of the cross section, $G$  the modulus of elasticity in the shear,
$A$  the cross-sectional area, $k$ a constant depending on the shape of the beam cross-section, and $K=kGA$, all assumed constant. The right hand sides $f(x,t)$,
$g(x,t)$ are used to represent external and internal damping.

The beam is clamped on the left and vibrates freely on the right, with the possibility to use feedback acting at the tip
on bending moment  and shear to stabilize and attenuate disturbances causing vibrations.  This leads to the boundary conditions
\begin{align}
\label{bdry}
\begin{split}
w(0,t)=0, \; \phi(0,t) &= 0 \\
K\frac{\partial w}{\partial x}(L,t) -K \phi(L,t) &= U_1(t) \\
EI \frac{\partial \phi(L,t)}{\partial x} &= U_2(t)
\end{split}
\end{align}
where $L$ is the length of the beam.
Control $U_1$ of shear and $U_2$ of bending use as
measured outputs the speed of rotation and distortion  at the tip position
\begin{align}
y_1(t) = w_t(L,t), \quad y_2(t) = \phi_t(L,t).
\end{align}
For the source terms we consider three scenarios,
the undamped case
\begin{align}
\label{undamped}
f=0, g=0, 
\end{align}
external (or viscous) damping, caused for instance by friction with a surrounding medium like air resistance,
\begin{align}
\label{viscous}
    f(x,t) = -d_w w_t(x,t), \quad g(x,t) = -d_\phi \phi_t(x,t),
\end{align}
and internal, or Kelvin-Voigt, damping
\begin{align}
\label{kv}
\begin{split}
    f(x,t) &= -D_w \left( \phi_{tx}(x,t) - w_{txx}(x,t) \right), \\
    g(x,t) &= -D_w \left( \phi_t(x,t)-w_{tx}(x,t)  \right) + D_\phi \phi_{txx}(x,t).
    \end{split}
\end{align}

In the undamped case
the  system is unstable with infinitely many poles arranged vertically in the right half plane, so that preliminary stabilization is required to start our 
method. The authors of \cite{kim_renardy} show that a simple proportional control law of the form
\begin{align}
\label{prestab}
U_1(t) = -\alpha w_t(L,t) + {u}_1(t), \quad U_2(t) = -\beta \phi_t(L,t) + {u}_2(t)
\end{align}
with $\alpha >0$, $\beta >0$
leads to a stable systems with new controls ${u}_1,{u}_2$. In the terminology
of the algorithm this is the pre-stabilizing controller $K_0$ of step 3,
and the pre-stabilized system is $G_0=(G,K_0)$.
$G_0$ is now amenable to our technique, which means that if a structured control law
$K \in \mathscr K$ is sought, then we can apply the algorithm to the class
$\mathscr K'= \{K': {\rm diag}(-\alpha,-\beta) + K' \in \mathscr K\}$. This is used to optimize performance and robustness of the controller.

\begin{figure}[ht!]
\includegraphics[scale=0.21]{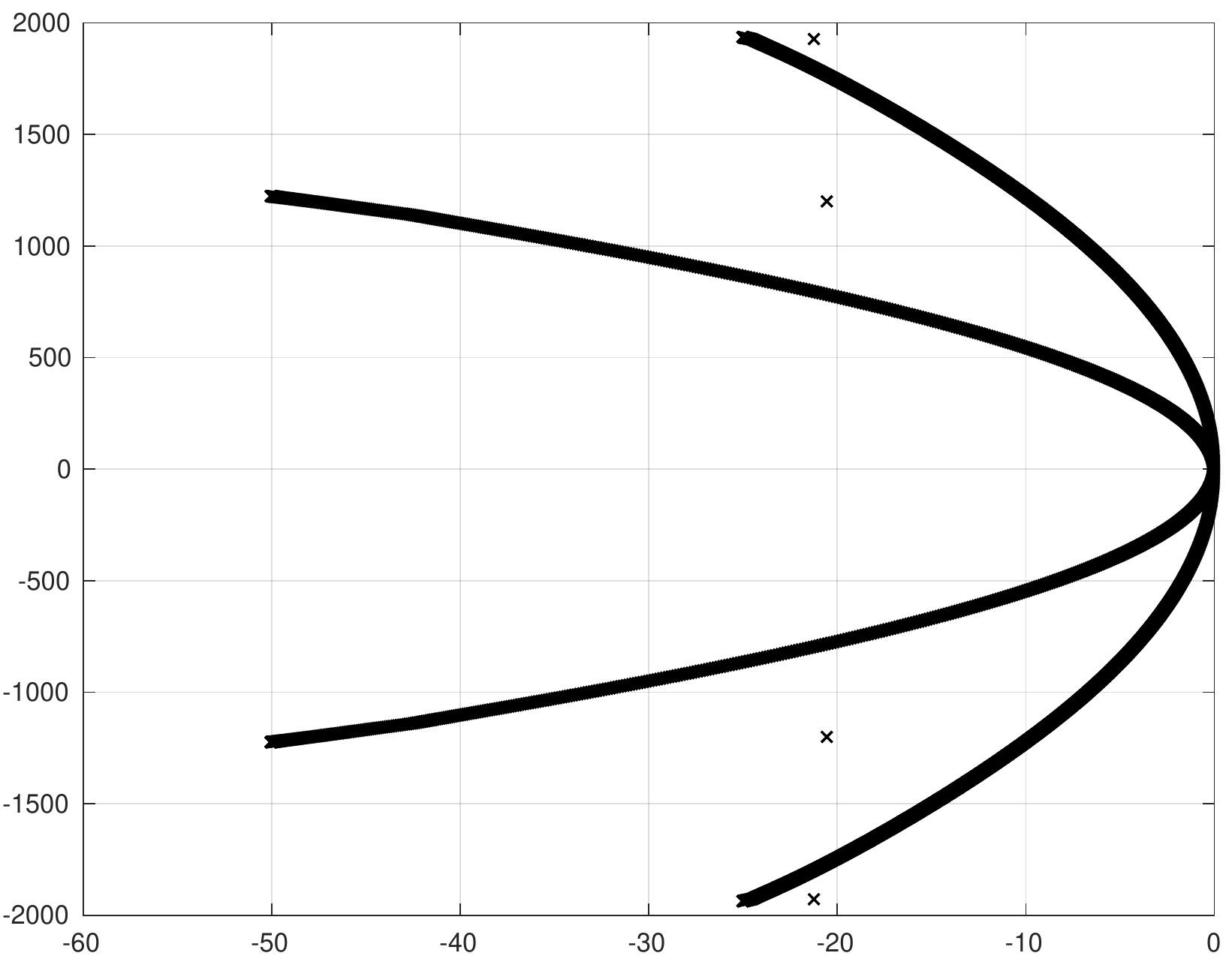}
\includegraphics[scale=0.21]{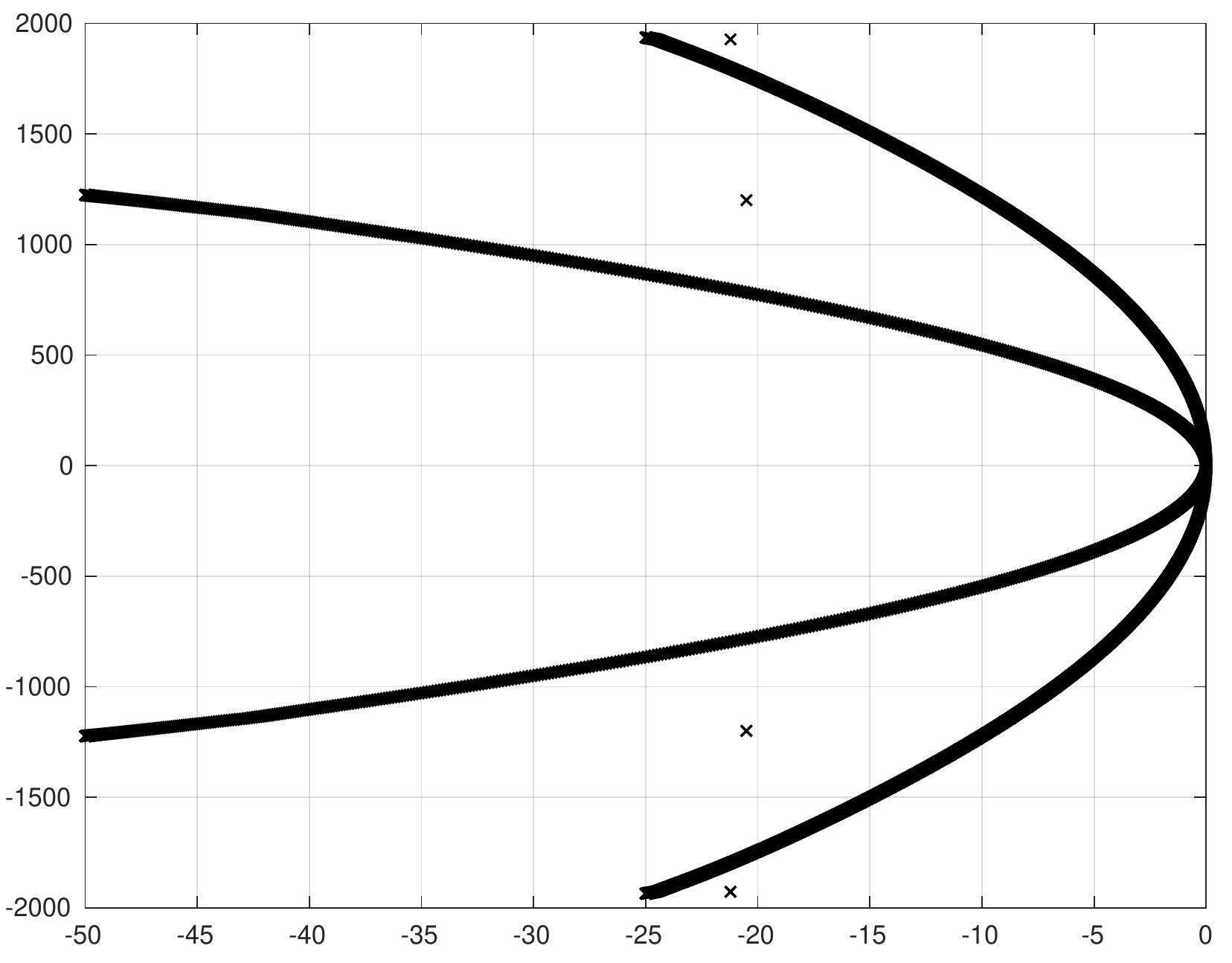}
\includegraphics[scale=0.21]{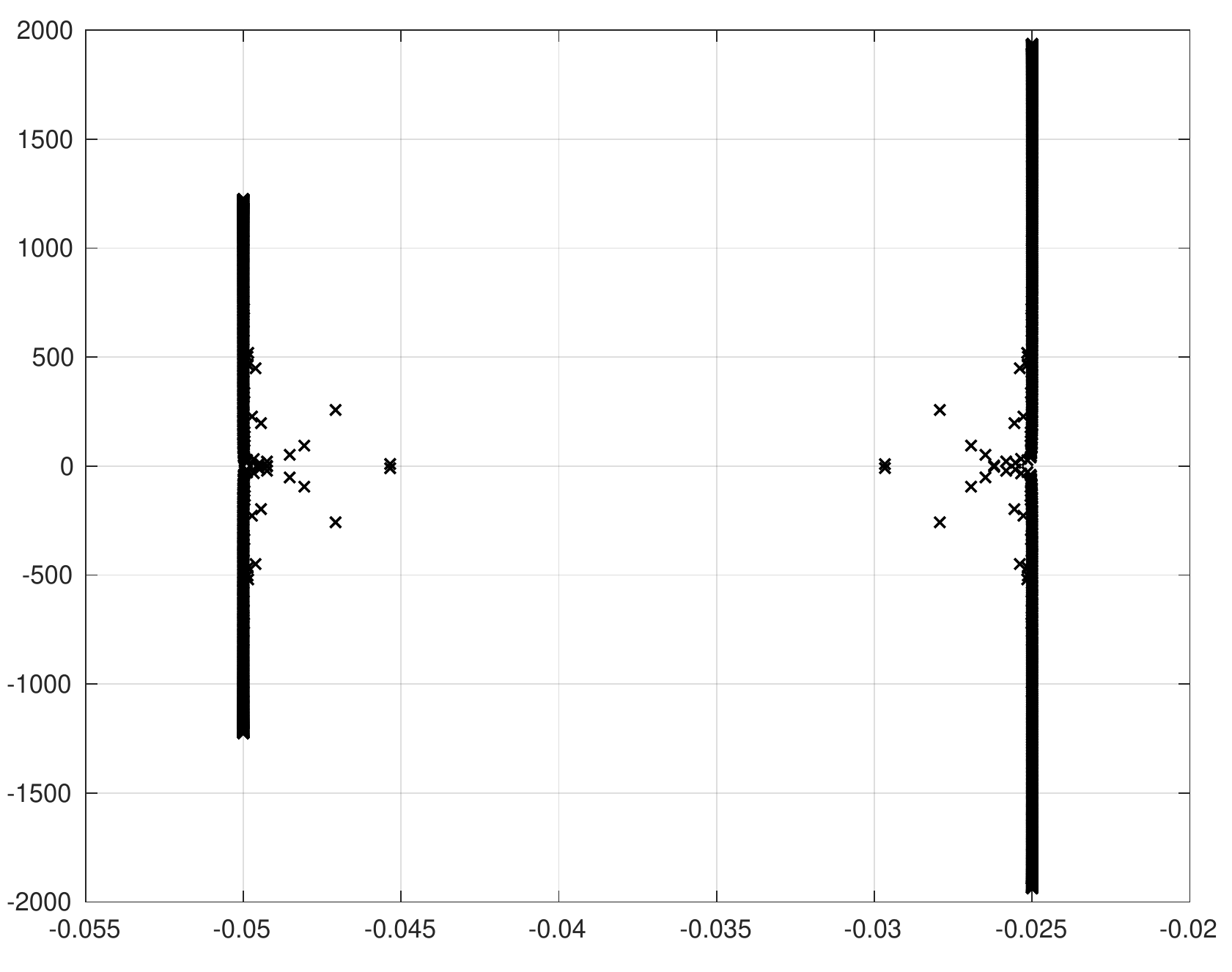}
\includegraphics[scale=0.21]{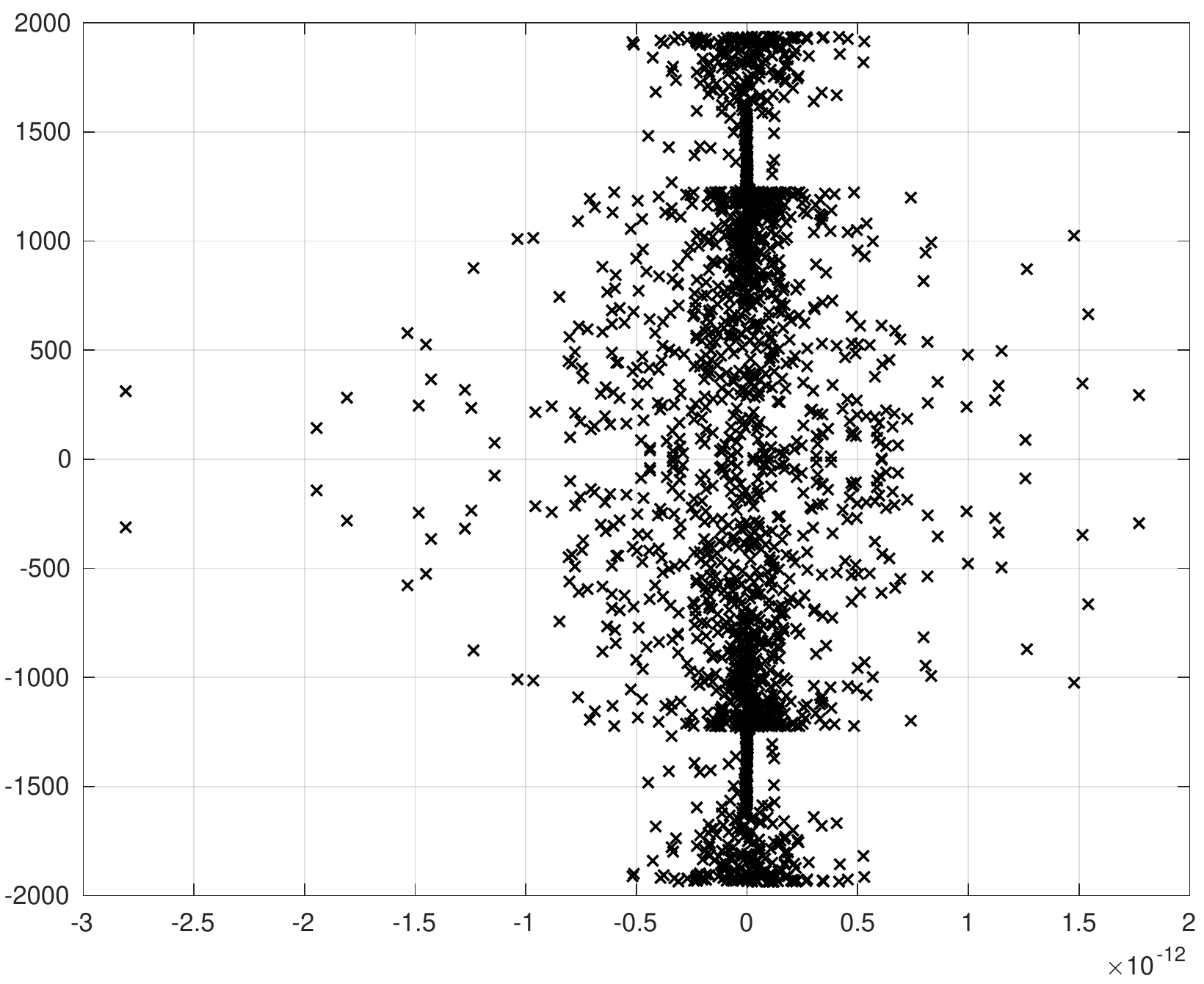}
\caption{Open-loop poles of (\ref{timoshenko}), (\ref{bdry}), computed via the real-rational model $G_r$ with $N=500$. From left to right: both types of damping, KV-damping only, viscous damping only, no damping.}
\end{figure}

According to our general procedure
the transfer function of the pre-stabilized system
(\ref{timoshenko}), (\ref{prestab}) is now computed by solving a succession of elliptic boundary value problems
\begin{align}
\label{laplace}
\begin{split}
    Kw'' &= \rho s^2w + K\phi' - f\\
    EI\phi'' &= I_\rho s^2 \phi + K(\phi-w') -g\\
    w(0,s)&=0, \phi(0,s) =0\\
    K w'(L,s) &- K\phi(L,s) + \alpha s w(L,s) = u_1(s)\\
    EI \phi'(L,s) &+ \beta s\phi(L,s) = u_2(s)\\
    y_1(s) &= sw(L,s), y_2(s) = s\phi(L,s)
    \end{split}
\end{align}
which for given $s=j\omega$ has to be solved twice, with $u_1=1, u_2=0$, and $u_1=0, u_2=1$, using a solver like
{\tt bvp4c} of \cite{MATLAB}.

For the purpose of comparison and assessment of our method
this potentially  infinite-dimensional $2 \times 2$ transfer matrix function $G(s)$ is matched with
a formally computed transfer function $G_f(s)$, which we obtain by exploiting the special structure of the Timoshenko beam
system. Starting with case (\ref{undamped}),
eliminating $\phi$ from the Laplace transformed system (\ref{laplace})
leads to a fourth order equation
$$
w'''' = p(s) w'' + q(s) w
$$
with primes denoting spatial derivatives, 
where in the undamped case
$$
p(s) = s^2/a+s^2/b, \quad q(s) = -(s^2+c)/b \cdot s^2/a
$$
so that eigenvalues are obtained as
$$
\lambda_i(s) = \pm \sqrt{\textstyle\frac{1}{2}p(s) \pm \frac{1}{2}\sqrt{p(s)^2+4q(s)}}, \;\; i=1,2,3,4,
$$
leading to $w(x,s) = \sum_{i=1}^4 A_i(s) e^{\lambda_i(s)x}$. Going back gives
$$\phi(x,s) = \frac{(c-bs^2/a)w'(x,s) + bw'''(x,s)}{s^2+c},$$ hence
$\phi(x,s) = \sum_{i=1}^4 \frac{(c-bs^2/a)\lambda_i + b\lambda_i^3}{s^2+c}A_i(s) e^{\lambda_i(s)x}$. 
Substituting the four boundary conditions at $x=0$ and $x=L$ leads to a $4 \times 4$ linear system
for the coefficients $A_1(s),\dots,A_4(s)$, to be solved for every $s$ with right hand sides $[0,0,1,0]$ and $[0,0,0,1]$. 
Similar computations are used for the two types of damping, where the expressions of $p(s), q(s)$ are suitably adapted. 

\begin{figure}[ht!]
\begin{center}
\includegraphics[height=0.2\textheight,width=0.4\textwidth]{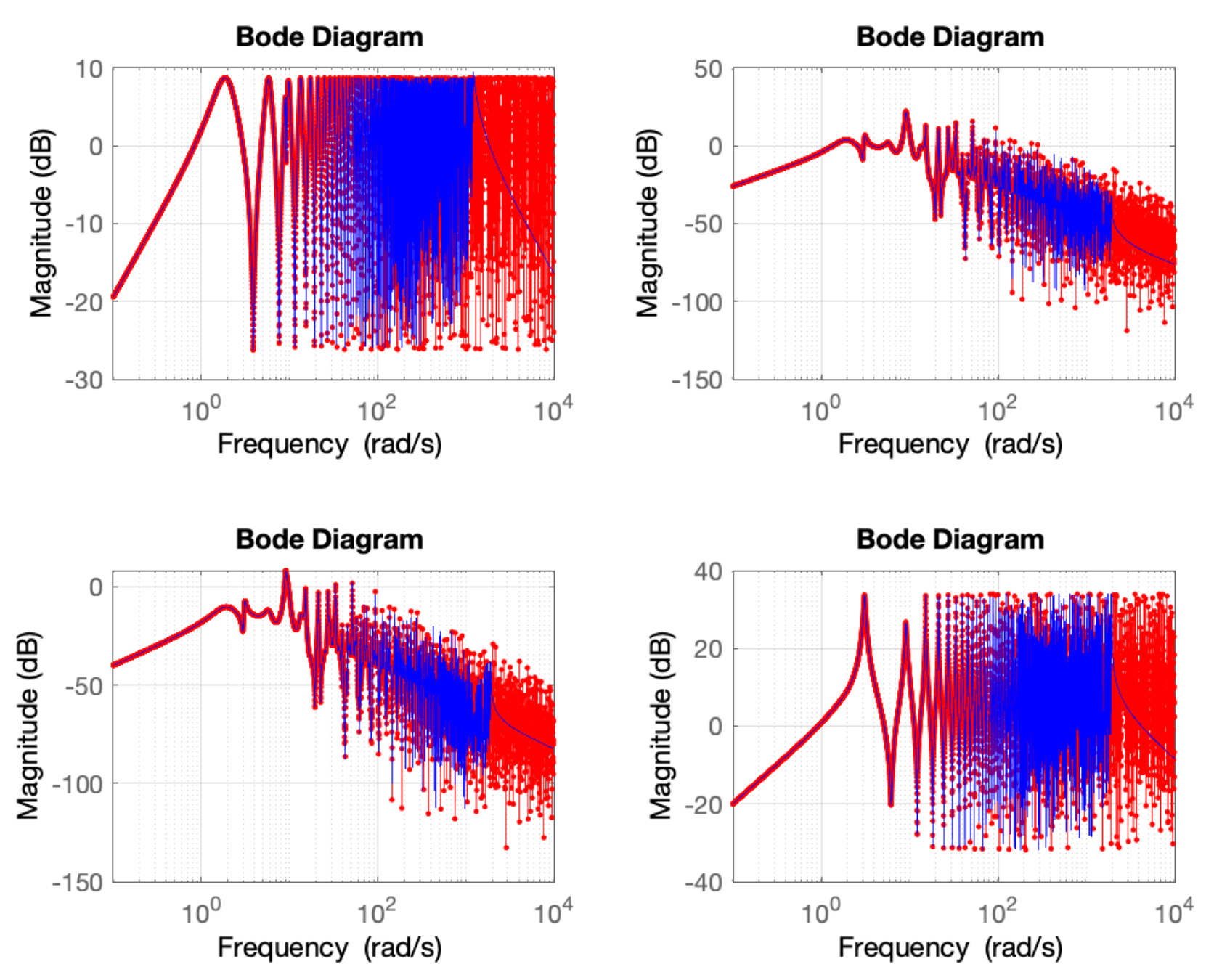}
\includegraphics[height=0.2\textheight,width=0.4\textwidth]{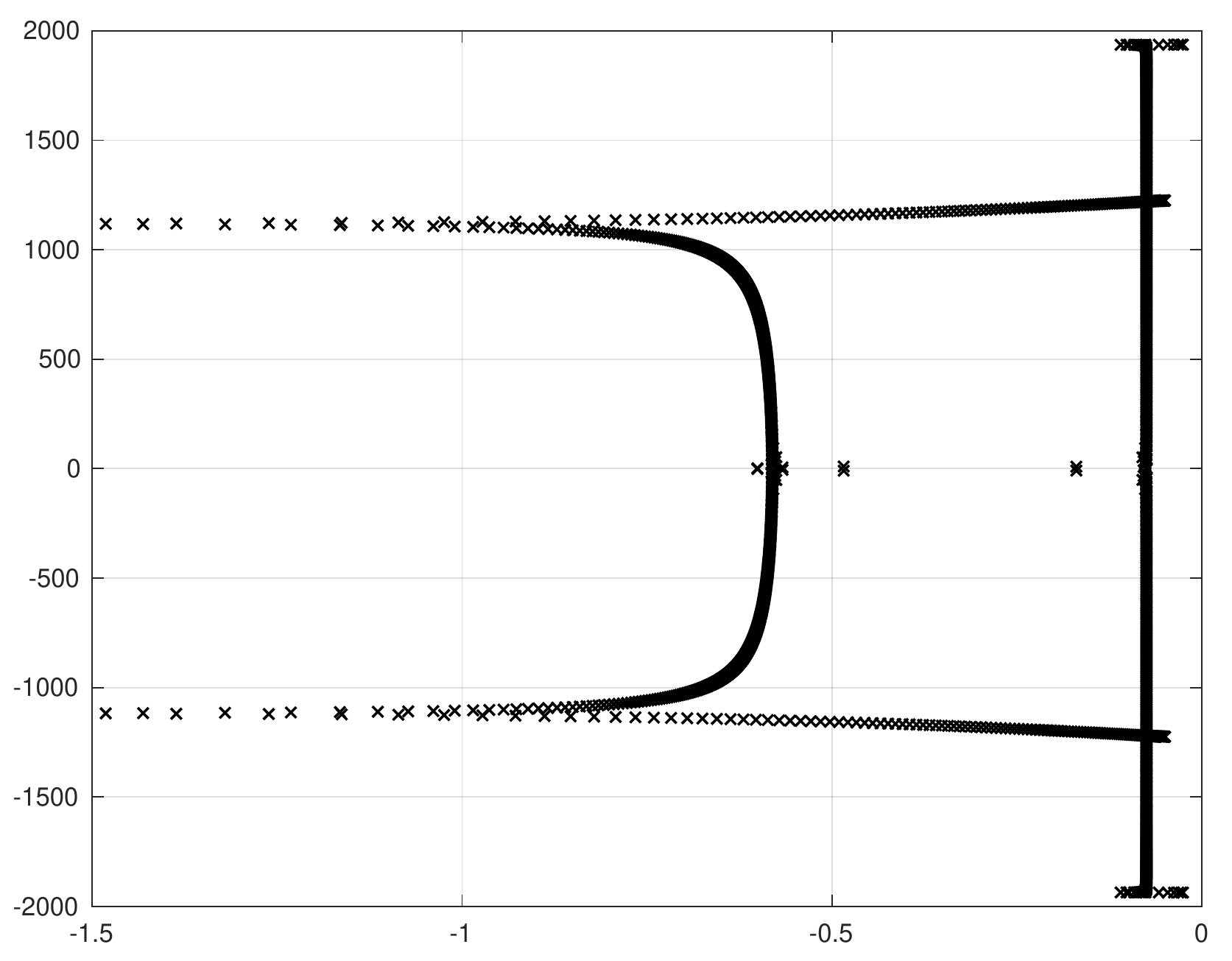}
\includegraphics[height=0.2\textheight,width=0.4\textwidth]{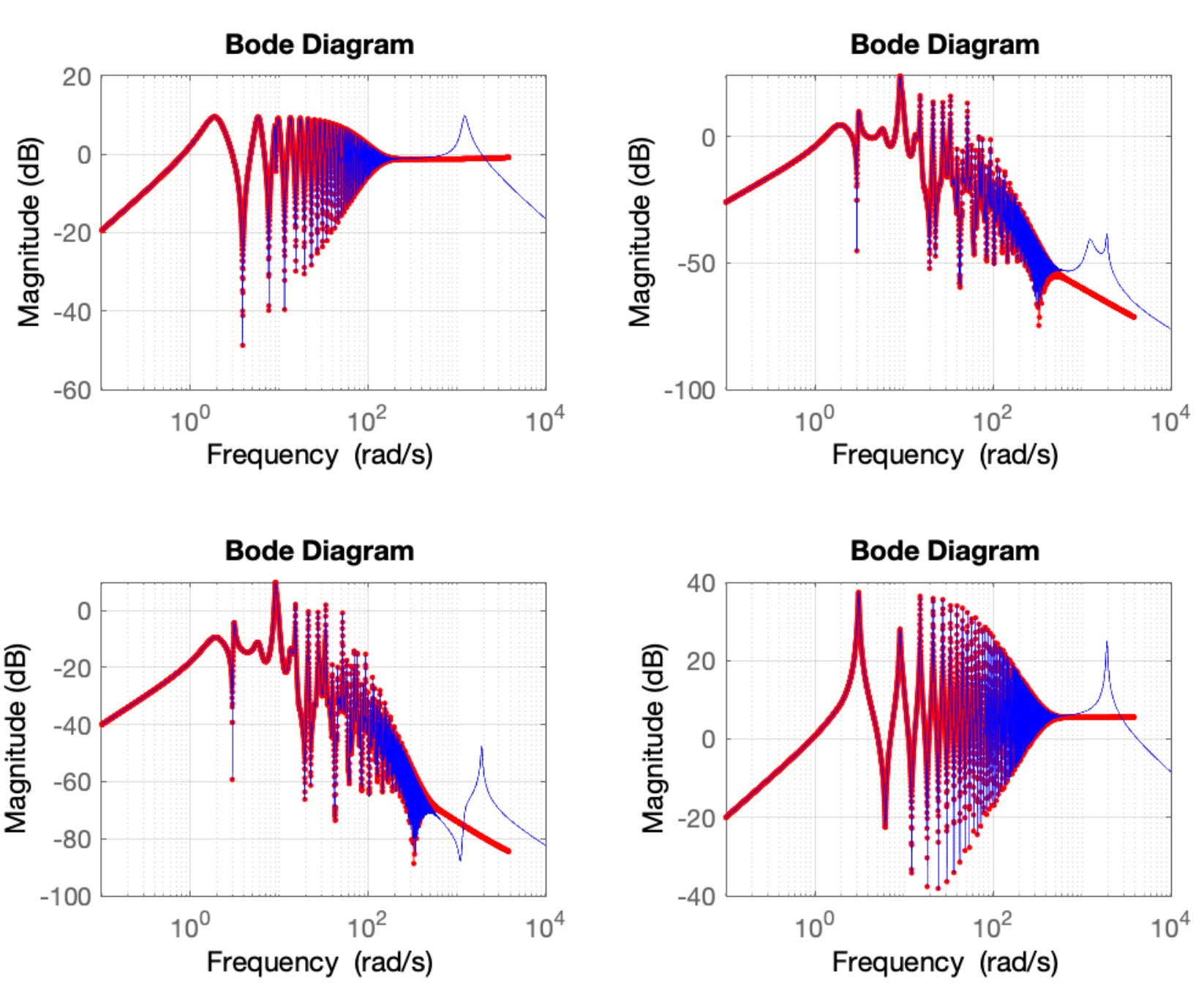}
\includegraphics[height=0.2\textheight,width=0.4\textwidth]{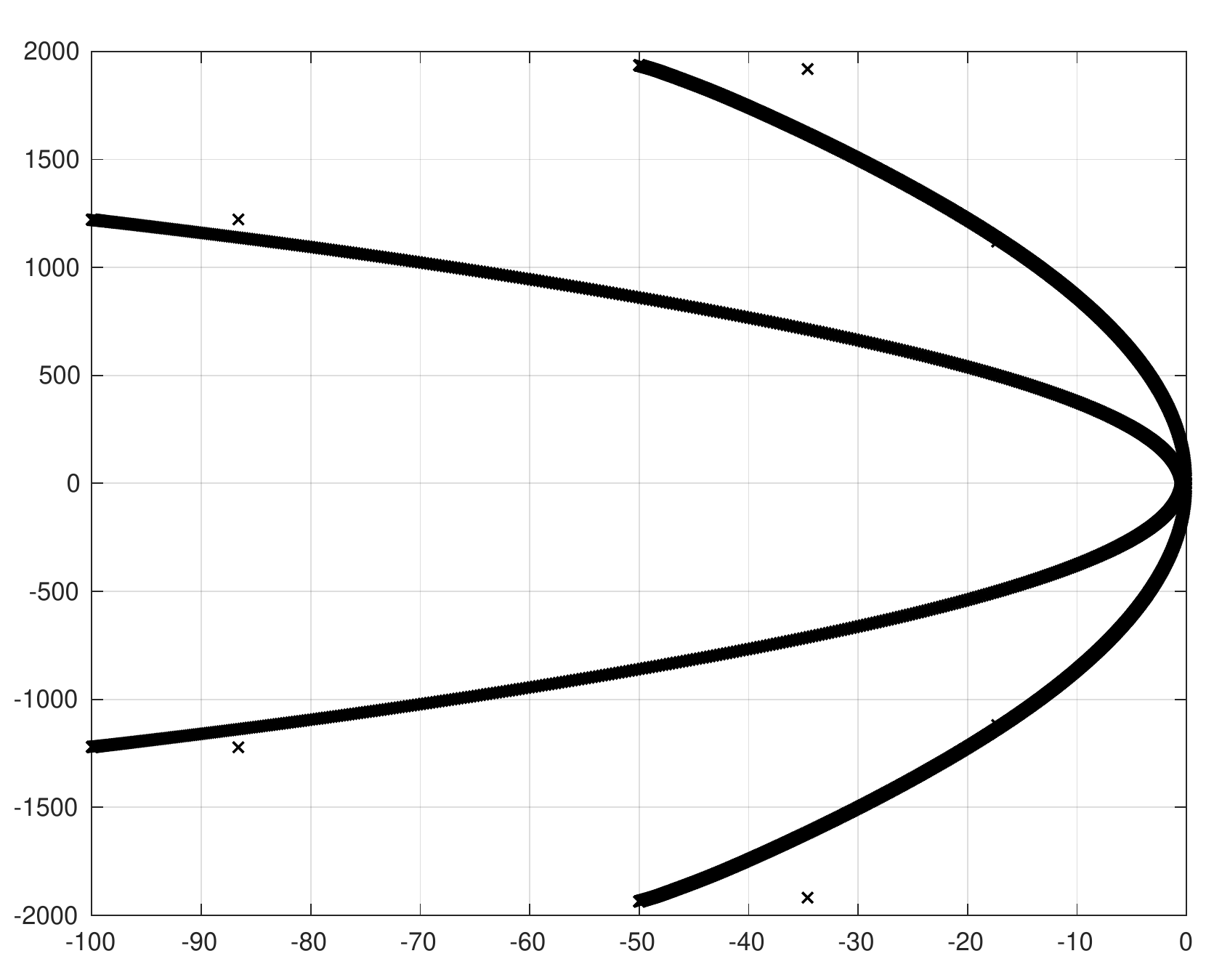}
\end{center}
\caption{Timoshenko beam. $2\times 2$ Bode magnitude plot. 
Formal method (red), 
finite-dimensional approximation with $N=500$ (blue). Viscous damping (upper part)
$d_w=d_\phi = 0.1$, $\alpha=0.5,\beta=0.1$.
Closed loop spectral abscissa $-0.0280$.
KV-damping (lower part), $D_s=D_b=0.0002$, $\alpha=0.5,\beta=0.1$. 
Closed loop spectral abscissa $-0.0487$. 
\label{bode_timo_v}}
\end{figure}

\begin{figure}[ht!]
\begin{center}
\includegraphics[height=0.2\textheight,width=0.4\textwidth]{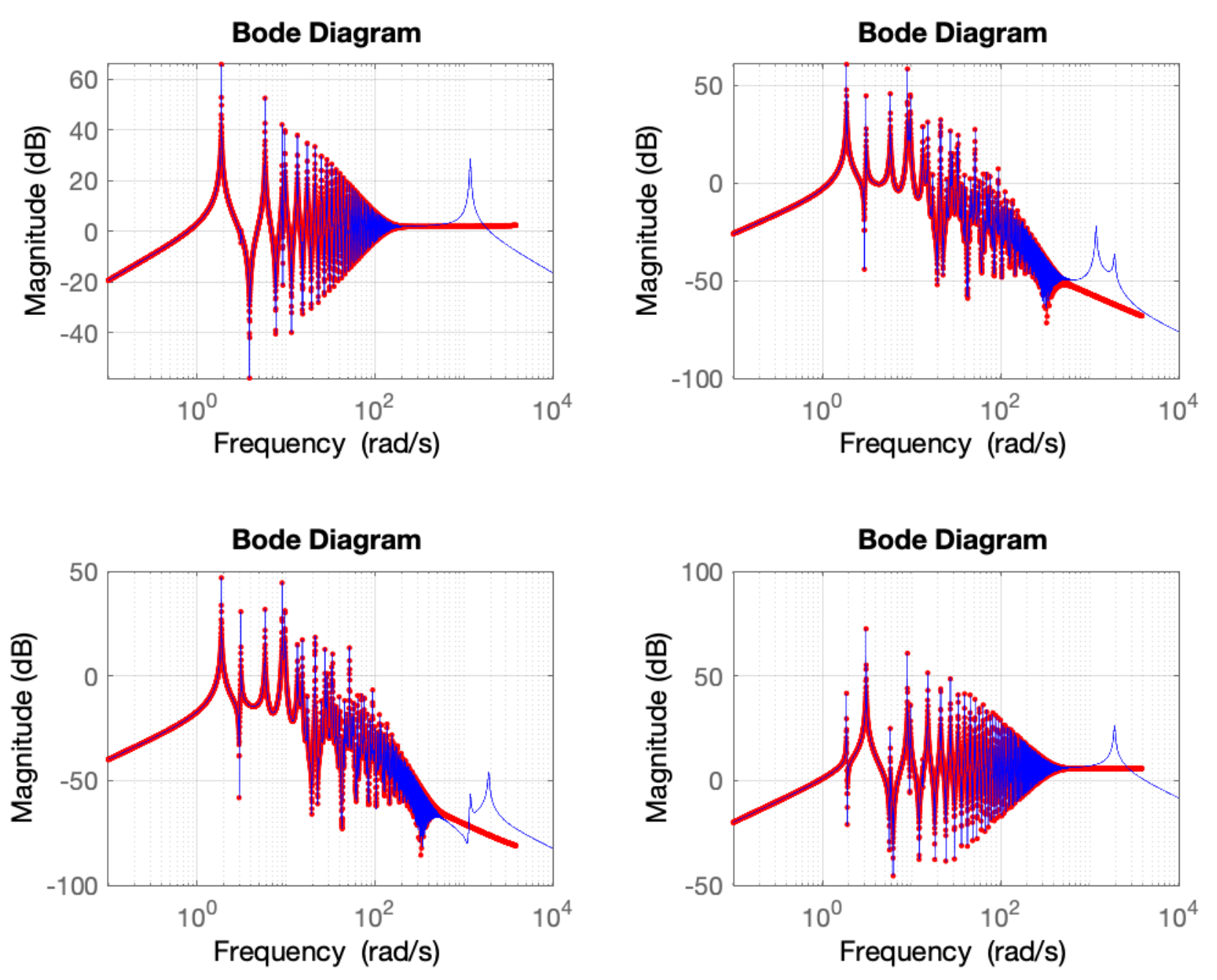}
\includegraphics[height=0.2\textheight,width=0.4\textwidth]{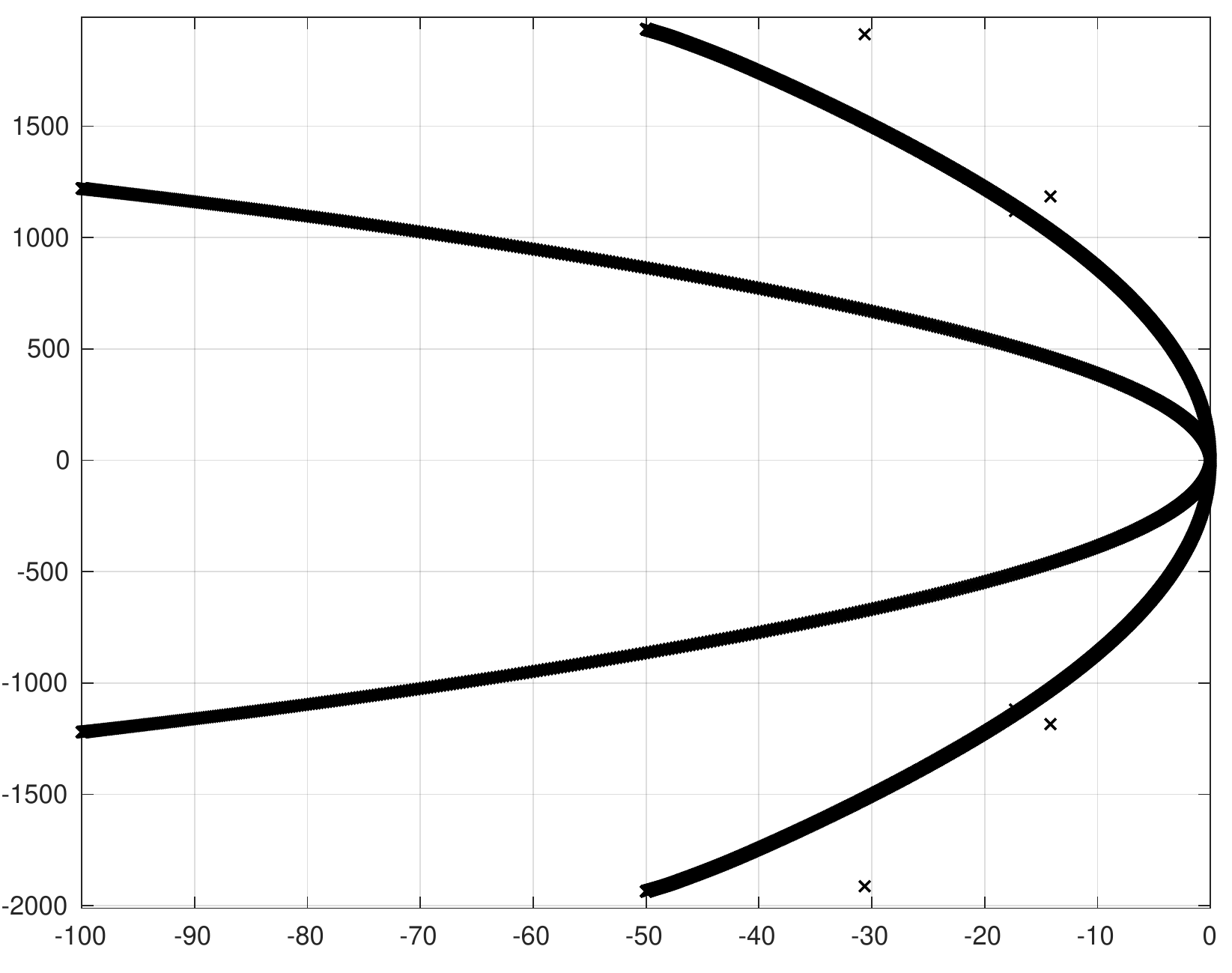}
\includegraphics[height=0.2\textheight,width=0.4\textwidth]{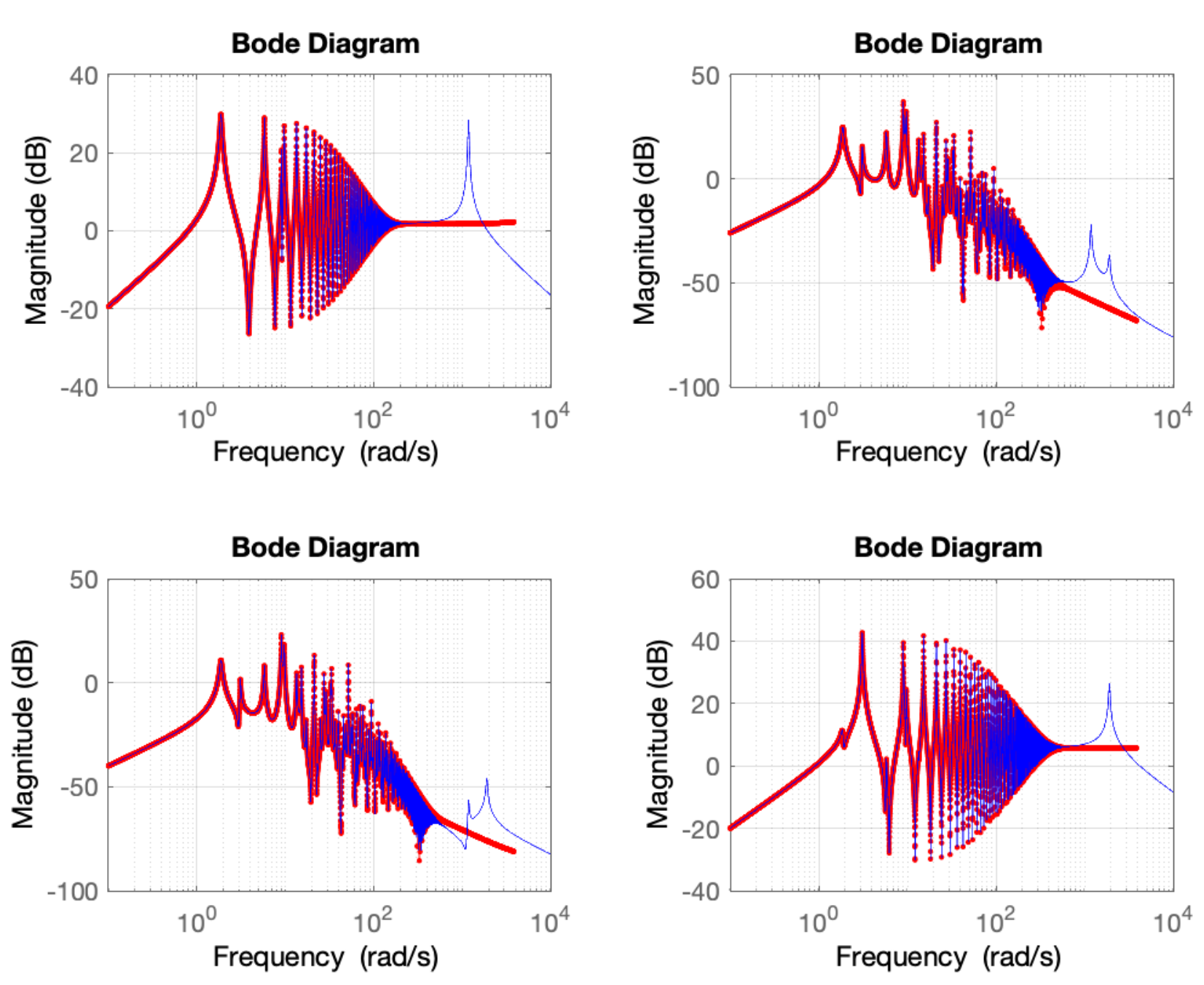}
\includegraphics[height=0.2\textheight,width=0.4\textwidth]{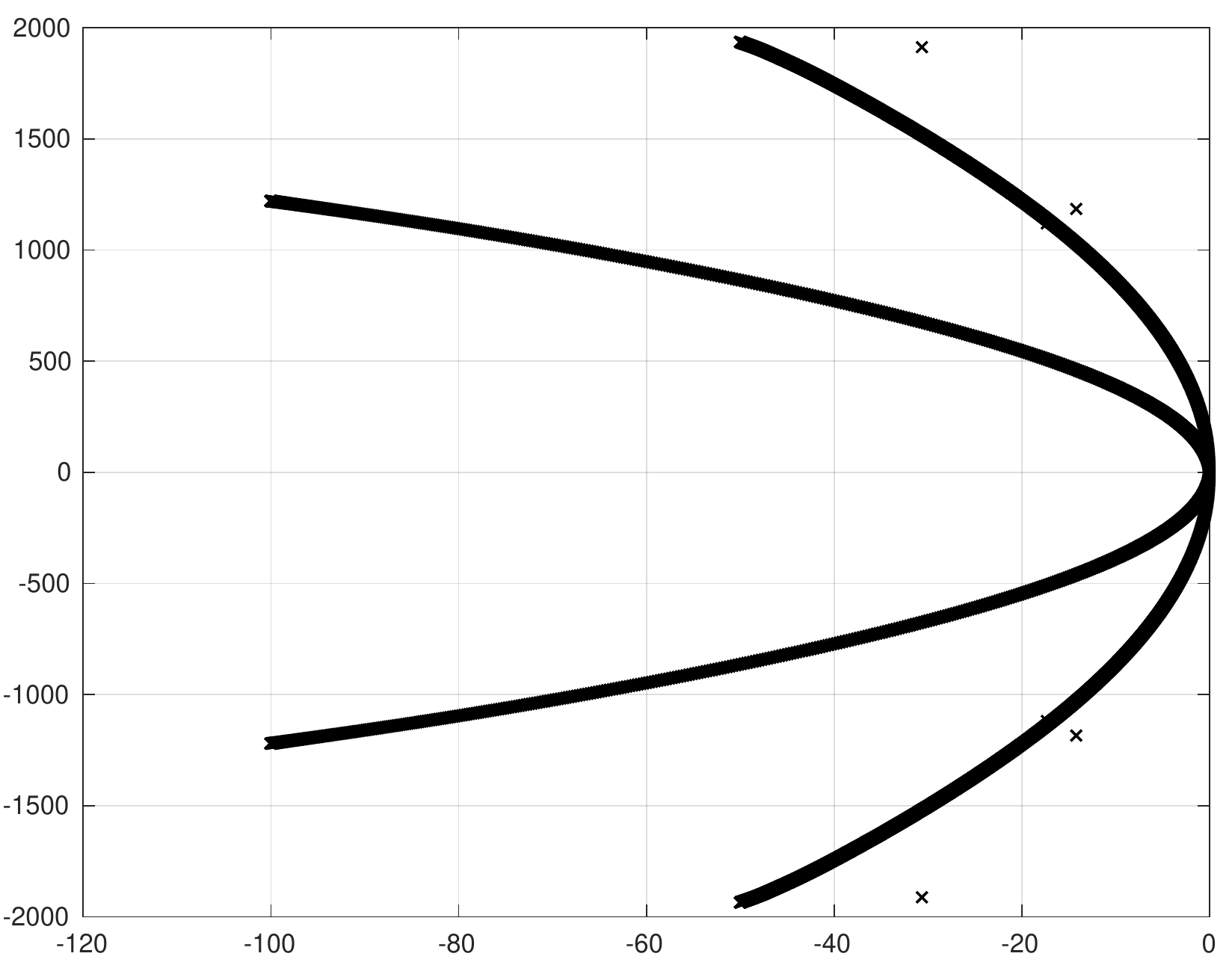}
\end{center}
\caption{Timoshenko beam. $2\times 2$ Bode magnitude  plot. Formal method (red), 
finite-dimensional approximation with $N=500$ (blue). KV damping only (upper part) $D_s=D_b=0.0002$, $\alpha=0,\beta=0$. 
Closed loop spectral abscissa $-1.7090e-04$. KV and viscous damping (lower part)  
$\alpha=0,\beta=0$. 
Closed loop spectral abscissa $-0.0264$. 
\label{bode_timo_kv1}}
\end{figure}

\begin{figure}[ht!]
\begin{center}
\includegraphics[scale=0.4]{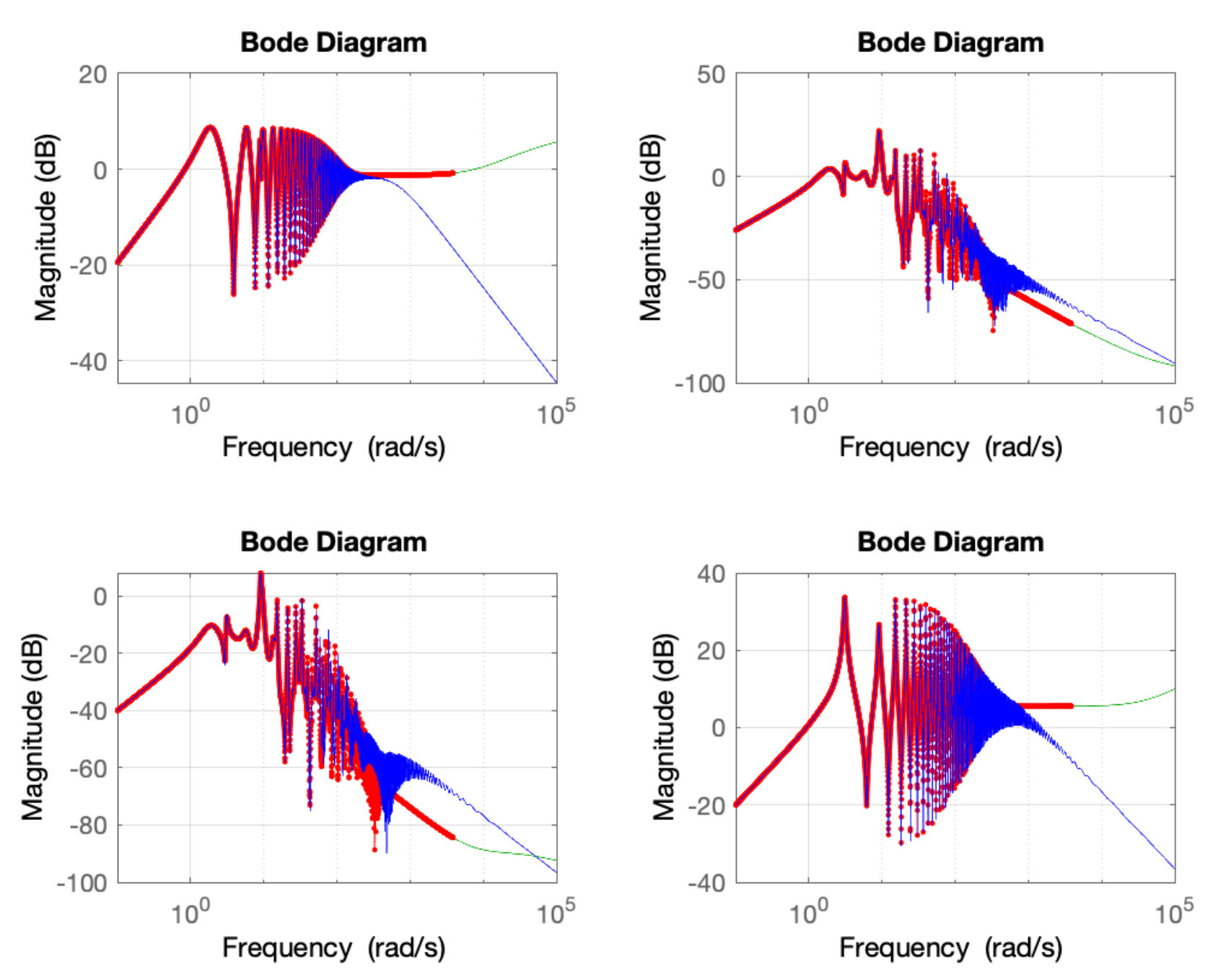}
\end{center}
\caption{Bode magnitude plot of Kelvin-Voigt damped Timoshenko beam. Finite-dimensional approximation with $N=100$ blue,
succession of boundary value problems (green), formal method (red). 
 \label{bode_timo_kv}}
\end{figure}

For comparison, we also use a finite-difference approximation based on a finite-volume
approach \cite{sun2007finite}, which leads to a descriptor system $(E,A,B,C,0)$. 
Note that straightforward finite-difference schemes fail here, as they are not stability preserving
and develop spurious unstable modes for large $N$.
Comparison
of the various transfer functions are shown in Figs. \ref{bode_timo_v}-\ref{bode_timo_kv}.

\begin{remark}
Using numerical data from \cite{kim_renardy}, computation of the spectrum
of $A$ shows good agreement with the results of a model analysis applied directly to the PDE. Presently the approximation $G_r=(E,A,B,C,0)$
turns out stable, which as we know is indispensable for its use in control design. Comparison of the three ways to compute the transfer matrix
are shown in Fig. \ref{bode_timo_kv}. We observed that all methods agree over a wide range of frequencies, with deviations located in the high frequency range which, to some extent, must be taken into account when designing the controllers.
\end{remark}

We end this section by addressing briefly the stability aspects. A good survey is given in 
\cite{liu2018stability}. The undamped case is already covered by \cite{kim_renardy}, and their technique
is easily seen to extend to the viscous damping case. We have the following

\begin{proposition}
Let $K$ be a finite-dimensional controller which stabilizes the Timoshenko beam
$G$ in {\rm (\ref{timoshenko})} with boundary conditions {\rm (\ref{bdry})} and no damping,
or with viscous damping {\rm (\ref{viscous})}, in the $H_\infty$-sense. Then the loop $(G,K)$ is even exponentially stable.
\end{proposition}

\begin{proof}
It follows from \cite{kim_renardy} that $G$ can be exponentially stabilized by the
proportional control $K_0= {\rm diag}(\alpha,\beta)$ in (\ref{prestab}), and we write $G_0 = (G,K_0)$ for the loop. The proof is given for the undamped case, but is seen to carry over to the case of
a global viscous damping (\ref{viscous}).
Let $K = K'+K_0$, then we have to show that the loop $(G_0,K')\simeq (G,K_0+K')=(G,K)$
is exponentially stable. Since $K'$, as a finite-dimensional system, is exponentially stabilizable and detectable, it remains to prove that $G_0$ is exponentially stabilizable and detectable, as the result will then follows from \cite{morris} in tandem with \cite[Lemma 8.2.7]{Staffans_book}.

Since $G_0$ is exponentially stable, it is also exponentially stabilizable, so it remains to prove that $G_0$ is exponentially detectable. 
As $G_0$ now satisfies the spectrum decomposition assumption \cite[Theorem 5.2.6]{CurtainZwart1995}, it follows from
\cite[Theorem 5.2.11]{CurtainZwart1995} that
in order to check exponential detectability, it suffices to
show that ker$(sI-A) \cap {\rm ker}(C) =\{0\}$, where $A$ is the generator, $C$ the output operator 
of $G_0$. This can be done in the frequency domain. We have to show that the system
(\ref{laplace}) 
with $u_1=0$, $u_2=0$ and observed outputs $y_1(s)=sw(L,s)=0$ and $y_2(s)=s\phi(L,s)=0$ for every
$s\in \overline{\mathbb C}_+$ has only the trivial solution $(w,\phi)=(0,0)$. This leads to
an overdetermined  linear system for the coefficients $A_1(s),\dots,A_4(s)$, 
where we have to satisfy not only the boundary conditions $w(0)=0$, $\phi(0)=0$, $u_1=0$, $u_2=0$,
but also the two conditions on the outputs
$\sum_{i=1}^4 sA_i(s) e^{\lambda_i L}=0$ and
$\sum_{i=1}^4 s \frac{(c-bs^2/a)\lambda_i(s)+b\lambda_i(s)^3}{s^2+c} A_i(s) e^{\lambda_i(s)L}=0$.
These 6 conditions for the 4 unknowns can only be satisfied when $A_i(s)=0$, $i=1,2,3,4$.
\hfill $\square$
\end{proof}

Stability under KV-damping is discussed in \cite{liu1998exponential} and \cite{Zhao2005stability}, where the authors allow parts
of the material to be elastic, others visco-elastic, giving natural conditions under which the system is exponentially open loop stable. From the point of view of 
robust control design these cases represent the same level of difficulty,
so in our experiments we concentrate on the case of global KV-damping, where verification by alternative methods remains easier. In this particular case, the system is 
guaranteed open-loop stable.

\section{Second study: Cantilever Euler-Bernoulli beam}
\label{sect_Euler}
Our second study considers piezo-electric control of a thin cantilever beam, where the Euler-Bernoulli beam model may be considered adequate. 
As a variety of methods for this problem have been collected over the years, this is again an instance where our method can be evaluated. 
The equation is of the form
\begin{equation}
    \label{euler}
EI \frac{\partial^4 w(x,t)}{\partial x^4} + \rho A \frac{\partial^2 w(x,t)}{\partial t^2} + c_v \frac{\partial w(x,t)}{\partial t}
+ c_{kv} \frac{\partial^5 w(x,t)}{\partial t \partial^4x}
=
K_a \left( \delta'(x-x_2) - \delta'(x-x_1) \right) u(t)
\end{equation}
with
boundary conditions
\begin{equation}
    \label{euler_bdry}
w(0,t)=0, w_x(0,t)=0, w_{xx}(L,t)=0, w_{xxx}(L,t)=0,
\end{equation}
where $L$ is the length of the beam, which is clamped at $x=0$ and free at $x=L$, $E$ is Young's modulus, $I$ the
moment of inertia of the beam, $\rho$ its density, $A$ the beam cross section, and where the constant $K_a$
depends on width, thickness and piezoelectric strain constant of the actuator, with $u(t)$ representing the applied voltage. Coefficients $c_v,c_{kv}$ stand for viscous and Kelvin-Voigt damping. 
The values
$0 < x_1 < x_2 < L$ indicate the left and right end of the sensor/actuator pair. The measured output is
$$
y(t) = K_s \left(  w_x(x_2,t) - w_x(x_1,t) \right)
$$
for a constant $K_s$ now depending on properties of the sensor.
The goal is to compute a finite-dimensional simply structured controller $u = Ky$ which stabilizes the structure and, in addition, allows to
attenuate the induced vibrations.

\begin{figure}[ht!]
\includegraphics[scale=0.27]{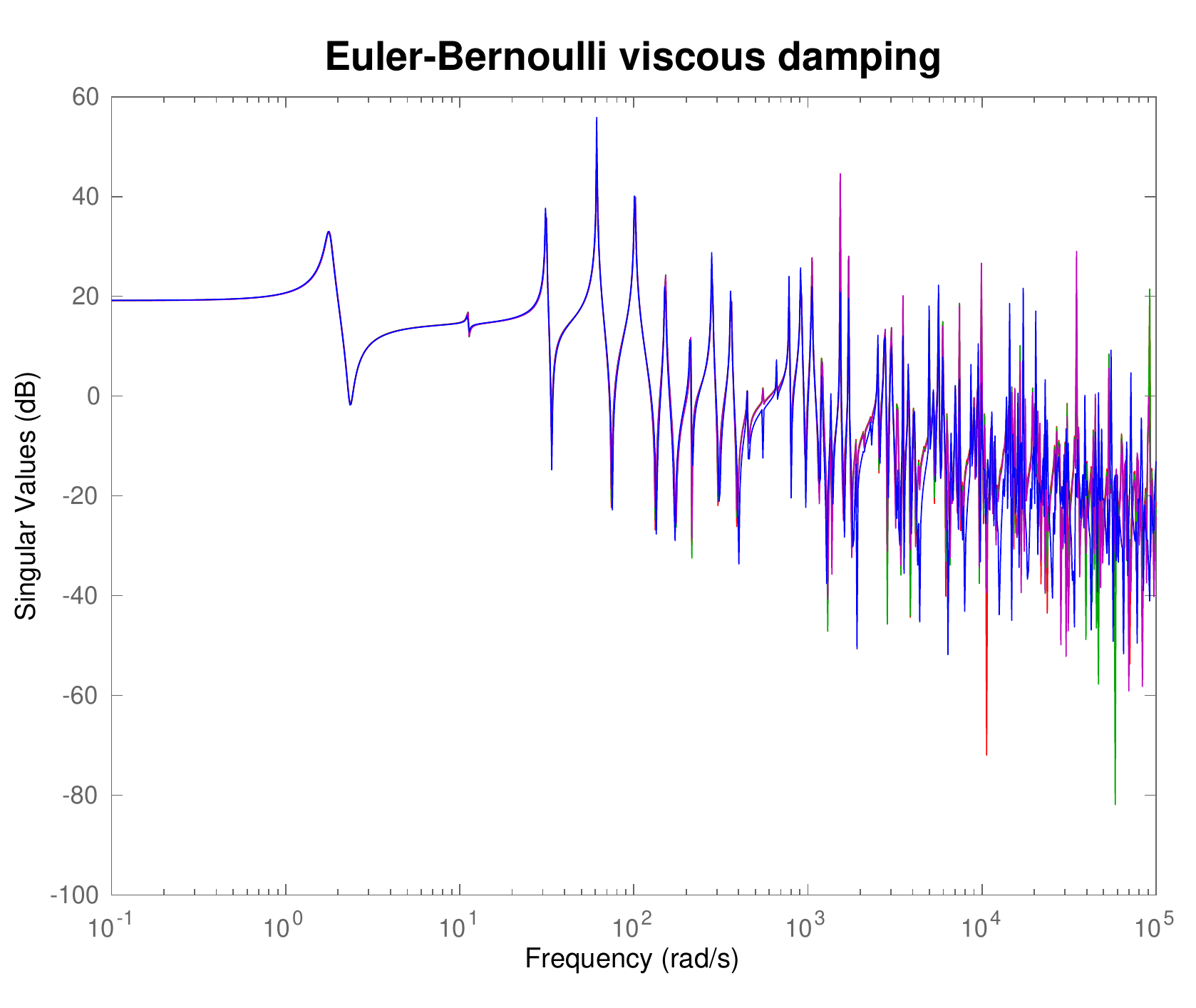}
\includegraphics[scale=0.27]{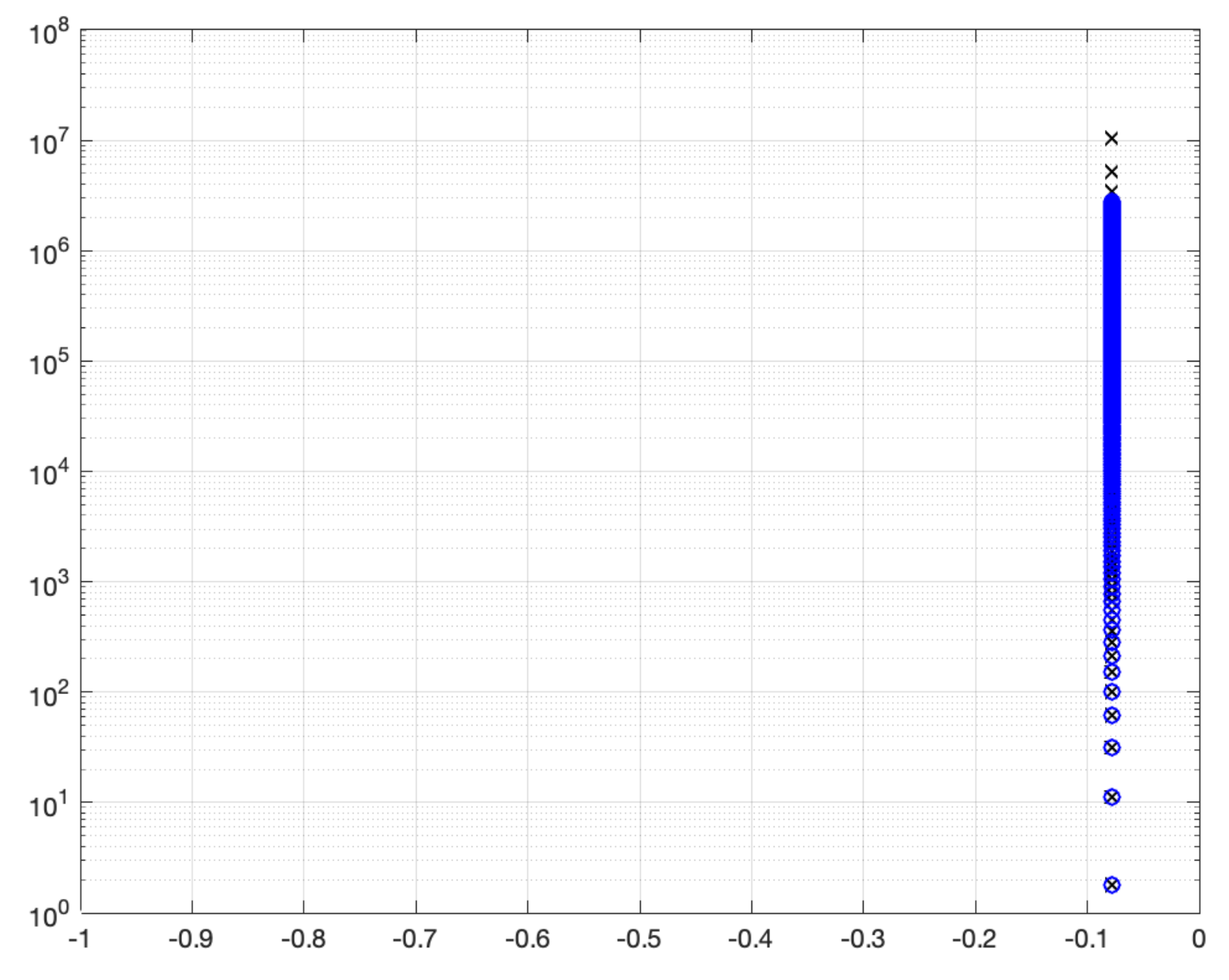}
\includegraphics[scale=0.27]{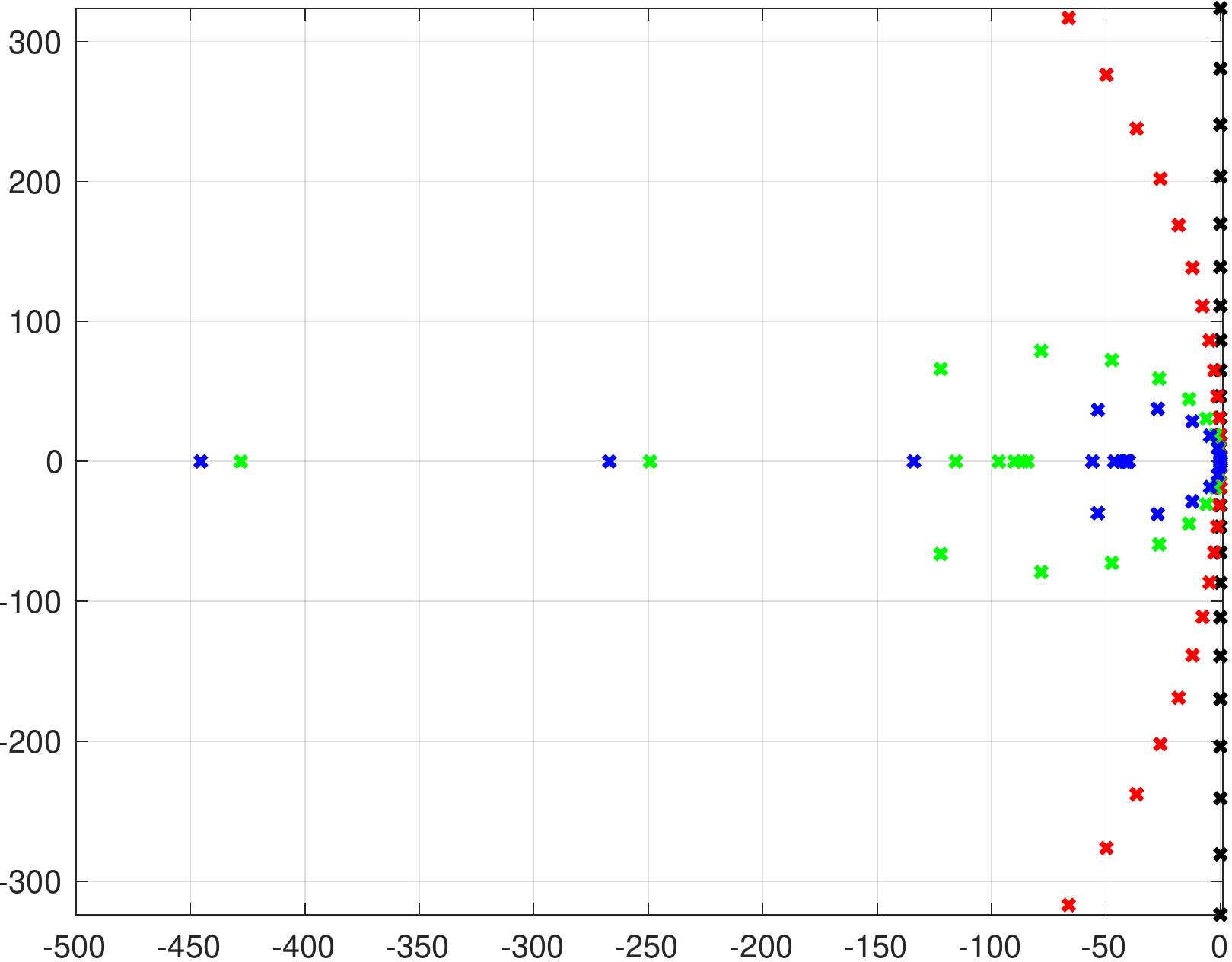}
\caption{Euler-Bernoulli beam. Left: Transfer functions computed via elliptic boundary value problems (magenta), semi-formal method (red), diagonalization (green) and descriptor system (blue) ($N=400$). Middle: 
Poles of descriptor system compared to theoretical poles (blue) $c_v=0.156,c_{kv}=0$.
Right: Poles for $c_v=0.156, c_{kv} \in \{0,0.0001,0.001,0.002\}$
\label{TFeuler}}
\end{figure}

While the transfer function $G(s)$ can again be computed to arbitrary precision as $G(s)=y(s)/u(s)$ by
solving a succession of complex elliptic boundary value problems
\begin{equation}
\label{elliptic_euler}
(EI+c_{kv}s) w''''(x,s) + (\rho A s^2 + c_vs) w(x,s) = K_a \left( \delta'(x-x_2) - \delta'(x-x_1) \right)u(s)
\end{equation}
with boundary conditions
$$
w(0,s)=0, w'(x,s)=0,w''(L,s)=0, w'''(L,s)=0
$$
and 
$$
y(s) = K_s(w'(x_2,s)-w'(x_1,s)),
$$
we evaluate the results using alternative ways to get $G(s)$. This includes diagonalization
and expanding $G(s)$ into a series of eigenfunctions, and a classical {\em Ansatz}
with $w(x,s) = \sum_{i=1}^4 A_i(s) e^{\lambda_i(s)}$, $\lambda_i(s)$ the eigenvalues
of (\ref{elliptic_euler}), fitting the $A_i(s)$ through the boundary conditions;
see \cite{gurgoze2006dynamic}.

A finite-dimensional stability preserving approximation
based on a descriptor system obtained from a finite-volume type discretization is discussed
in \cite{liu2019novel}. On the right hand side we have approximated $\delta(\cdot-x_i)$ 
using a Gaussian $\phi(\cdot-x_i)$ centered at $x_i$, 
which gives $\delta'(\cdot-x_i)$ as $\phi'(\cdot-x_i)$. 
Comparisons of these transfer functions are shown in Fig. 
\ref{TFeuler}.
Note again that
the straightforward second-order difference scheme fails, as it introduces spurious modes which
get unstable as $N=hL$ increases. 

Realistic models for the Euler-Bernoulli beam should
include
damping, and
in our experiment we concentrate on the viscous damping case, 
as from the automatic control point of view 
this is the severest case.
As stressed in \cite{awada2022optimized}, 
design of an appropriate control law in tandem with an optimal choice of the sensor/actuator
positions is the key to a successful vibration suppression in a smart structure. 
The authors of \cite{awada2022optimized}  give a good overview on previous attempts 
based on various controller structures such as LQG, PID or simple proportional control laws. 
Presently we address this problem via structured $H_\infty$-control.

The viscous damped system is open-loop stable,  \cite{gurgoze2006dynamic}, while 
stabilization of the undamped beam is for instance discussed in \cite{Le2007output,guo2001riesz,paunonen2022finite}, where the authors
use a non-realizable D-controller to stabilize the loop, which corresponds to the tip load damping
of \cite{gurgoze2006dynamic}.
In \cite{wang2008stability} the authors use an  infinite-dimensional observer.

Open-loop poles $s_k$ of (\ref{euler}) are obtained from standard semi-group theory,
which gives the relation
$$
r_k^4 = \frac{\rho A s_k^2 + c_v s_k}{EI+c_{kv} s_k}, \quad 1+\cos(r_kL) \cosh(r_kL)=0
$$
where the $r_kL$ are pre-computed with arbitrary precision.
In tandem with \cite[Thm. 5.2.6]{CurtainZwart1995} this formula
shows that the undamped Euler-Bernoulli beam cannot be 
exponentially stabilized by a finite-dimensional controller.

\section{Feedback control of  Timoshenko beam}
\label{sect_timo_control}
In this section, we discuss how a finite-dimensional controller 
for the  Timoshenko beam can be designed. 
Physical parameters for (\ref{timoshenko}) are chosen as
$$ L = 1, \rho =1, K = 1.5, I_\rho = 2, E = 2.5, I = 3\,
$$
adopted from \cite{kim_renardy}, where for the
pre-stabilizing control (\ref{prestab}) the values
$\alpha = \beta = 0.1$ are chosen. Frequency responses of all  possible combinations are shown in 
Fig. \ref{sigmaGinfAll}. 
We study the cases undamped, viscous damping ($d_w = d_\phi = 0.5$) with pre-stabilizer $K_0$, and Kelvin-Voigt damping ($D_s = D_b = 1\text{e-}4$) with and without pre-stabilization. Cases with both types of damping (fourth of Fig. \ref{sigmaGinfAll})  are  the easiest 
and will not be treated to save space.

\begin{figure}[h!]
\begin{center}
\includegraphics[scale=1.0]{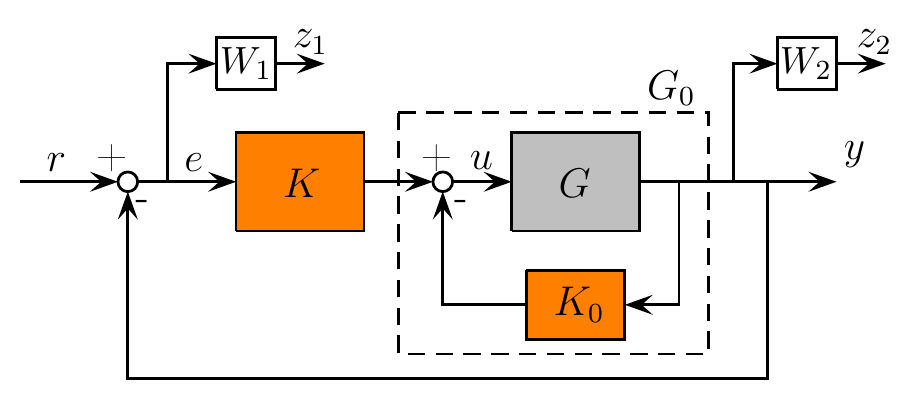}
\end{center}
\caption{Feedback design interconnection. \label{beam_scheme}}
\end{figure}

The design problem is formulated as a classical $2\times 2$ reference tracking problem, where the measured outputs,
shear and bending moment, should follow in a decoupled fashion  reference inputs such as steps or ramps. 
Tracking performance requires 
minimizing a weighted sensitivity function $W_1(s)S(s)$, where $S = (I+GK)^{-1}$ and
$W_1(s)$ is a low pass filter specifying the frequency range on which tracking should be achieved. With
tracking error $e = r-y$, reference input $r$, and measured outputs $y$, this is shown in Fig. \ref{beam_scheme}.
Good tracking requires a high-gain $W_1$ at $s=0$ to limit or eliminate errors in steady state, as well as in a range $\omega \in [0, \omega_b ]$, where $\omega_b$ is the bandwidth to meet rise, settling time and overshoot constraints. 

In addition to tracking, the controller should also attenuate resonant modes in the high frequency range. This is addressed by minimizing the transfer function from $r$ to $y$.  
With $y = GK(I+GK)^{-1} r$ this leads to
minimizing $W_2(s) T(s)$, where $T(s):= GK(I+GK)^{-1}$ is the complementary sensitivity function
and $W_2$ is a high-pass filter, whose role is to 
specify the frequency range where resonant modes are critical. Minimizing $T$ has the additional beneficial effect to
improve robustness of the design 
against a loss of model fidelity. As discussed in Section \ref{sect-Timo}, exact and approximative methods agree at low frequencies, but significant discrepancies are observed in the high frequency range. This clearly suggests 
 a high-pass filter $W_2$. Namely, minimizing $T$ in the high frequency range then
generates roll-off in the loop gain $GK$, thereby mitigating the effect of the 
infinitely many resonant high frequency modes. 
This also facilitates computations, both for the Nyquist stability test and for performance and robustness criteria.
Further details on how the Nyquist stability can be put to work are given in \cite{AN:18} and references therein.

\subsection{Specific measures for the infinite dimension}
\label{special}
For finite-dimensional systems $G(s)$,
specifications based on $S$ and $T$ in tandem with closed-loop stability are 
usually highly effective \cite{KWAKERNAAK2002}, especially when filters $W_1,W_2$ are properly chosen.

In infinite dimensions, $G(s)$ is generally only known on a finite set
of sample frequencies $s_\nu = j\omega_\nu$, even though this set can potentially be enlarged and refined when required. 
This introduces an additional difficulty,
to which we respond by adding further constraints, which we now discuss.

Putting a bound on
the disk margin \cite{PS1993,dullerud:00,ZDG:96} is a means to prevent the Nyquist plot 
from getting too close to the critical point. 
This is expressed as $\| S \|_\infty \leq  \gamma$, where $1/\gamma$ is the disk margin. This constraint has
a repelling effect against stable poles crossing the imaginary axis and becoming unstable. 

Two further requirements are recommended. First, one aims at realistic implementable controllers, which may be achieved 
by selecting classes $\mathscr K$ of simple low-order controllers. 
For controllers $K\in \mathscr K$  to be practical, however, 
we also enforce constraints on its dynamics. 
This may be translated into so-called disk constraints on eigenvalues 
$|\lambda_i(K)| \leq \delta$, $\delta >0$, $i = 1,\ldots,n_K$, where $n_K$ is the order of $K$. 
More compactly, we write $\rho(K) \leq \delta$, where $\rho(.)$ denotes the spectral radius. 

In the same vein, selecting sample frequencies $\omega_\nu$ clearly requires minimal 
knowledge of the plant $G$, and the various methods discussed in Sections \ref{sect-Timo} and \ref{sect_Euler} 
should help finding suitable discretizations. 
As we will see later  in our experiments, 
the choice of frequency samples must be based on the characteristics of the open-loop models and must also cope
with control design constraints in programs (\ref{problem}) or (\ref{pgST}) mentioned above. 
But even when the system $G$ is ideally discretized, distortions may still occur outside the sets $\Omega$. 
Those are referred to as {\it inter-sample} distortions, inherent to any discretization approach.
To mitigate inter-sample effects, we propose  to rule out too small damping in the controllers dynamics. 
This may be written as a constraint $\zeta(K) \geq \mu$, where $\zeta(.)$ is a max-function on the 
damping of controller eigenvalues. 

Putting all those together leads to the following optimization program,
which is a special case of (\ref{problem}):
\begin{eqnarray}
\label{pgST}
\begin{array}{ll}
\mbox{minimize} & \displaystyle\max  \{  \| W_1 S\|_{\infty,\Omega_S},\, \| W_2 T\|_{\infty,\Omega_T} \}\\
\mbox{subject to} & \displaystyle \|S \|_{\infty,\Omega_D} \leq \gamma\\
 & \rho(K) \leq \delta \\
 & \zeta(K) \geq \mu \\
& K \in \mathscr K \mbox{ stabilizes $G$}
\end{array}
\end{eqnarray}
Here $\| M(s)\|_{\infty,\Omega}$ is short for $\max_{\omega \in \Omega} \max_i \sigma_i (M(j\omega))$ with  $\sigma_i$ denoting singular values. Sample frequencies $\Omega_S,\Omega_T$ are adapted to $W_1S$ and $W_2T$,  $\Omega_D$ is adapted to
the disk margin constraint, and sampling $\Omega_N$ to the Nyquist test. This gives extra flexibility and can be exploited to reduce execution times and to address specific frequency bands.

Program (\ref{pgST}) is solved using nonsmooth bundle or trust region optimization techniques, \cite{noll2019cutting, ANR:17,ANR:2015,AN:2007,Noll10}. However,
a new element is required for infinite-dimensional closed-loop stability.
In finite dimensions 
this may be assured by 
a constraint $\alpha(A_{cl}(K)) \leq -\epsilon$ on the closed-loop spectral abscissa.
In infinite dimension, particularly for non-sectorial systems, 
this in no longer possible, as $\alpha(A_{cl})$ is not reliably computable. 
Instead we use a workaround based on the Nyquist stability criterion, which we now explain.

As soon as $G$ is pre-stabilized, closed-loop stability is maintained during 
optimization over $K\in \mathscr K$
by keeping the winding number of $f(s)=\det(I+G(s)K(s))$ at the correct
value.
This requires that the number of unstable poles of iterates $K\in \mathscr K$
does not change either.  Presently this is assured since the damping constraint $\zeta(K) \geq \mu$ 
restricts the search $K\in \mathscr K$  to stable controllers $\alpha(K) < 0$, which means that
the correct winding number for closed-loop stability is 0.
As soon as a step $K^+ = K + dK$ leads to an unstable loop, recognized by a
non-zero winding number,
backtracking $K^+ = K + \alpha dK$, $0 < \alpha < 1$, or tightening proximity control, are used to maintain stability. In the next step it is then necessary to
prevent the optimizer from  tempting the direction $dK$ again. As a rule this may be achieved
by using the closed-loop sensitivity function as a barrier, as has been pointed out in
\cite{AN:18}. One has just to be aware that $\|S(K)\|_\infty$, while having
a large peak at points $K$ where the closed-loop turns unstable, may show misleading small values
for points $K^+$ farther beyond the point of instability, so does not behave in the same way as say a log-barrier function. Our experiments indicate that this difficulty is avoided through the combined use of the barrier and the Nyquist test.

The Nyquist test gives a stability certificate as soon as it
can be proved that between two sample frequencies
$\omega_\nu$ and $\omega_{\nu+1}$, the closed curve $\gamma_\nu$ obtained by concatenating the true Nyquist arc   
$\{f(\omega): \omega\in [\omega_{\nu},\omega_{\nu+1}]\}$
with the segment $[f(\omega_{\nu+1}),f(\omega_\nu)]$ does not encircle the origin. This can for instance be arranged
when a prior bound on the variation $f'(\omega)$ of
$f(\omega)={\rm det}(I+K(j\omega) G(j\omega))$ can be provided. We call a function $L[\omega^-,\omega^+]$ a first-order bound
if $L[\omega^-,\omega^+]\geq |f'(\omega)|$ for all $\omega\in [\omega^-,\omega^+]$. Then
any sampling $\omega_\nu$ satisfying $L[\omega_\nu,\omega_{\nu+1}](\omega_{\nu+1}-\omega_\nu) \leq |f(\omega_\nu)|+|f(\omega_{\nu+1})|$
gives a provable certificate of stability via the Nyquist test. See \cite[Lemma 3]{AN:18},
and Algorithm 1 of that reference how to construct the $\omega_\nu$. 
Even when no rigorous bound $L[\cdot,\cdot]$ is available, numerical bounds usually give excellent results. Note that the bound comes into play only for those parts of the Nyquist contour 
which are relatively close to the origin.

The Nyquist test has precedence over all other functions in (\ref{pgST}). Clearly, objectives and constraints need not be computed when the Nyquist test fails at $K^+=K+dK$, which 
speeds up computations. 

Appropriate sampling with certificate
$\Omega_S$, $\Omega_T$ and $\Omega_D$ for $W_1S$, $W_2T$ and $S$ required for the constraints
and max objective of (\ref{pgST})
is also described in \cite[Lemma 3]{AN:18}. We also write $\Omega_N$ 
for sampling for the Nyquist criterion. In practice, it makes sense to take the same sampling grid $\Omega_N = \Omega_D$, because the disk margin and the Nyquist test are  equivalent measures of singularity, that is, $\underline{\sigma} (S(j\omega)^{-1}) = \underline{\sigma} (I+G(j\omega)K(j\omega)) = 0$ if and only if $\det(I+G(j\omega)K(j\omega)) = 0$. 

Note that undersampling $\Omega_S$ and  $\Omega_T$ may lead to underestimating the values of the $H_\infty$-norms in the cost function of (\ref{pgST}), 
but is less critical than undersampling $\Omega_N$, which may
put stability at stake.

In our experiments, 
we have adopted two strategies to deal with the beam tracking problem. The first one is based on deriving a reduced-order approximation of the beam dynamics giving a descriptor state-space model
$(E,A,B,C,D)$. 
The second uses the idea of \cite{AN:18}, which avoids approximations 
and relies on frequency domain data of the beam obtained from (\ref{laplace}). This is more in line
with data-driven control techniques \cite{hast2013pid,karimi2017data,kergus2017frequency}, even though we have the possibility 
to supplement additional data if required. It leads  to sampled data $\{G(j\omega): \omega \in \Omega\}$ for a finite $\Omega\subset [0,\infty]$.  Both strategies have been discussed in Section \ref{sect-algo}.

\subsection{Undamped Timoshenko beam\label{sect-TimoU}}
\begin{figure}[ht!]
\includegraphics[height=0.25\textheight, width = 0.3\textwidth]{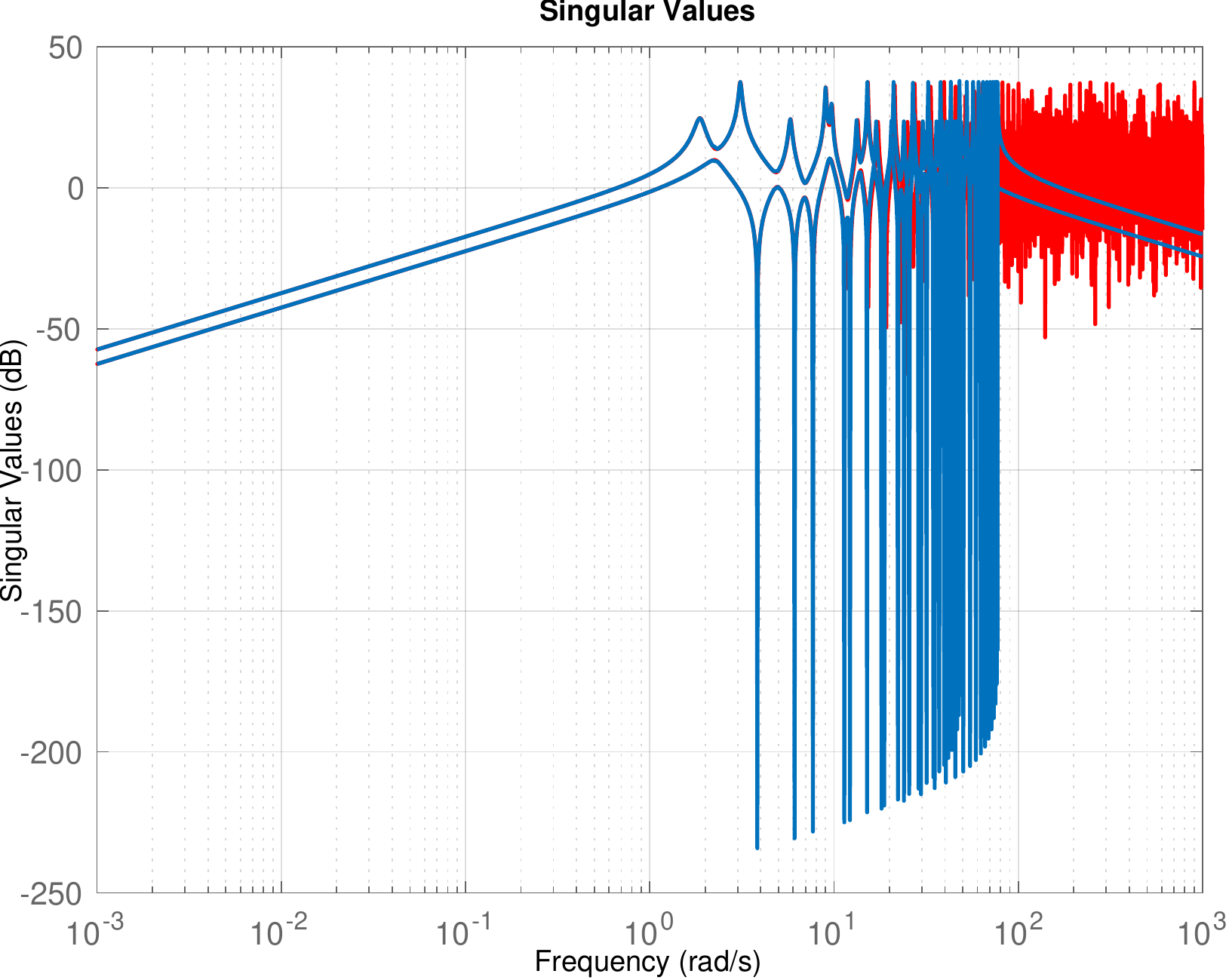}
\includegraphics[height=0.25\textheight, width = 0.3\textwidth]{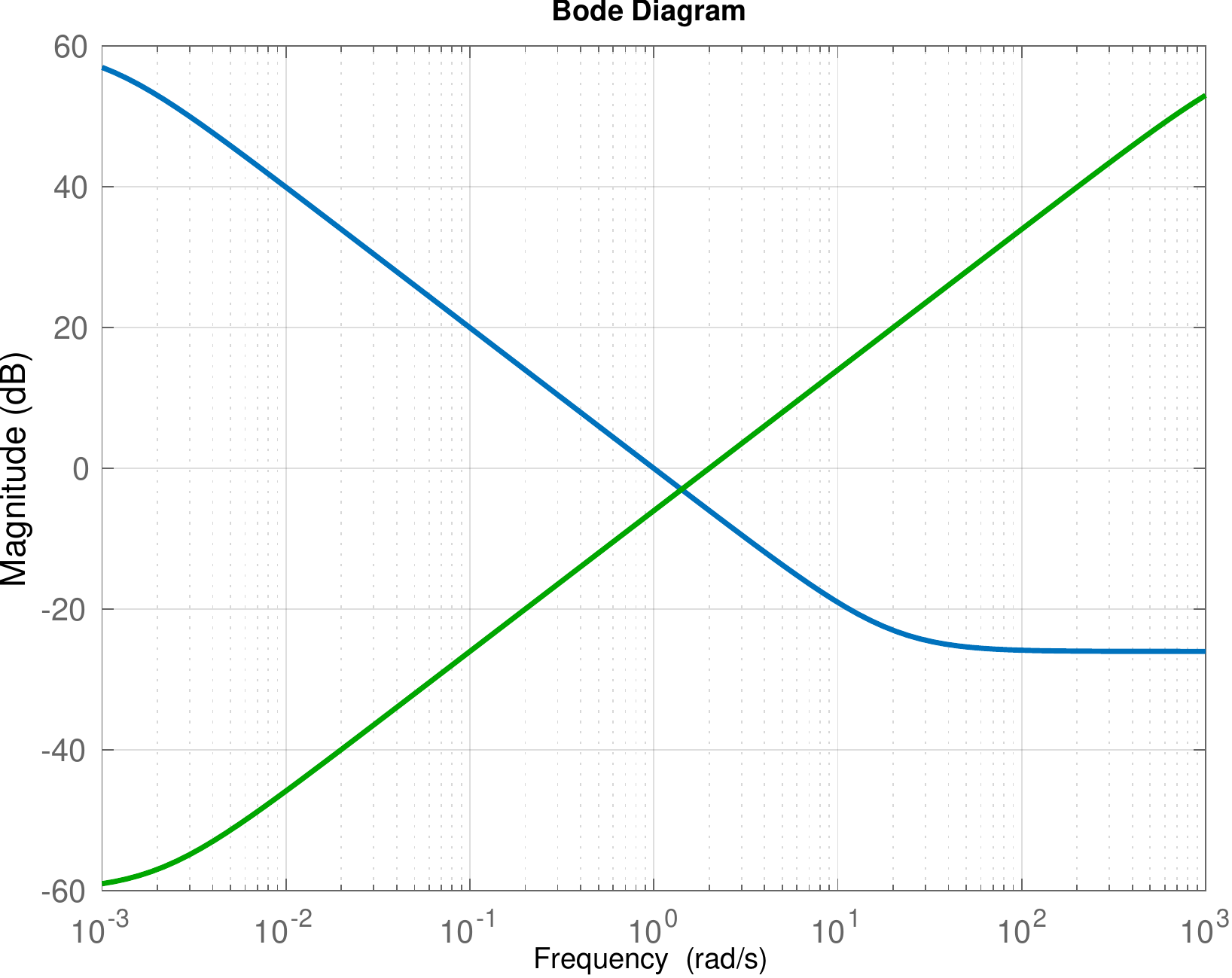}
\includegraphics[height=0.25\textheight, width = 0.3\textwidth]{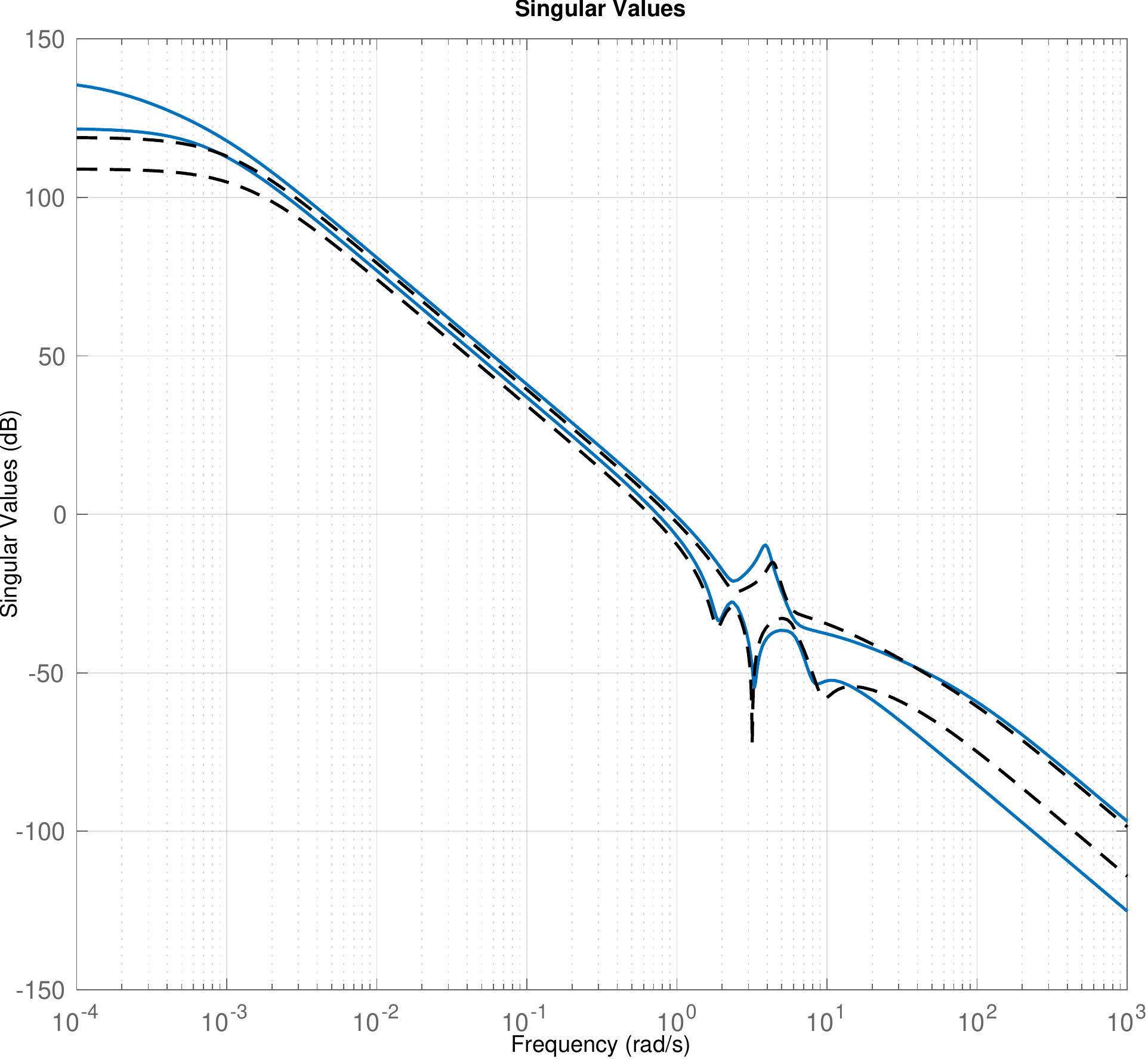}
\caption{Undamped Timoshenko beam. Left: singular values of infinite-dimensional model $G$
and reduced LTI model $G_r$ of order $40$. Middle: weighting functions used in (\ref{pgST}). Right: singular values of controllers $K_r$
based on reduced-order model (solid), and $K$ based on infinite-dimensional-model (dashed).  \label{figG}}
\end{figure}

We start with the case of the undamped Timoshenko beam, which is certainly the most challenging, numerically and in terms of feedback design. 
A reduced-order model retaining $20$ resonant modes, that is, of order $40$ and exact data are compared in Fig.
\ref{figG} left. These models $(E,A,B,C,0)$ and  $\{G(j\omega): \omega \in \Omega\}$ essentially differ in the high frequency range, where the reduced-order model disregards resonant modes. 

Good tracking and decoupling properties would dictate using integral action in the controller. This is not possible here, however, because measurements are of derivative type. In response we use a pseudo-integrator $1/(s+\varepsilon)$, directly included into the structure of
the controller as $K := 1/(s+\varepsilon) \widetilde{K}$, where $\widetilde{K}$ is the optimized portion, while $\varepsilon = 1\text{e-}3$ remains fixed. 
This can also be directly embedded into the cast (\ref{pgST}) by 
the substitutions $G \leftarrow G /(s+\varepsilon) $. Constraints on $K$ now turn
into constraints on $\widetilde{K}$. As structural  
constraint we choose $\widetilde{K} \in \mathscr K_8$, the set of $8$th-order controllers, a choice
motivated by  preliminary testing. 
This leaves $54$  tunable parameters or optimization variables $\x$ corresponding to tridiagonal state-space realizations. Those realizations do not limit generality and are less costly than full state-space forms. Note that including a pseudo-integrator which plays the role of a pre-compensator into the plant 
$G$ is reminiscent of the loopshaping approach \cite{MC90,McFG:1992}. 

Designs based on a reduced-order state-space model and on the 
infinite-dimensional model use the same cast (\ref{pgST}) 
to allow comparison. With a state-space model $G$, program (\ref{pgST}) can be solved using {\tt systune} \cite{RCT2021b,AN2006disk,ANtac:05}, whereas its data-based variant has to rely on the more recent technique in \cite{AN:18}. 

Weightings for tracking $W_1$ and roll-off $W_2$ shown in Fig. \ref{figG} middle are
$$
W_1(s): = \frac{0.05 s + 0.9987}{s + 0.0009987},\qquad W_2(s): = \frac{1000 s + 2}{s + 2000}\,.
$$
A disk margin of $1/\gamma$ with $\gamma =0.7$ is specified, which leads to the constraint 
$\|S\|_\infty \leq 1/0.7$.  Due to plant zeros at the origin, it is not possible to minimize the sensitivity function $S$ at the zero frequency. The $H_\infty$ norm $\|W_1S\|_\infty$ 
in (\ref{pgST}) is therefore restricted to the frequency band $[10^{-3}, \, \infty)$.

The controller radius constraint in (\ref{pgST}) is set to $\delta = 100$. Its minimal  damping is set to $\mu = 10^{-2}$. Taken together with the order limit $\widetilde{K} \in \mathscr K_8$ these constraints secure practical and implementable controllers. 

Design (\ref{pgST}) based on frequency-domain data has to trade accuracy against execution times. 
Frequency-domain sampling a system with an infinite number of resonances as in Fig. \ref{figG} may seem daunting, because peak frequencies will almost certainly be missed. 
As an example, $1\text{e}4$ points were used to plot Fig. \ref{figG}, and some plots in Sections 
\ref{sect-Timo} and \ref{sect_Euler} required even more points.  In fact, as far as stability is concerned, it is sufficient to increase the sampling density only
in those frequency ranges where the controller gain is significant. Clearly, frequency ranges where the loop gain $GK$ is small do not contribute much to the Nyquist stability criterion, 
as there $\det(I+GK) \approx 1$. This is easily devised by taking into account the dynamics of the beam in Fig. \ref{figG} along with the design constraints.

In the present application, we distinguish $3$ frequency intervals for the disk margin, 
$\Omega_D = \Omega'_D\cup \Omega_D''\cup \Omega_D'''$ in (\ref{pgST}) and $\Omega_N:= \Omega_D$ for the Nyquist criterion.  In the low frequency range $\Omega_D'=[1\text{e-}8,\, 1]$ rad/s, we have high gain control but dynamics are rigid, thus $300$ linearly spaced samples suffice. In the mid frequency (crossover) range 
$\Omega_D''=[1,\, 15]$, we have moderate gain control but resonant dynamics and it is enough to use $500$ samples. For higher frequencies $\Omega_D'''=[15,\, 1\text{e}3]$, we have low  roll-off control  with resonant dynamics. We therefore cut down density to only $100$ samples. These figures convey the general idea that is followed repeatedly in the sequel. 

\begin{remark}
In order to apply the certificates in  \cite{AN:18}, it is necessary to adapt the sampling grids
$\Omega_S,\Omega_T, \Omega_D, \Omega_N$ not only to $G$, but also to the controller $K$. 
Since $K$ changes during
optimization, this could appear to be a major difficulty due to a large number of updates. Theoretically this can be avoided by
considering optimization with inexact data \cite{noll_inexact}, or by fixing the grid and
re-checking the result in the end, using re-starts if the grid turned out insufficient.
\end{remark}

Fortunately, this is more a theoretical quest than a challenge in practice, as updates 
due to varying $K$ are rare.
Presently
we never had to stop-and-restart (\ref{pgST}) with refined sampling. 

More significantly,
it can be observed that running times fall from $2$ hours for blind sampling with $1\text{e}5$ nodes to less than $2$ minutes using the outlined strategy.

As mentioned earlier, soft objectives in (\ref{pgST}) do not necessarily require very dense sampling. In this application, we have used a linear spacing of $1\mathrm{e}{3}$ frequencies over  the range $\Omega_S = [1\text{e-}{3},\, 1\mathrm{e}{2}]$ for tracking and $\Omega_T = [1\text{e-}{1},\, 1\mathrm{e}{3}]$ for roll-off.

\begin{figure}[ht!]
\includegraphics[height=0.20\textheight, width = 0.4\textwidth]{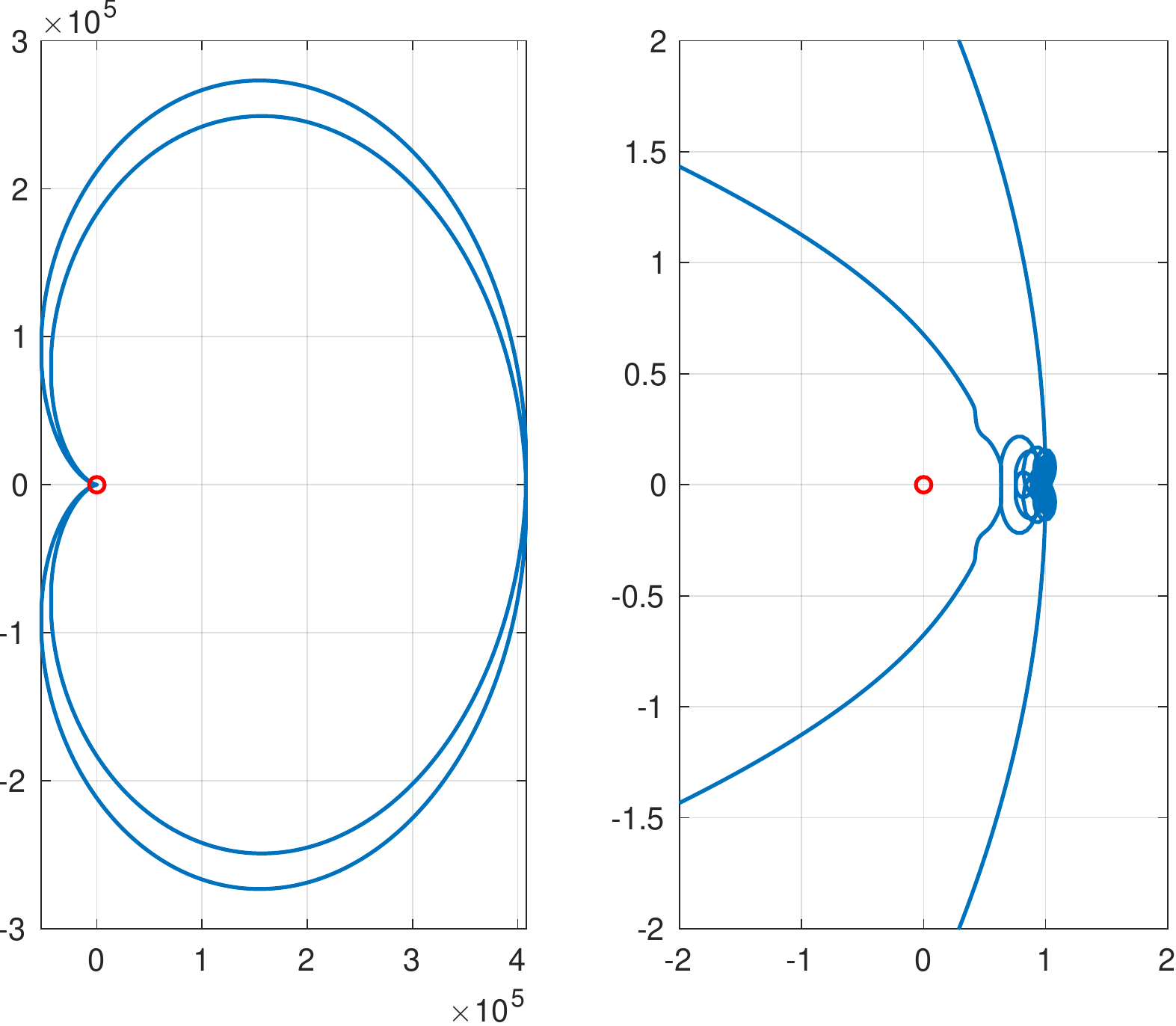}
\includegraphics[height=0.20\textheight, width = 0.4\textwidth]{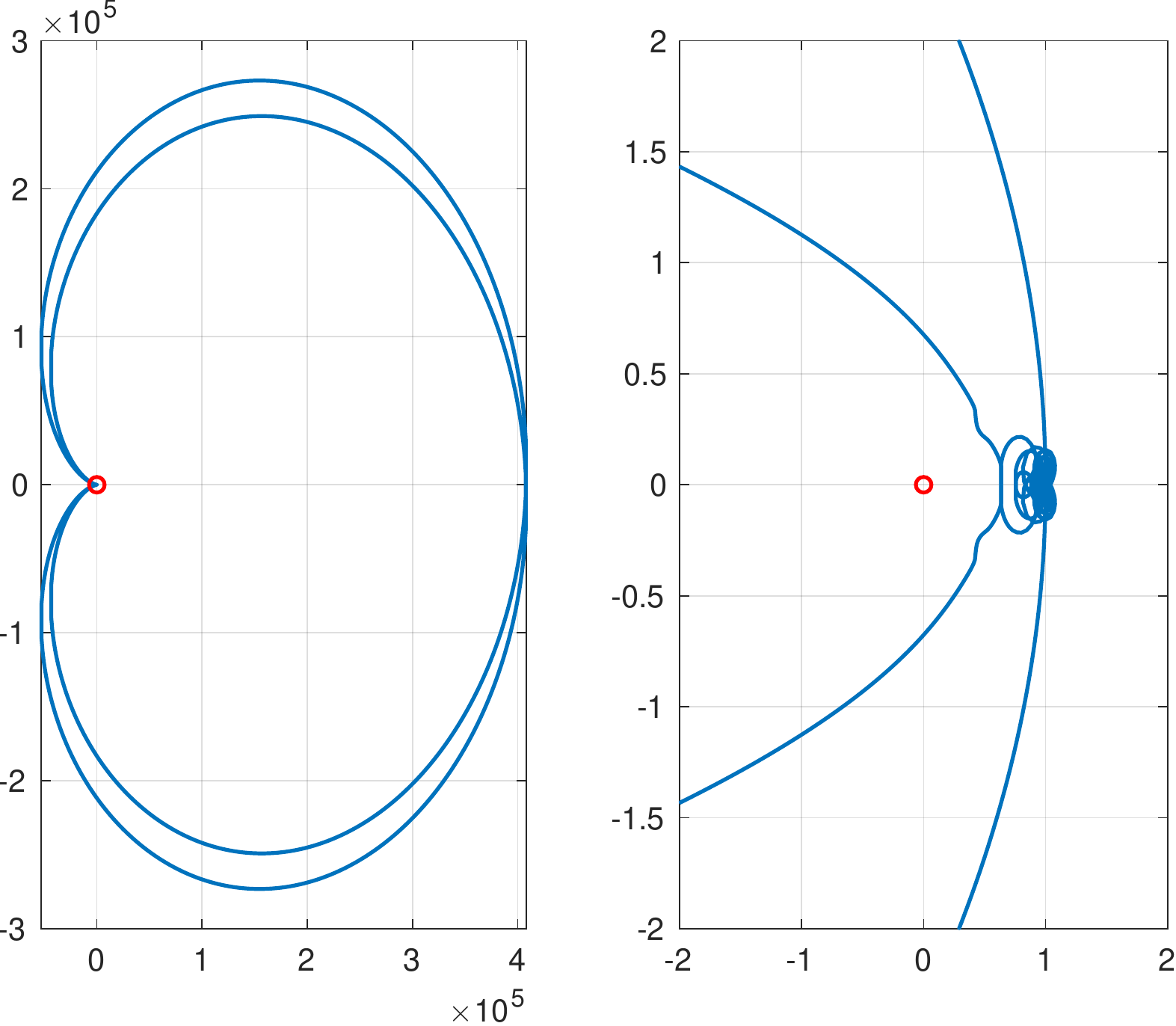}
\caption{Undamped Timoshenko beam. Multivariable Nyquist for reduced-order design. Two left:  reduced-order model. Two right: infinite-dimensional model. Critical point, red dot, at the origin. 
\label{nyquistGred}}
\end{figure}

\begin{figure}[ht!]
\includegraphics[height=0.25\textheight, width = 0.4\textwidth]{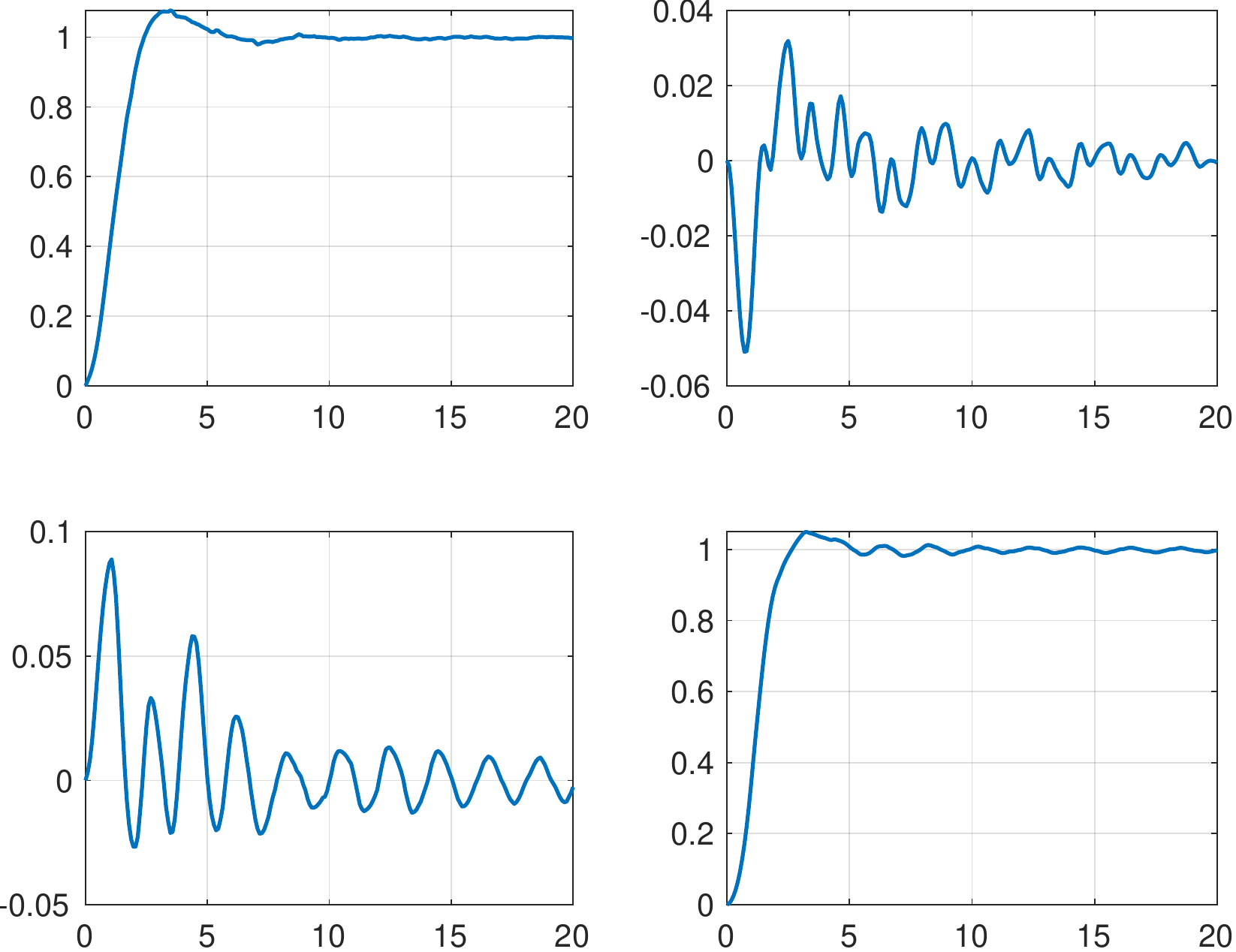}
\includegraphics[height=0.25\textheight, width = 0.4\textwidth]{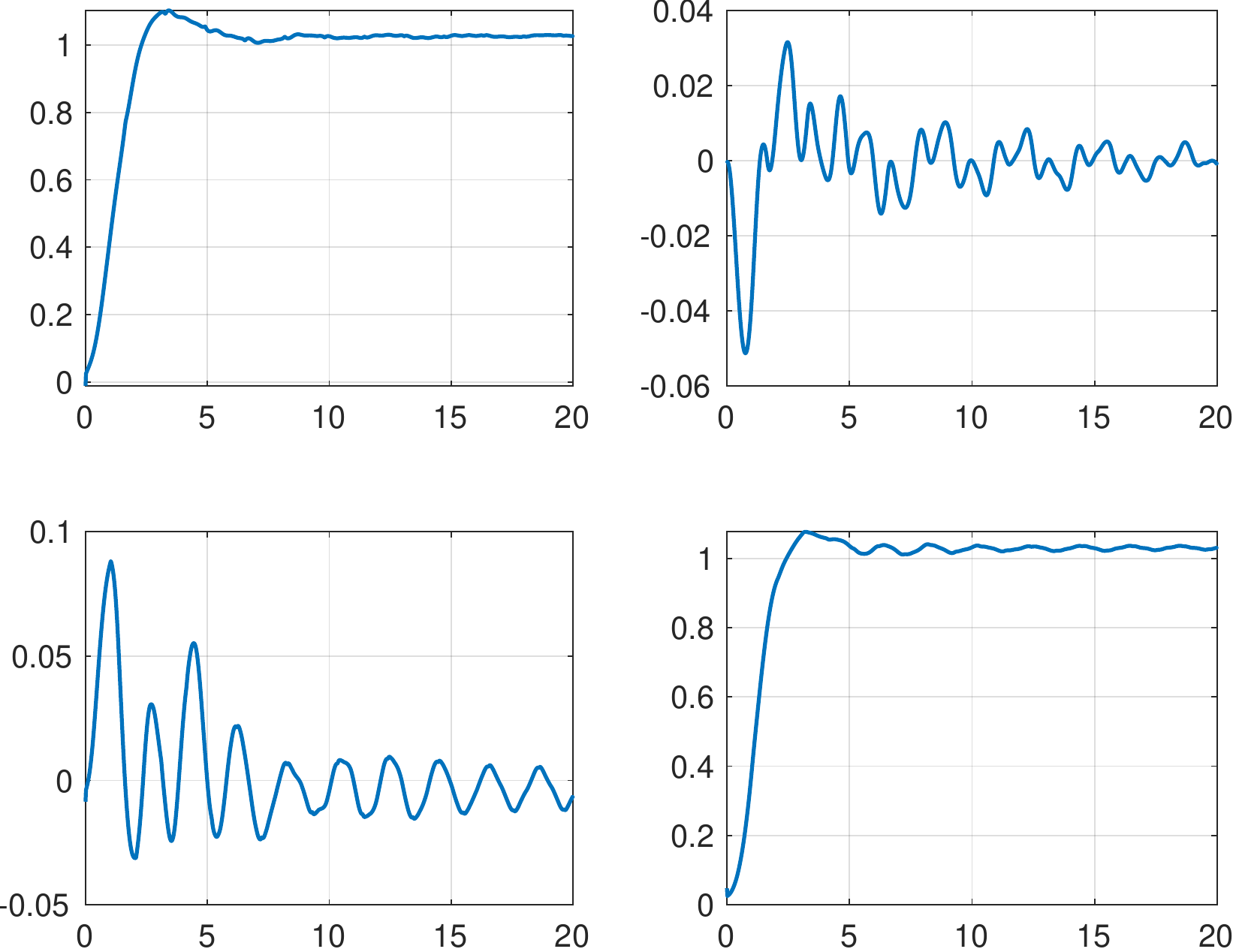}
\caption{Undamped Timoshenko beam. Closed-loop step responses with controller based on reduced-order model.
Two columns left: reduced-order model. Two columns right: infinite-dimensional model \label{simuGred}}
\end{figure}

\begin{figure}[ht!]
\includegraphics[height=0.20\textheight, width = 0.4\textwidth]{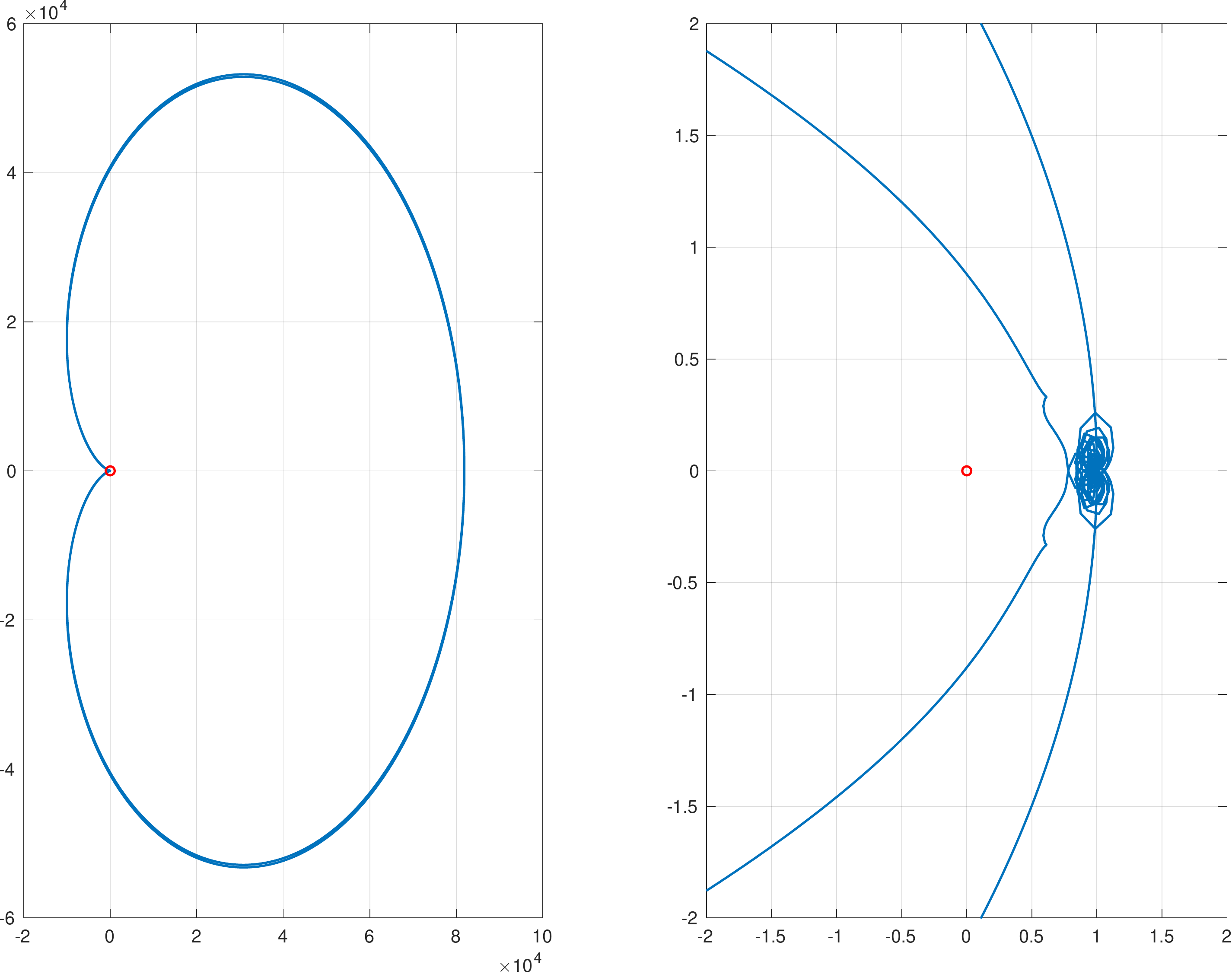}
\includegraphics[height=0.20\textheight, width = 0.4\textwidth]{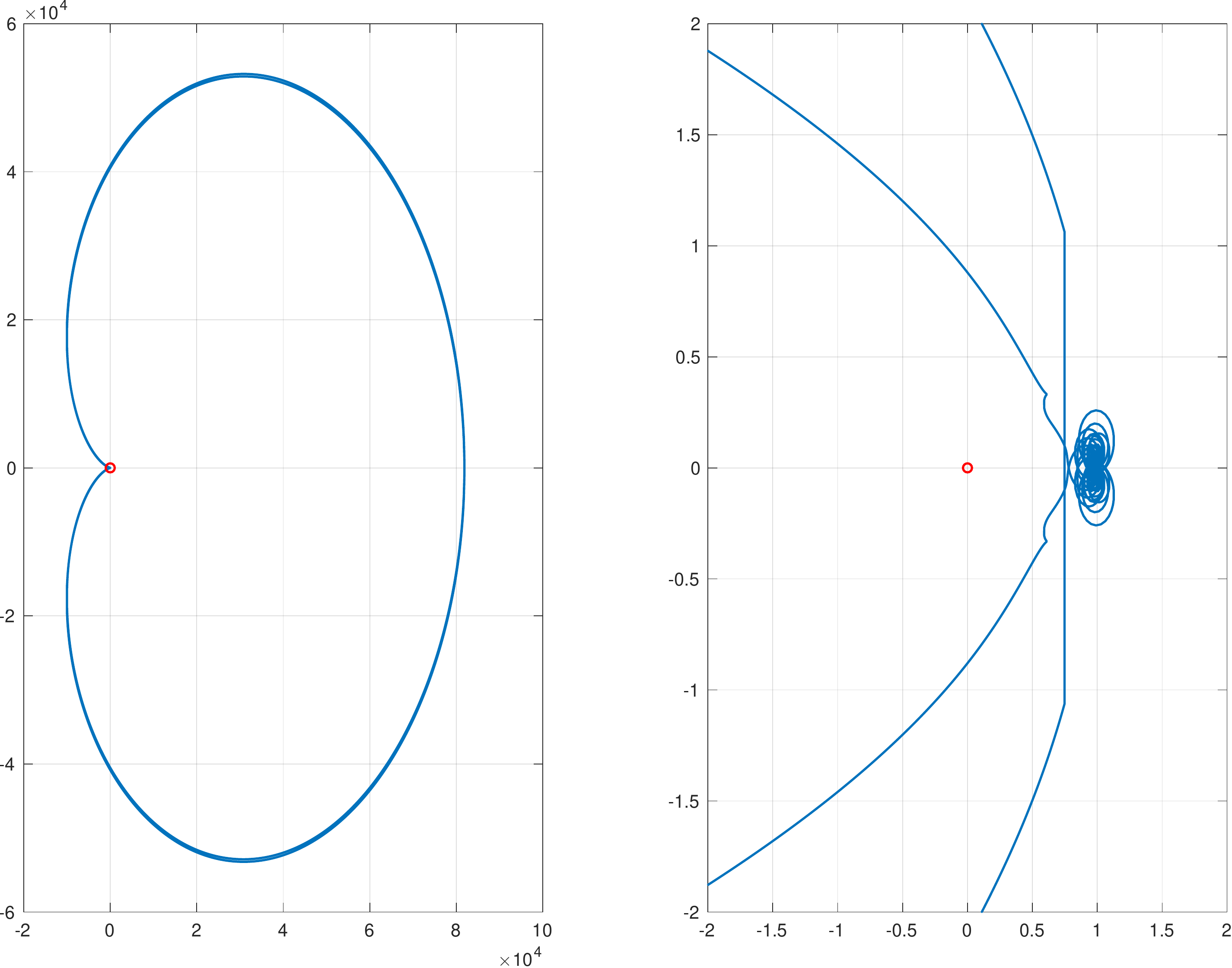}
\caption{Undamped Timoshenko beam. Multivariable Nyquist for infinite-dimensional design. Two left:  infinite-dimensional model. Two right:  reduced-order model. \label{nyquistGinf}  }
\end{figure}

\begin{figure}[ht!]
\includegraphics[height=0.25\textheight, width = 0.4\textwidth]{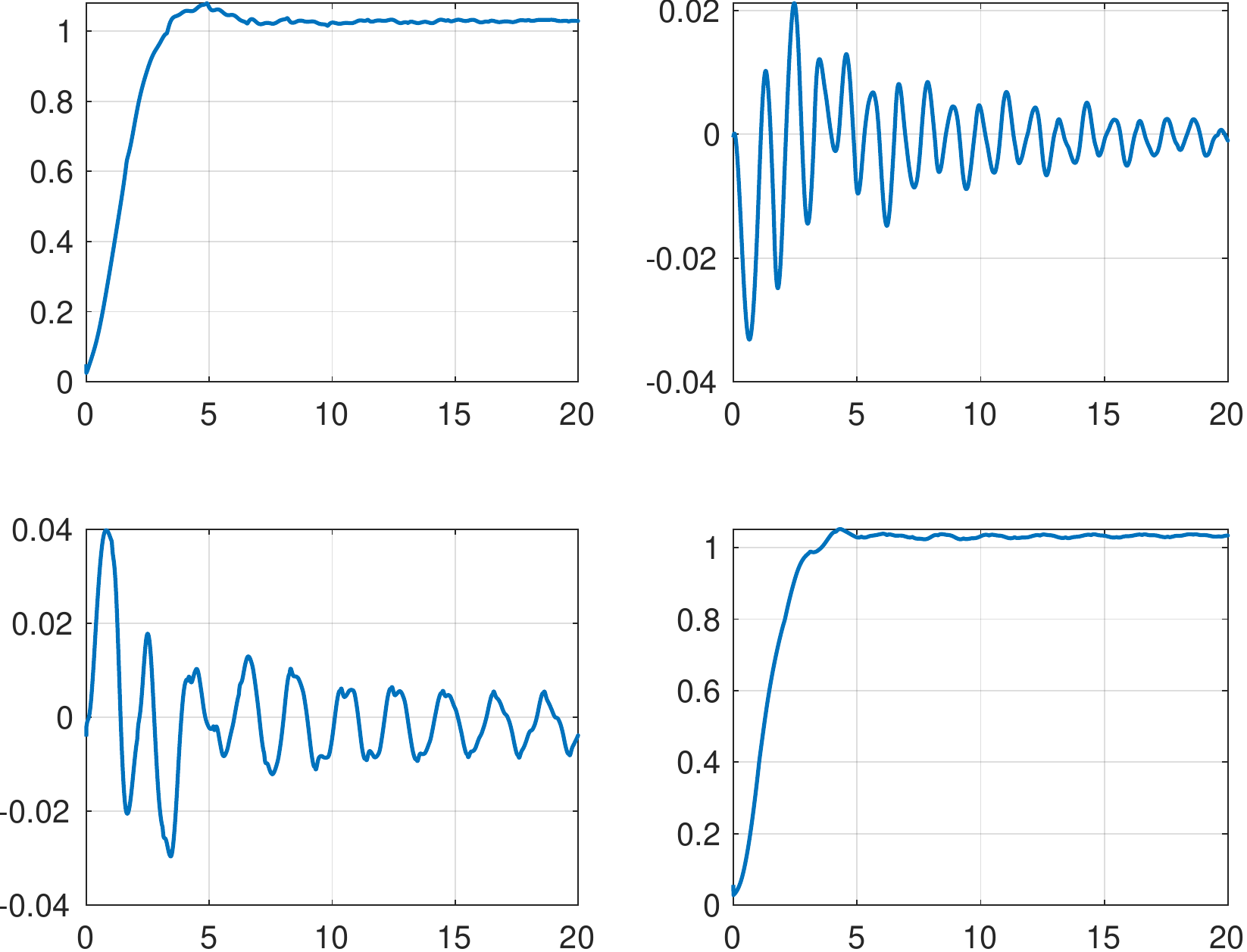}
\includegraphics[height=0.25\textheight, width = 0.4\textwidth]{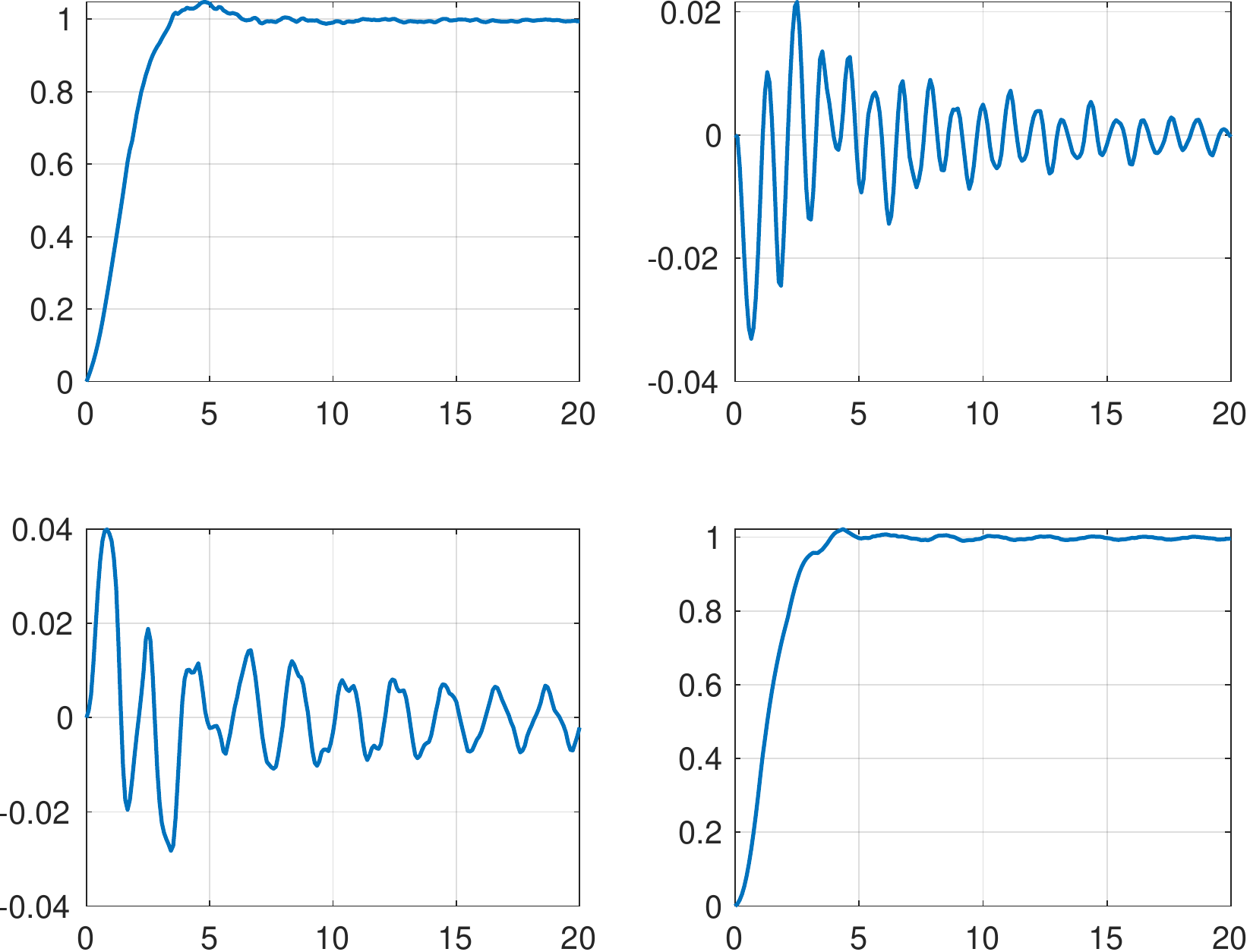}
\caption{Undamped Timoshenko beam. Closed-loop simulations for  infinite-dimensional design. Two columns left: infinite-dimensional model. Two columns right: reduced-order model. \label{simuGinf}}
\end{figure}

 Reduced-order model $G_r=(E,A,B,C,0)$,  controller $K_r$ obtained from $G_r$, full frequency data model $G$ and  controller $K$ obtained through $G$ 
 are compared in Fig. \ref{figG}. For $K_r$ and $K$, high gain is observed at low frequencies, low gain at high frequencies, as required
 for good performance and roll-off, respectively. Multivariable Nyquist plots  associated with the reduced-order model and the infinite-dimensional model are shown  in Figs. \ref{nyquistGred} and \ref{nyquistGinf}. Stability is confirmed in all cases with zero winding numbers as required.  We observe little differences between
 the two approaches, despite the relatively low order of the state-space model, which may be attributed to the roll-off effect introduced by the controller. Step response simulations in a similar configuration are displayed in Figs. \ref{simuGred} and \ref{simuGinf} and indicate good agreement in terms of rise time and decoupling. 
 
 Simulations for the infinite-dimensional model were obtained using numerical Laplace inversion based on the conventional Bromwich contour. Talbot's idea of using  deformed contours \cite{dingfelder2015improved,trefethen2007} is not exploited, but we noticed that shifting the Bromwich contour along the $x$-axis, while preserving analyticity, can help  improving accuracy. 

Finally, it appears that both strategies compete on equal terms in this example.  Stability, performance and robustness are nearly indistinguishable. A major advantage of the direct frequency-domain approach being that the critical phase of devising a suitable reduced-order model along with post-certification in infinite dimension are entirely  bypassed. We will therefore favor this approach in the next applications.

\subsection{Timoshenko beam with damping\label{sect-TimoD}}
\begin{figure}[ht!]
\includegraphics[height=0.3\textheight, width = 0.9\textwidth]{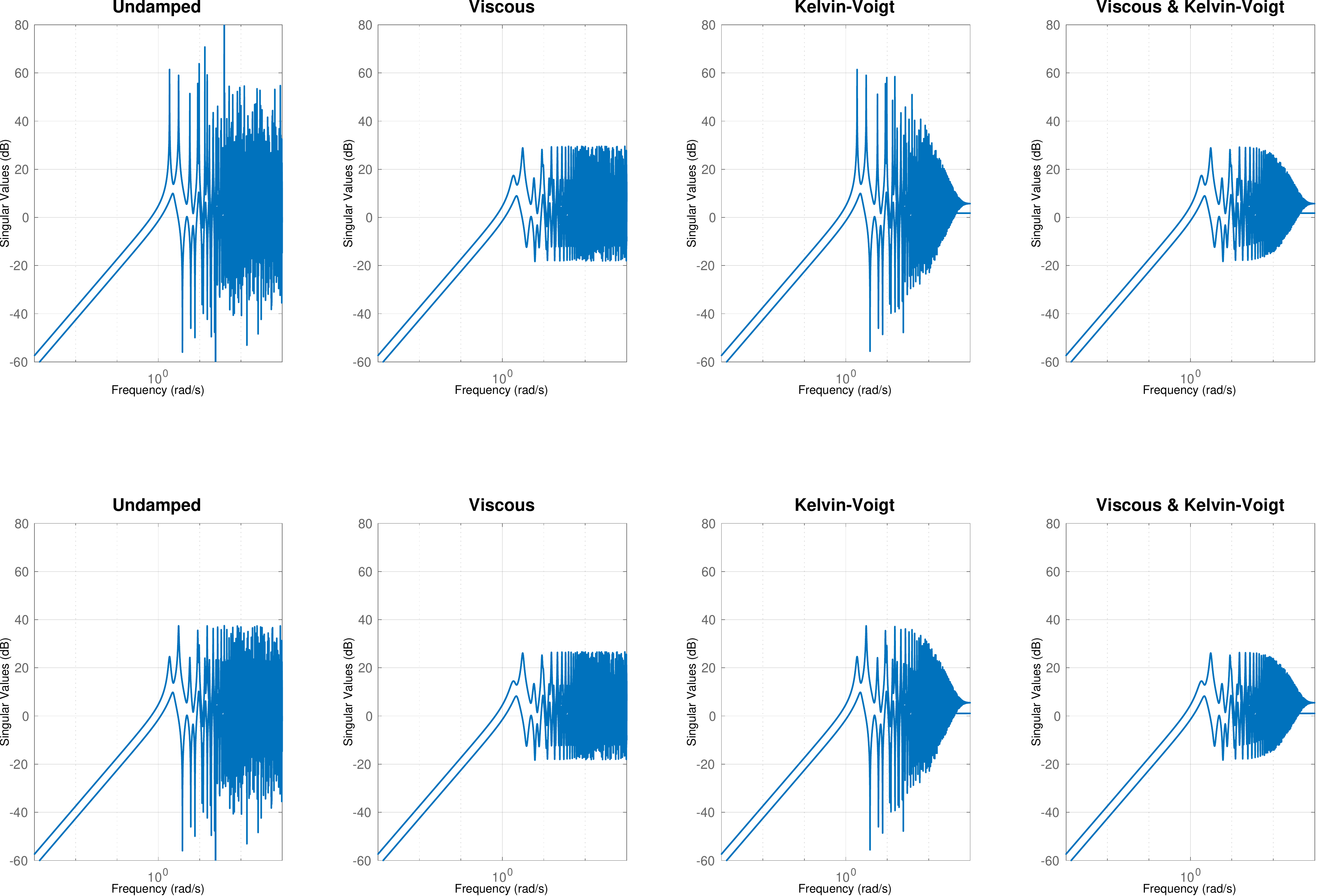}
\caption{Timoshenko beam. First row: without pre-stabilizer. Second row: with pre-stabilizer $K_0$. From left to right: Undamped, viscous  alone, Kelvin-Voigt  alone, viscous and Kelvin-Voigt.
\label{sigmaGinfAll}}
\end{figure}

In this section, we consider the Timoshenko beam with two types of damping, 
viscous, Kelvin-Voigt, and both together. Frequency responses of the beam are shown in Fig. \ref{sigmaGinfAll}. 

It is observed that viscous damping preserves the hyperbolic pole pattern 
of the undamped case.  The difference is that resonances are significantly attenuated,
in our example by a factor of $10$ dB when the pre-stabilizer $K_0$ is on and $30$ dB otherwise.  This is related to a uniform pole shift in the left half plane, see Section \ref{sect-Timo}.  With Kelvin-Voigt damping, the hyperbolic pattern disappears, and poles in the left half plane take a parabolic form. 
Peaks   fade out 
with increasing frequencies and vanish beyond $1\text{e}3$ rad/s. High frequency modes are much more  damped than lower ones. See also Section \ref{sect-Timo}.  With both viscous and Kelvin-Voigt damping a combination
of both effects are observed in the pole pattern.

\subsubsection{Kelvin-Voigt damping with pre-stabilizer}
\begin{figure}[ht!]
\includegraphics[height=0.25\textheight, width = 0.32\textwidth]{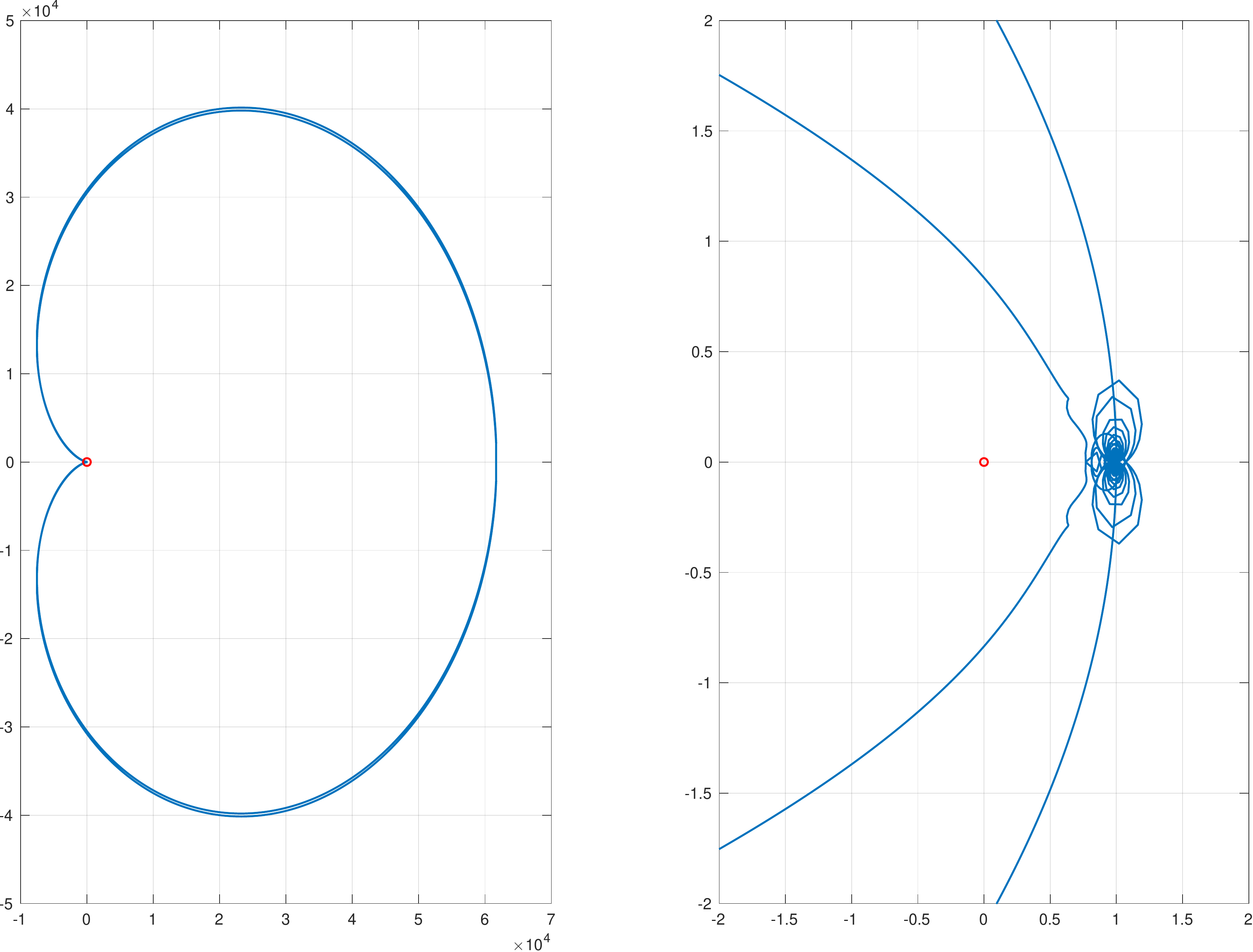}
\includegraphics[height=0.25\textheight, width = 0.32\textwidth]{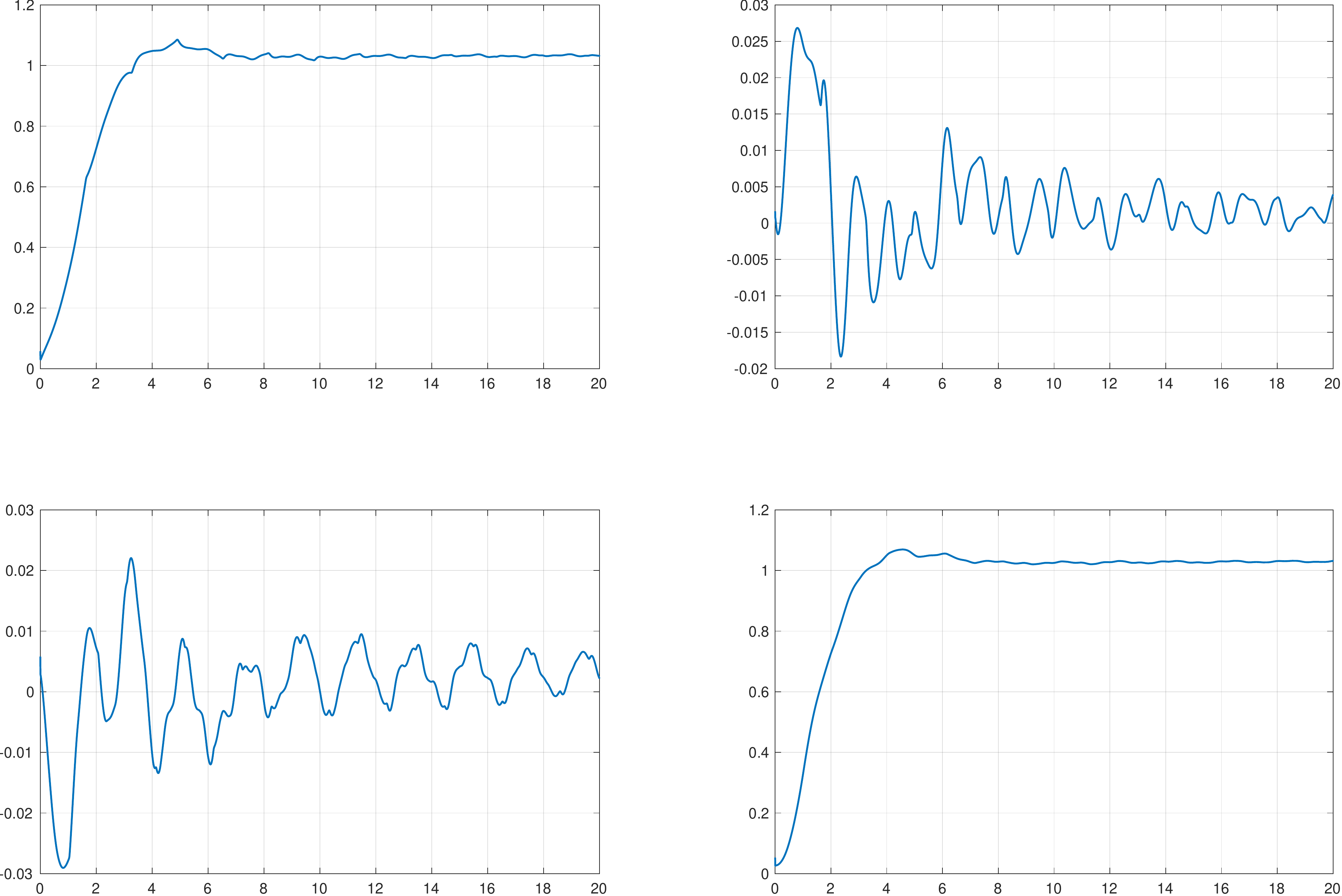}
\includegraphics[height=0.25\textheight, width = 0.32\textwidth]{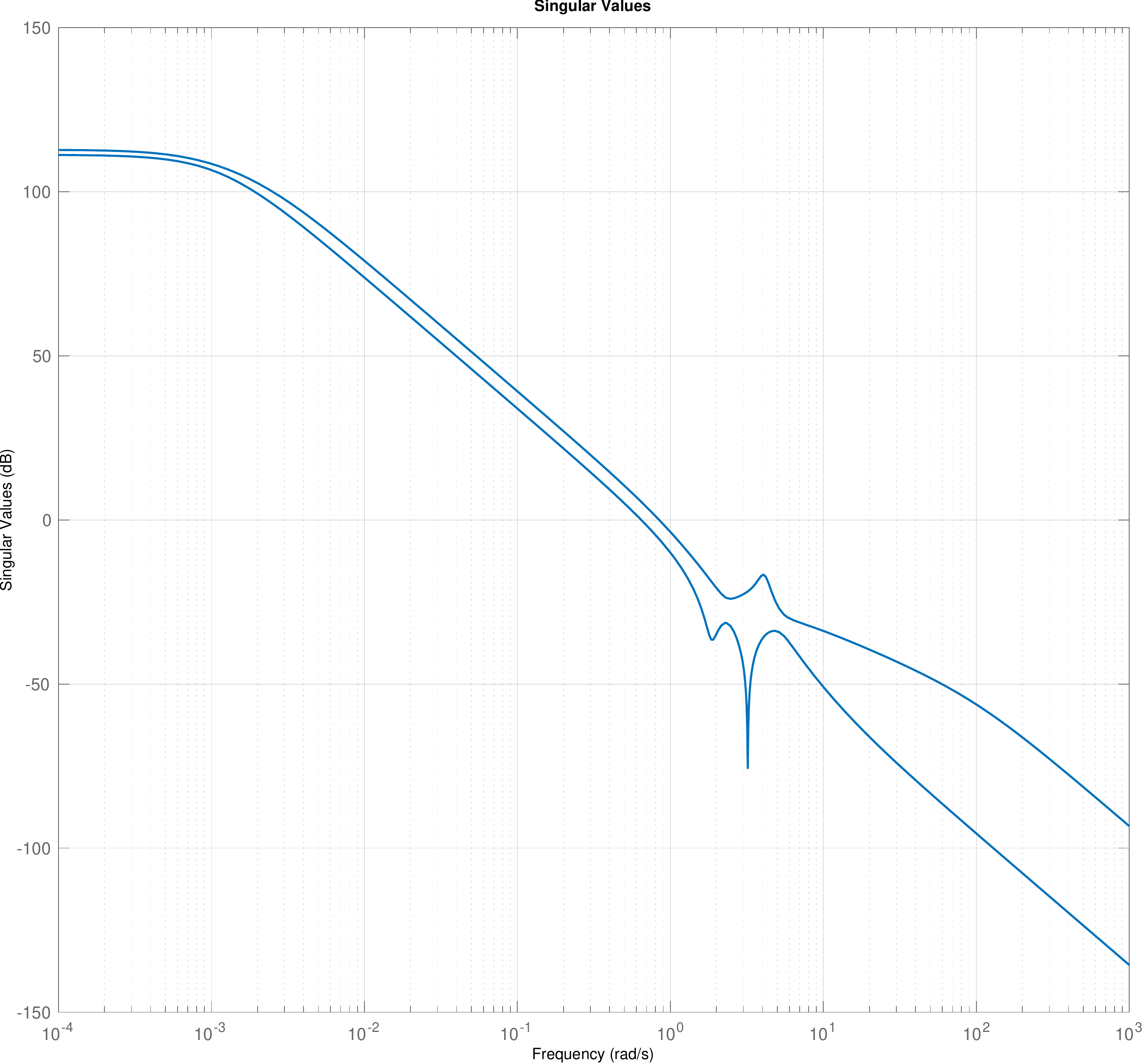}
\caption{Timoshenko beam with Kelvin-Voigt damping alone and pre-stabilizer. Two left: Nyquist plot. Middle: step responses. Right: singular value plot of controller. \label{resultsKV}}
\end{figure}

We consider Kelvin-Voigt damping alone and re-use Program (\ref{pgST}) without change, 
because the low frequency peaks are identical to those of the undamped case. Results are shown in 
Fig. \ref{resultsKV}. Similarly to the undamped case, strong roll-off is achieved for $\omega > 1$ rad/s with notch filtering of low frequency resonances. Good rise times are attained with no more than $4\%$ coupling  of the responses.

\subsubsection{Kelvin-Voigt damping without pre-stabilizer}
\label{without}
Comparing frequency responses in the $3$rd column of Fig. \ref{sigmaGinfAll}, we see that resonances are amplified by $30$ dB when the pre-stabilizer is off. This leads now
to a much more complicated design problem, which requires a review of objectives and constraints,
the challenge being to enable the Nyquist test during optimization. The cut-off frequency of performance weight $W_1$ was reduced from $1$ to $0.5$ rad/s. For the roll-off weight $W_2$, we specified an attenuation of $15$ dB at $1.88$ rad/s corresponding to the $1$st resonant mode. 
This lead to
$$
W_1(s): = \frac{0.05 s + 0.4994}{s + 0.0004994},\qquad W_2(s): = \frac{1000 s + 0.3343}{s + 334.3}\,.
$$

\begin{figure}[ht!]
\includegraphics[height=0.25\textheight, width = 0.32\textwidth]{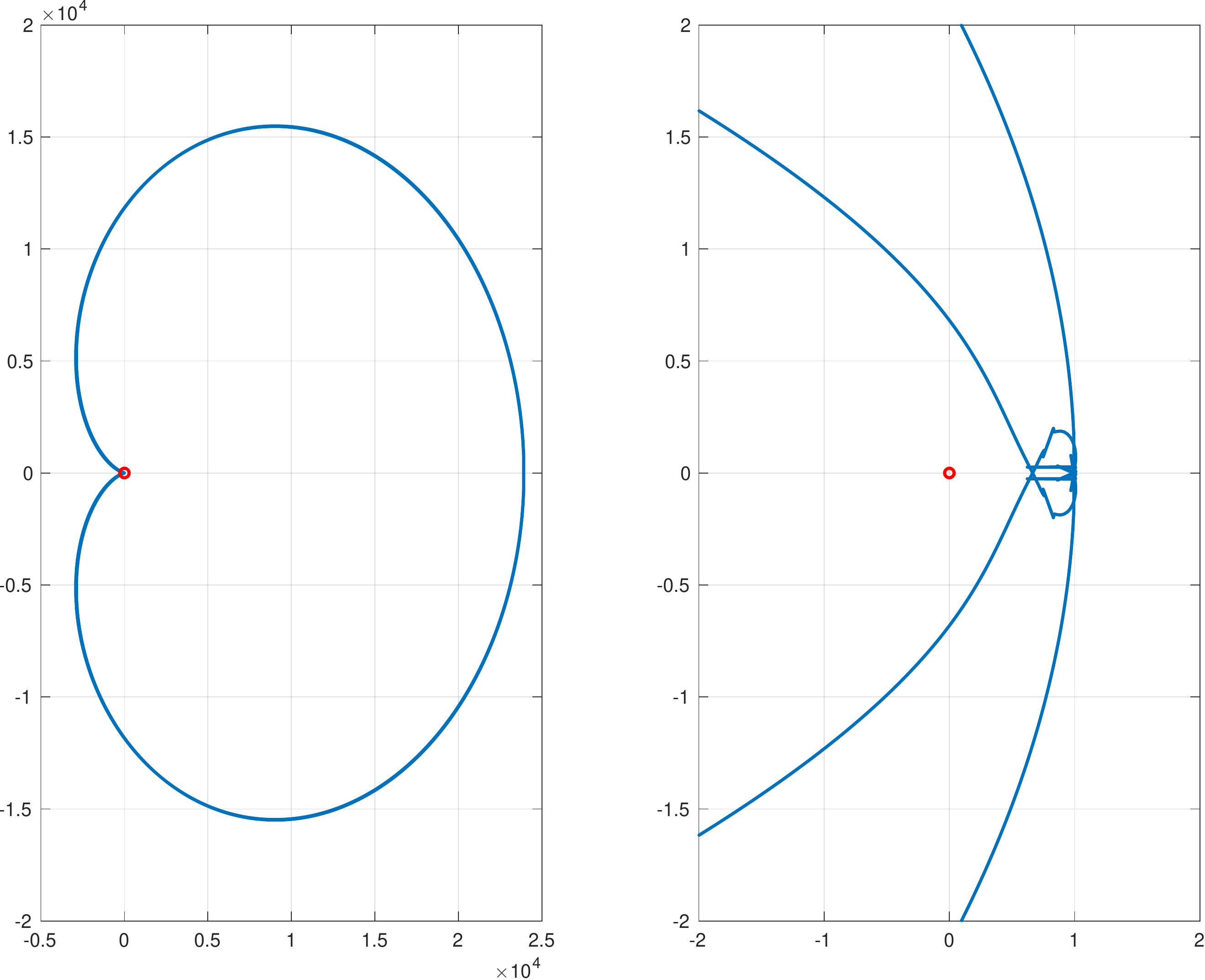}
\includegraphics[height=0.25\textheight, width = 0.32\textwidth]{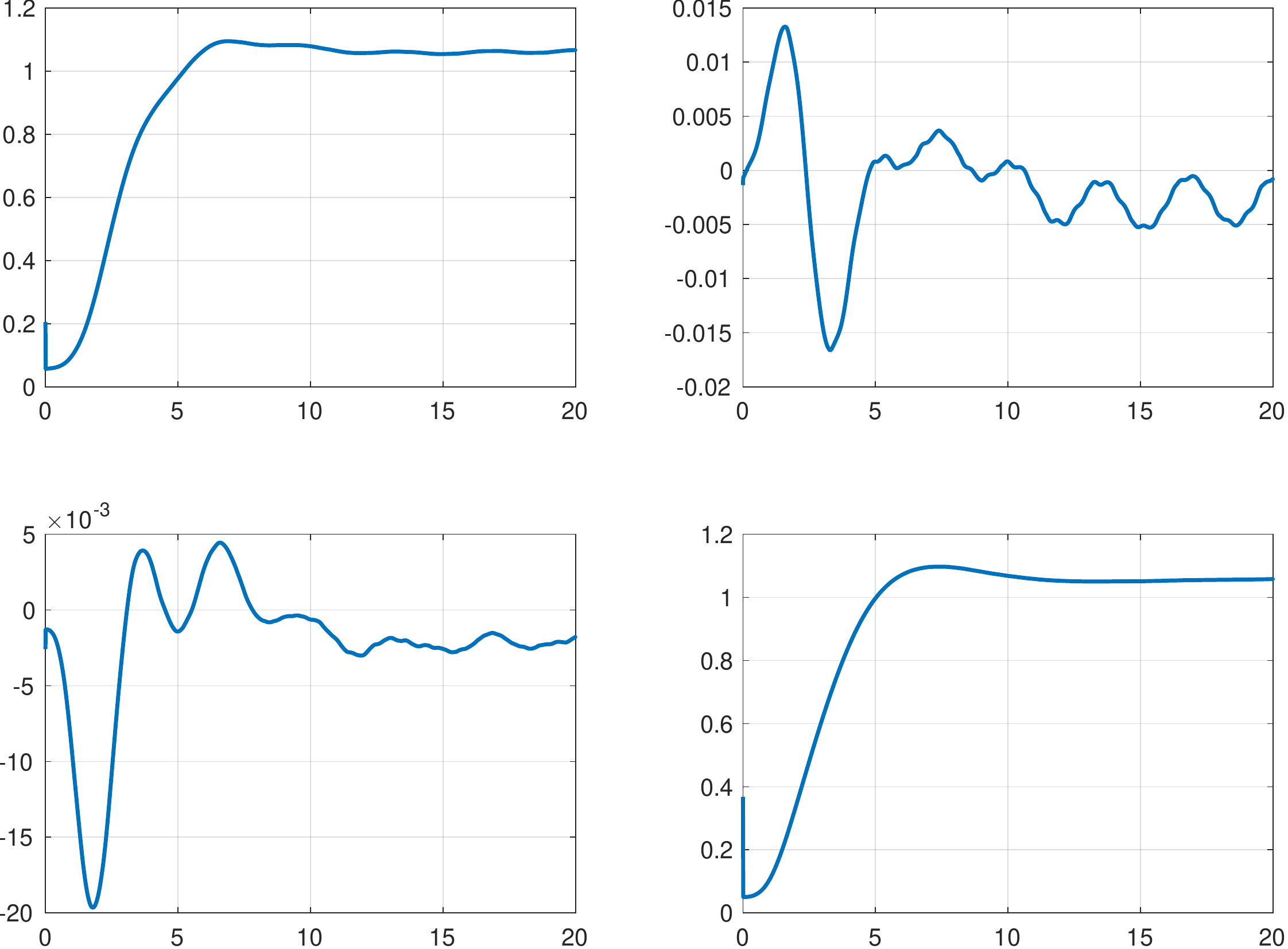}
\includegraphics[height=0.25\textheight, width = 0.32\textwidth]{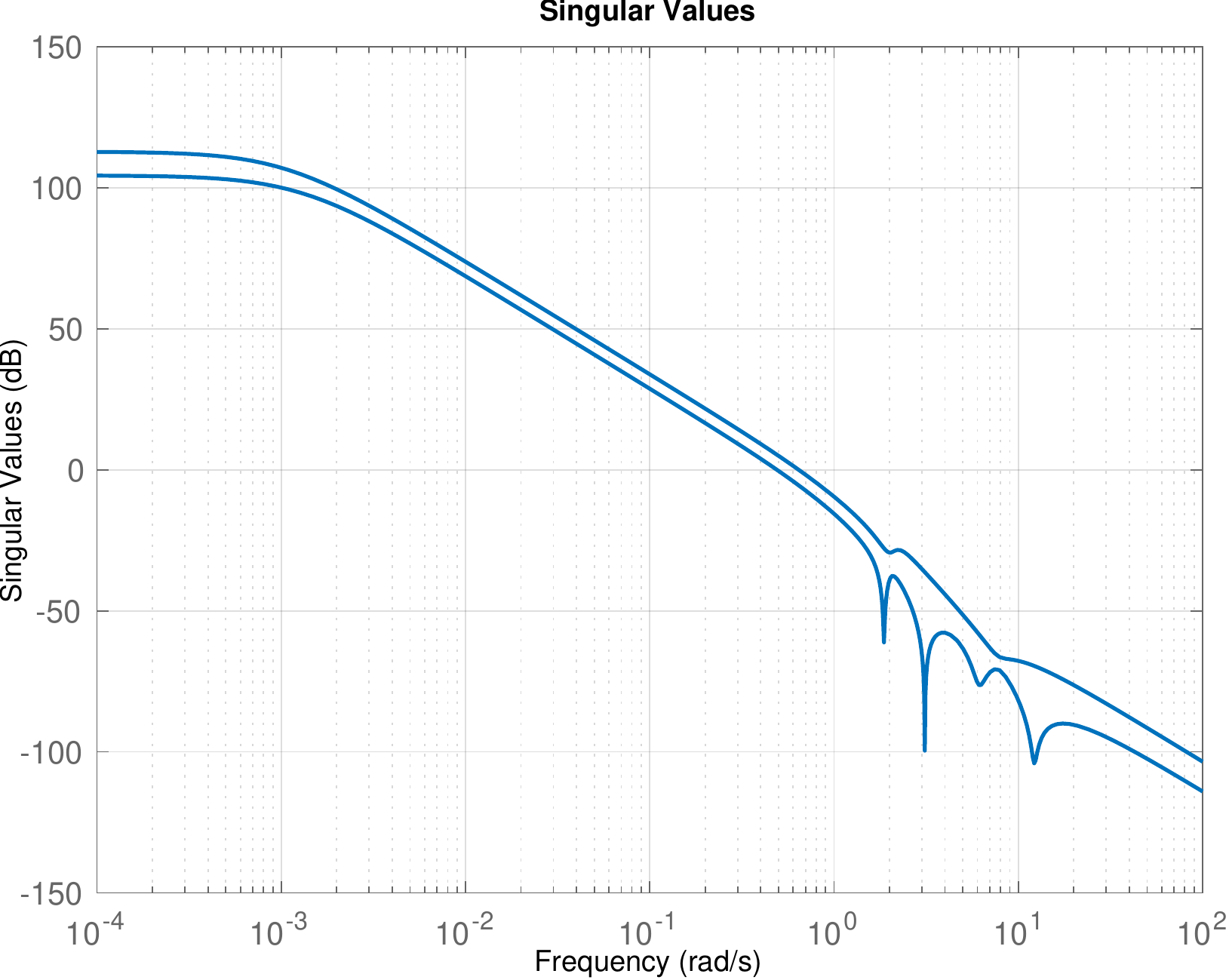}
\caption{Timoshenko beam with Kelvin-Voigt damping and no pre-stabilizer. Two left: Nyquist plot. Middle: step responses. Right: singular value plot of controller. \label{resultsKVwithout}}
\end{figure}

The changes, however, proved insufficient, as bending of poles to
parabolic shape in Fig. \ref{bode_timo_kv1} is fairly weak for the chosen $c_{kv}$. It was necessary to increase the order of the controller from $8$ to $10$ for $\widetilde{K}$, leading to the structure $K := 1/(s+\epsilon) \widetilde{K}$ with $\widetilde{K} \in \mathscr K_{10}$. Also, the disk margin requirement in (\ref{pgST}) had to be alleviated from $0.8$ to $0.7$ 
which is equivalent to setting $\gamma = 1/0.7$ in program (\ref{pgST}). Results are consistent with what we had obtained so far except that lower quality performances are achieved. See Fig. \ref{resultsKVwithout}.

\subsubsection{Viscous damping alone}
\begin{figure}[ht!]
\includegraphics[height=0.25\textheight, width = 0.32\textwidth]{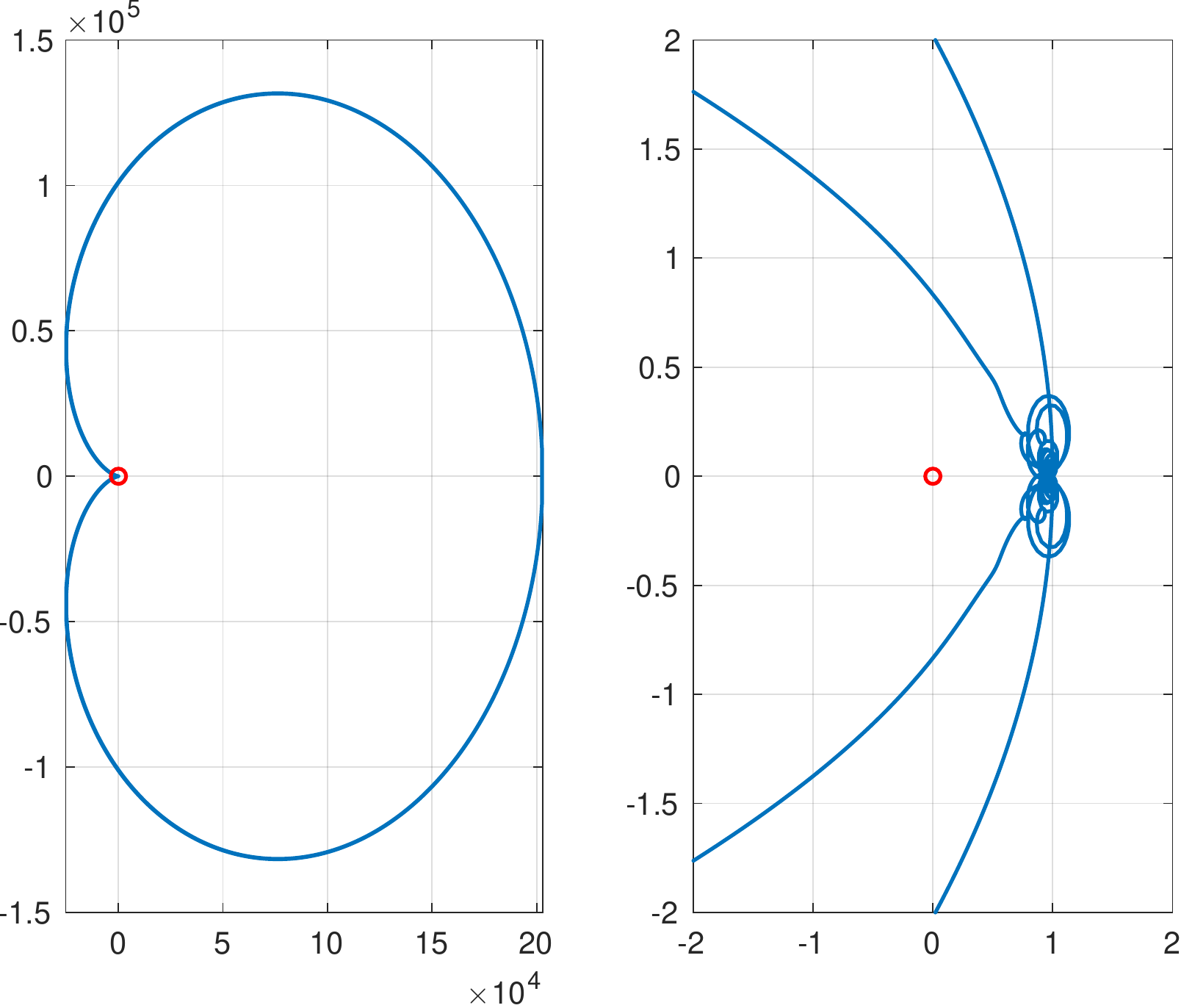}
\includegraphics[height=0.25\textheight, width = 0.32\textwidth]{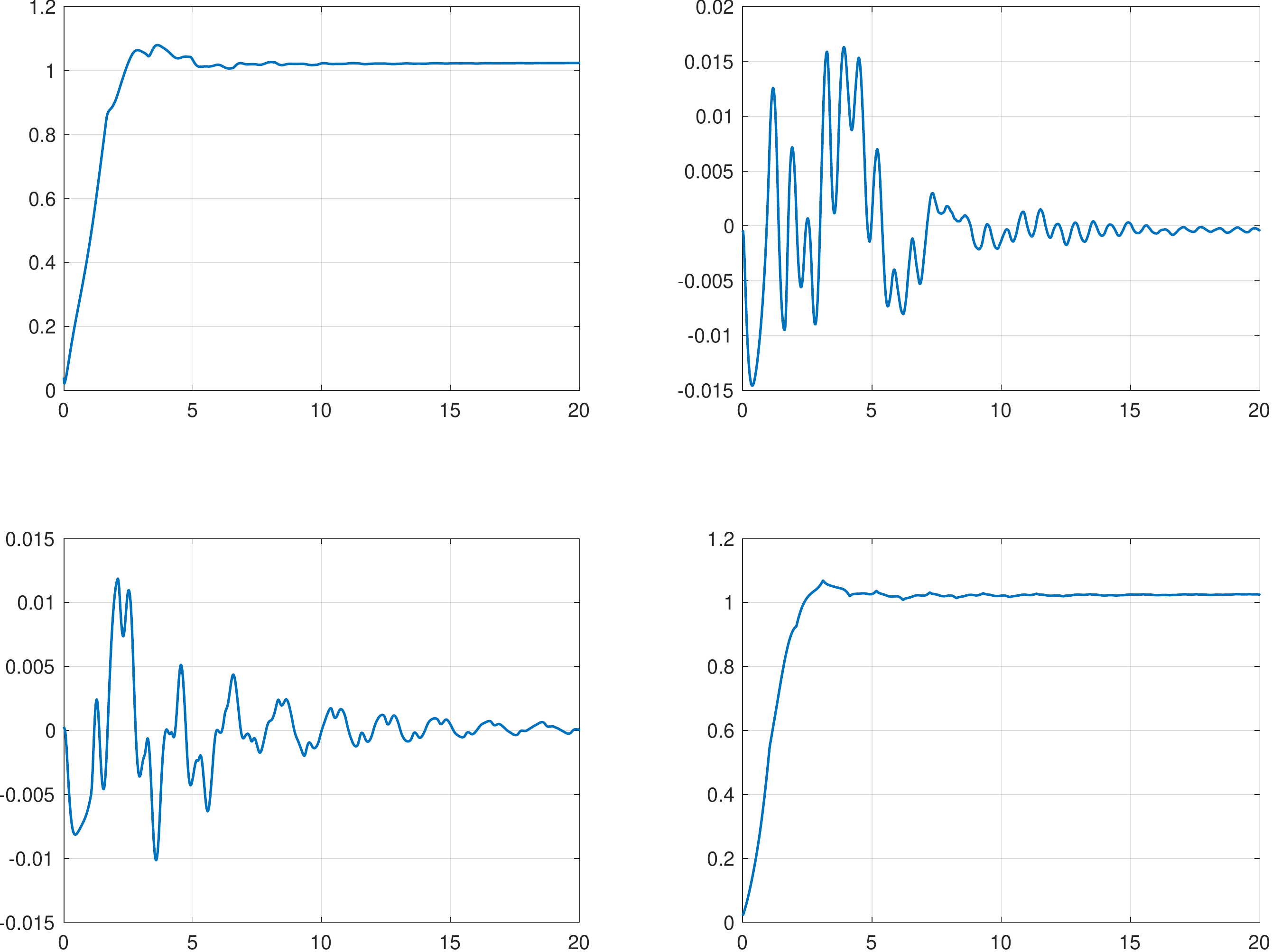}
\includegraphics[height=0.25\textheight, width = 0.32\textwidth]{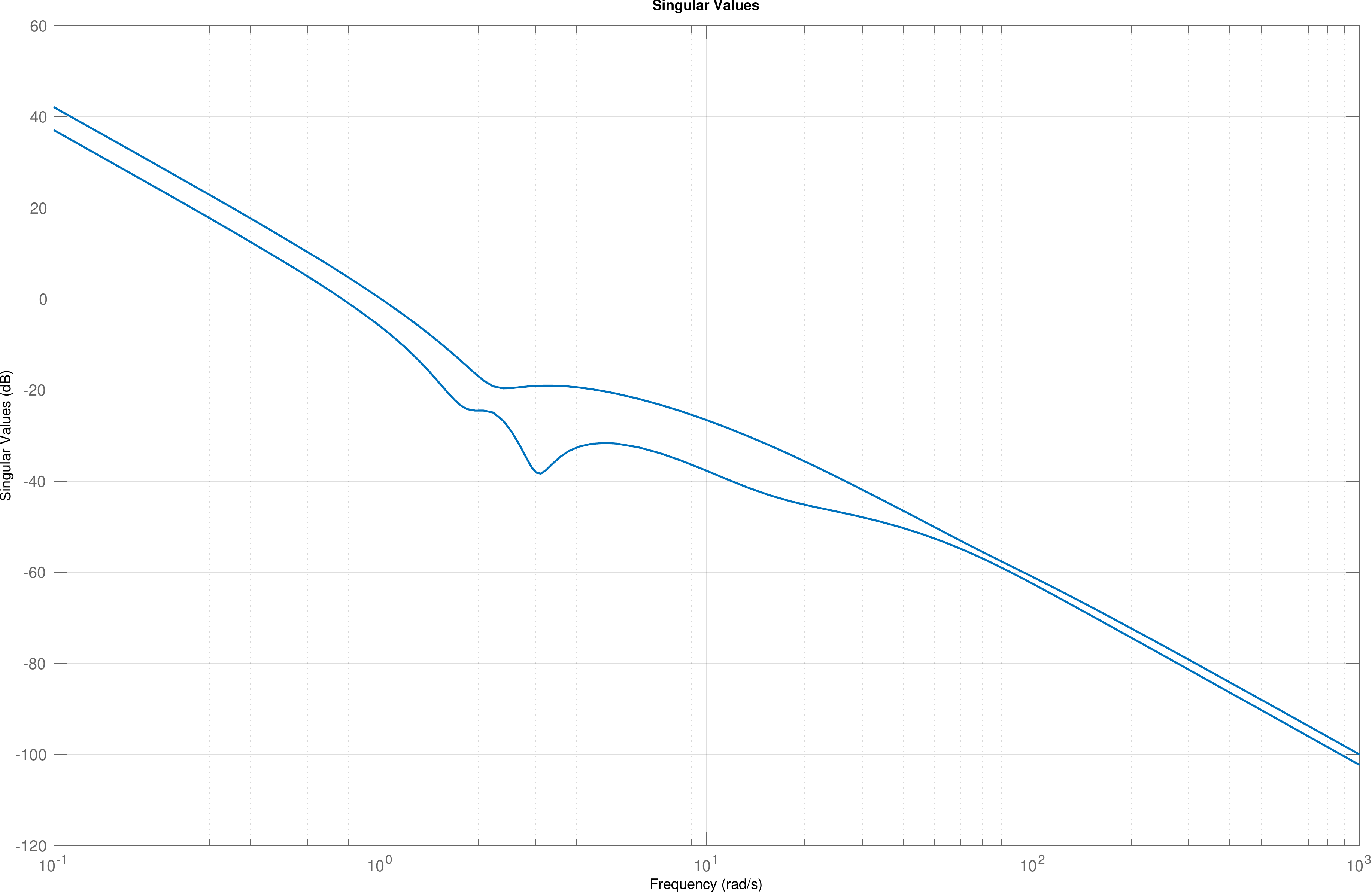}
\caption{Timoshenko beam with viscous damping alone. Two left: Nyquist plot. Middle: step responses. Right: singular value plot of controller. \label{resultsVISCOUS}}
\end{figure}

For viscous damping, frequency responses in Fig. \ref{sigmaGinfAll} suggest that better performance can be achieved. This is easily assessed by applying simple changes to the performance  filter $W_1$. The DC gain of the filter is increased from $1\text{e}3$ to  $3.1623\text{e}3$ corresponding to $10$ dB amplification. Its cut-off frequency at $0$ dB is increased from $1$ to $1.3$ rad/s. Running again program (\ref{pgST}) with the other constraints untouched leads to the results in Fig. \ref{resultsVISCOUS}. The Nyquist criterion certifies exponential stability, based on
Proposition \ref{morris}.
The improvement in the rise time is from about $4$ sec. (undamped and Kelvin-Voigt damping) to $2$ sec. with  viscous damping  $d_w = d_\phi = 0.5$. Better decoupling is also obtained. 

\section{Feedback control of Euler-Bernoulli beam}
\label{sect_euler_control}
Much the same approach is considered for the Euler-Bernoulli beam. We capitalize on  previous results and use the method based on  the infinite-dimensional model. The formulation is now a SISO tracking problem. The cast (\ref{pgST}) is unchanged in terms of ($\gamma,\delta,\mu$), whereas the weighting filters in the objective function require  adjustment. All models were computed by solving complex elliptic boundary value problems (\ref{elliptic_euler}) as in Sections \ref{sect-algo}-\ref{sect_Euler}. 

\subsection{Viscous damping alone} 
\begin{figure}[ht!]
\includegraphics[height=0.25\textheight, width = 0.8\textwidth]{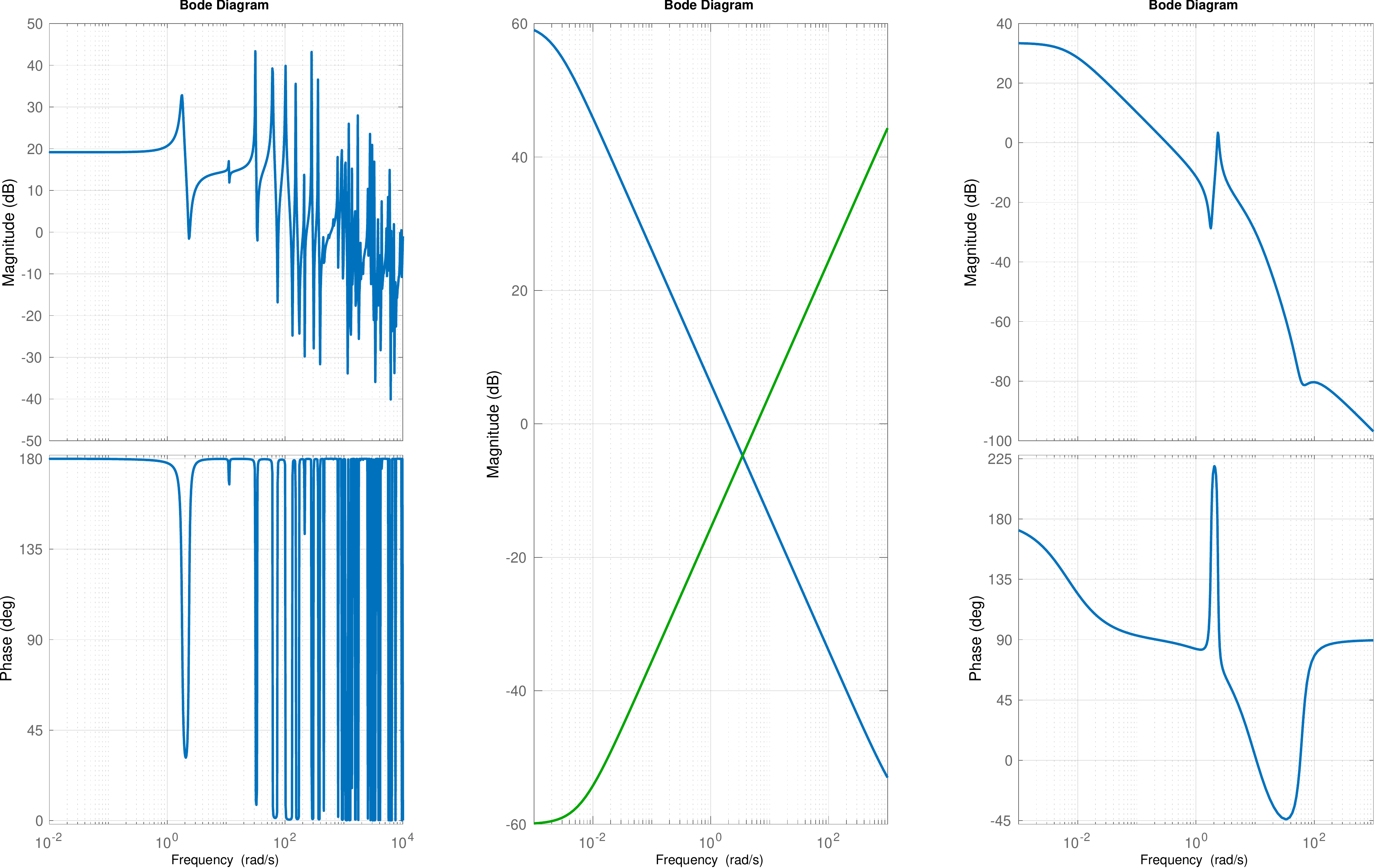}
\caption{Euler-Bernoulli beam. Left: Bode diagram of infinite-dimensional system. Middle: weighting functions used in program (\ref{pgST}). Right: Bode gain and phase diagrams of optimized controller.  \label{allEB}}
\end{figure}

We study the case of a strong viscous damping $c_v = 0.5079$ alone.  The  Euler-Bernoulli model is  to some extent less realistic and exhibits 
resonant modes at higher  frequencies, see Fig. \ref{allEB} left and middle.  Clearly,
this can be exploited to achieve a larger bandwidth and therefore better performance. Weights are shown in Fig. \ref{allEB} middle with transfer functions:
$$
W_1(s) := \frac{0.001 s + 2}{s + 0.002}, \qquad W_1(s) := \frac{1000 s + 6}{s + 6000}\,.
$$
Note that $W_1$ specifies a bandwidth of $2$ rad/s, which corresponds 
approximately to a settling time of $1.5$ seconds. 
Filter $W_2$ reflects a roll-off constraint with crossover at  $6$ rad/s.

\begin{figure}[ht!]
\includegraphics[height=0.25\textheight, width = 0.9\textwidth]{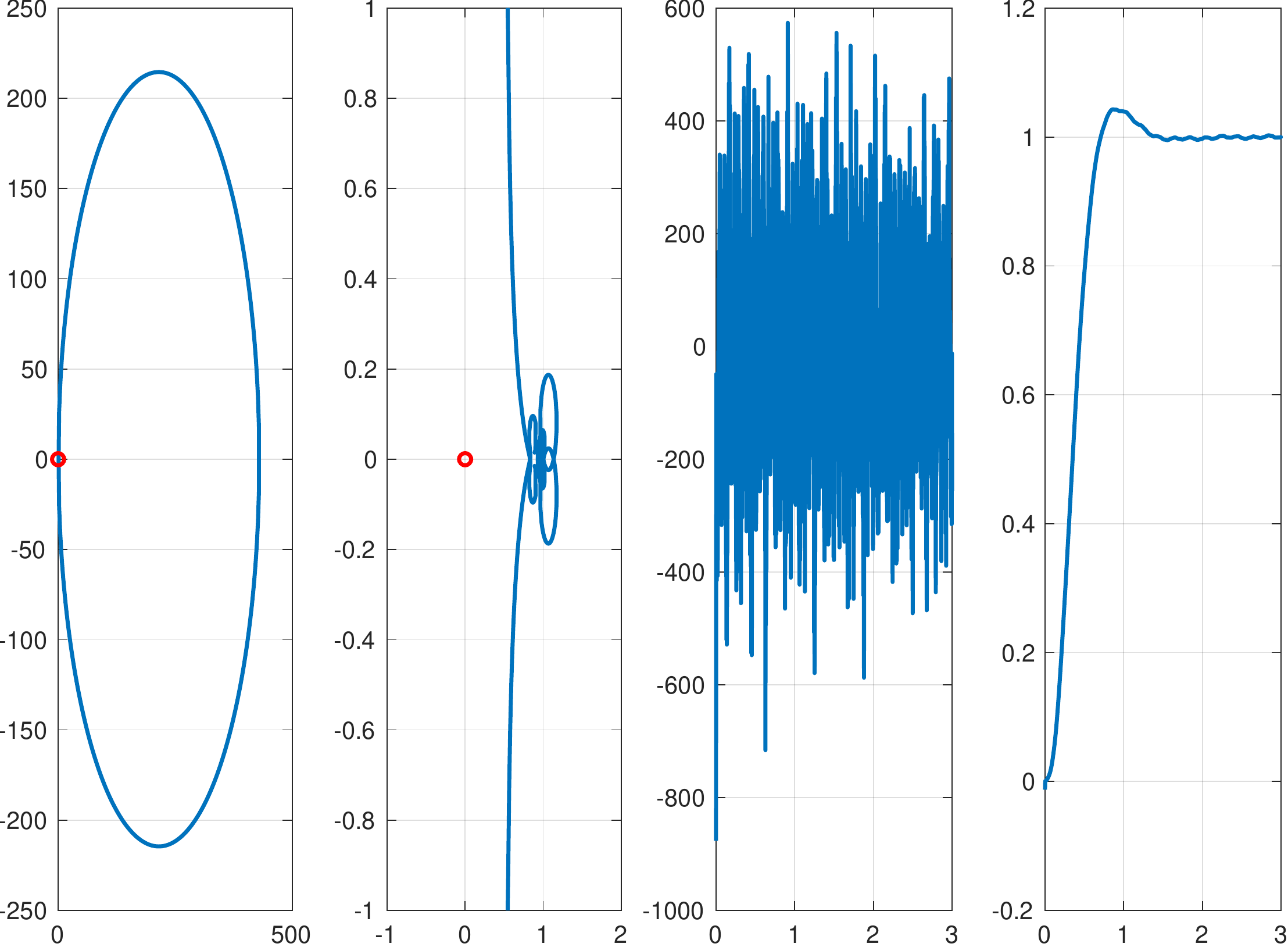}
\caption{Euler-Bernoulli beam. Two leftmost: Nyquist plot. Third: Open-loop step response. Right: Closed-loop step response.   \label{resultsEB}}
\end{figure}

Solving (\ref{pgST}) over $K\in \mathscr K_5$, the set of $5$th-order controllers, gives the controller with Bode diagrams in 
Fig. \ref{allEB} right and transfer function:
$$
K(s):= \frac{ -0.0144 s^4 - 0.3585 s^3 - 54.58 s^2 - 9.669 s - 173.8}{s^5 + 16.9 s^4 + 106.1 s^3 + 108.1 s^2 + 539.7 s + 3.678} \,.
$$
Closed-loop stability is checked using the Nyquist criterion in Fig. \ref{resultsEB}. Open- and closed-loop simulations are compared in the third and fourth plots of Fig. \ref{resultsEB}. Results all agree   with the design constraints. As before, program (\ref{pgST}) was solved with 
$\Omega_D = \Omega'_D\cup \Omega_D''\cup \Omega_D'''$  and $\Omega_N:= \Omega_D$ for the Nyquist criterion.  We have used $250$ samples in the low frequency range $\Omega_D'=[1\text{e-}1,\, 1]$ rad/s, $200$ samples in the mid frequency range $\Omega_D''=[1,\, 6]$ and $500$ samples in $\Omega_D'''=[6,\, 1\text{e}3]$. For $\Omega_S$ and $\Omega_T$, we have used $1\text{e}3$ frequencies over the range $[1\text{e-}1,\, 1e3]$.

\subsection{Kelvin-Voigt damping} 
\label{without2}
\begin{figure}[ht!]
\includegraphics[height=0.25\textheight, width = 0.9\textwidth]{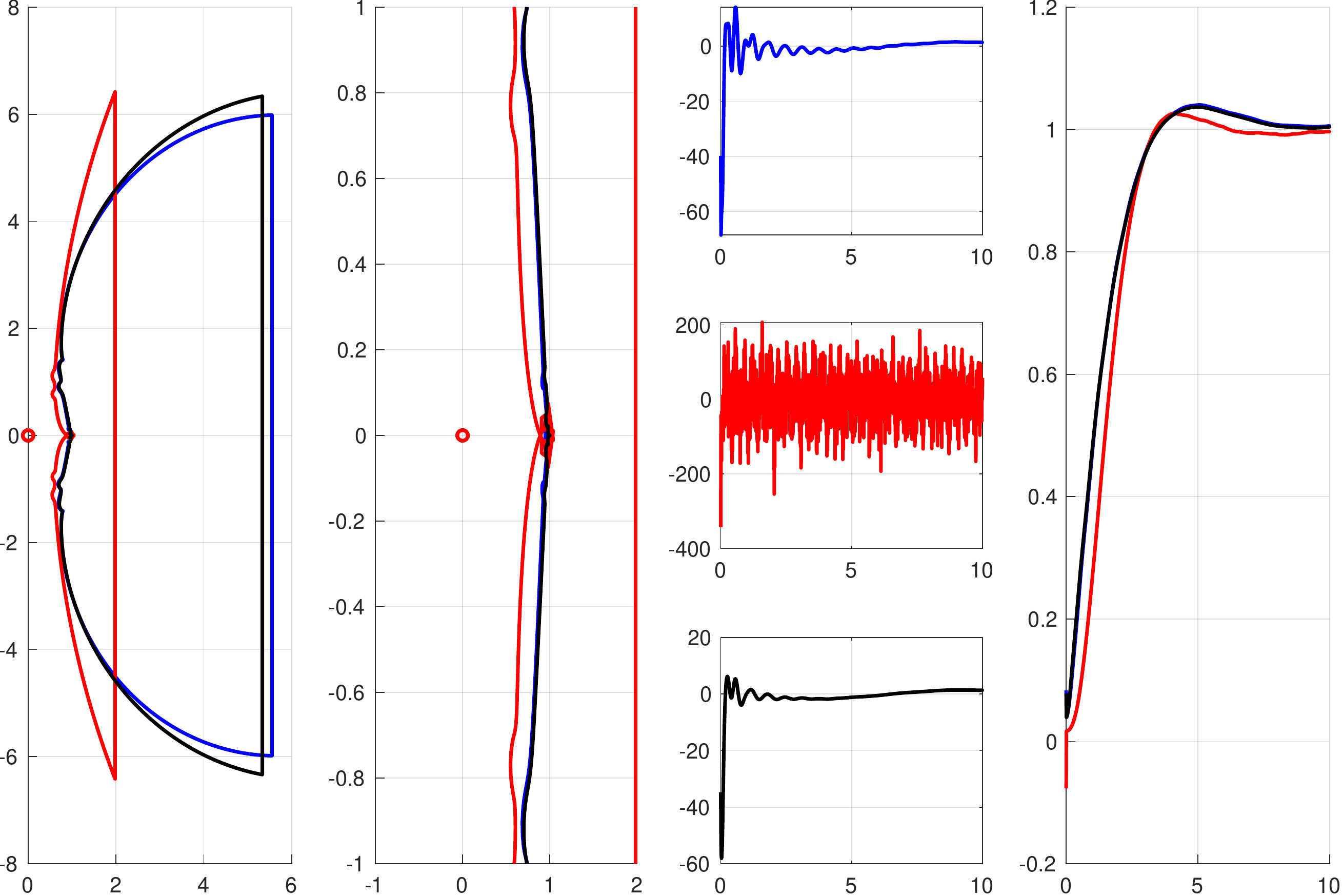}
\caption{Euler-Bernoulli beam. Viscous damping $c_v = 0.156$, Kelvin-Voigt dampings 
$c_{kv}  \in \{ 0, 1\text{e-}3, 2\text{e-}3\}$ with colors red, blue and black. Two leftmost: Nyquist plots. Third column: Open-loop step responses. Right: Closed-loop step responses.   \label{resultsEBKV234}}
\end{figure}

In this last study, we return to the cases discussed in Section \ref{sect_Euler} and displayed in Fig. \ref{TFeuler}. We have moderate viscous damping $c_v = 0.156$ and Kelvin-Voigt dampings 
$c_{kv}  \in \{ 0, 1\text{e-}3, 2\text{e-}3\}$. For the $3$ cases, the first resonant mode appears around  $0.5$ rad/s. Performance requirements should therefore be reduced. The cutoff frequencies at $0$ dB are set to $0.5$ and $1$ rad/s for $W_1$ and $W_2$, respectively. Running program (\ref{pgST})
for all $3$ cases leads to the results shown in Fig. \ref{resultsEBKV234}. As before, stability is assessed through  Nyquist plots. We observe in the $3$rd column that open-loop step responses vary widely with the Kelvin-Voigt damping. The undamped case $c_{kv}=0$, middle plot, exhibits almost persistent oscillations. This is consistent with the pole pattern in the right plot of Fig. \ref{TFeuler}. Step responses in closed loop (right plots) are  however very similar.  This is due to the fact that the pole pattern in Fig. \ref{TFeuler}
match in a low frequency horizontal band $[-10,\,10]$ rad/s.  This confirms, if such confirmation is still needed, that low frequency resonances are the main limiting factor in feedback tracking problems regardless of  pole asymptotics.  Controllers computed in this study are simple state-space systems of order $5$ and are available upon request. 
We note that our approach is a major progress over existing procedures, 
which either use PDE control techniques, where optimizing $K$ is impossible 
or controller structures are impractical, 
or rely on low-order approximations such as finite elements, 
where high order dynamics are typically ignored.

\section{Conclusion}
\label{sect_conclusion}
The main purpose of this note was to investigate whether, or to what extent, $H_\infty$-control methods remain applicable in
open loop
systems with an infinity of open loop poles arranged
in vertical strips close to the imaginary axis. This challenging situation
occurs when models include hyperbolic PDEs, or in neutral systems. It may be considered 
unrealistic or non-physical, because such systems have no natural roll-off and
respond substantially to arbitrary high frequency stimuli.  This extremal behavior
is to
some degree defanged when sensor and actuator models are included, or
when damping is added to make models more realistic. We studied these effects exemplarily for a cantilever Timoshenko and an
Euler-Bernoulli beam, where undamped or viscous damped systems exhibit such non-sectorial 
pole pattern, while Kelvin-Voigt damping bends poles to a more realistic sectorial 
shape. Surprisingly, our methods prove effective even in the challenging neutral case. We synthesize finite-dimensional implementable $H_\infty$-controllers, which attenuate high frequency stimuli in closed loop, introduce 
realistic roll-off, and still act sufficiently fast, as required in a technically feasible
controlled system. Synthesis is based on a recent infinite-dimensional frequency-based
optimization technique.

On closer look
the question whether, or to what extent, the infinite dimension of $G$ hampers the choice of synthesis strategies, calling for elements one would not be inclined to consider in a reduced-order
model, has the following partial answer.

We had to enforce significant disk margins  in closed loop,  constraints  on the dynamics of $K \in \mathscr K$ and suitable roll-off in the high-frequency range in order to enable the Nyquist test during optimization and to make our frequency-domain approach more reliable. This restricts the set of reachable controllers to some extent.  
However, in the present beam studies this type of constraint
did not seem totally unnatural and might be opportune even
in reduced-order models. Specifically, disk margins improve robustness of the loop and roll-off constraints attenuate measurement noise and high-frequency resonances. Constraints on $K$, on the other hand, are less conventional and are tied to our frequency sampling technique. 

A more severe restriction occurred in the Timoshenko study when the pre-stabilizer
$K_0$ was not available, as then the lack of roll-off in the system
strongly limited the choice of the synthesis strategy 
(see Section \ref{without}, and to a lesser degree, Section \ref{without2}). Altogether,
this quest may require further investigation with other non-sectorial open-loop
models. A difficulty in the assessment is  that in the literature comparison with
finite element or reduced-order synthesis,
where these restrictions might not be on the agenda, 
often
disregards high frequency effects of the final controllers, by conducting
simulations within the reduced-order beam models only.

A detail to be mentioned is that in beam models velocity measurements seem to alleviate
functional analytic Lyapunov-based  proofs for stabilization, but
render synthesis
harder from a practical control point of view, as e.g. integral action and low frequency high gain for good tracking are impeded. 

Finally, the frequency-domain approach taken in this work offers a myriad of possible extensions including the design of two-degree-of-freedom or multi-block controllers, as well as of controllers dealing with robustness against parameter uncertainties or  self-adjusting to variations in PDE dynamics.


\end{document}